\documentclass[12pt,twoside]{amsart}
\usepackage{amssymb}
\usepackage{amscd}
\usepackage[abbrev,alphabetic]{amsrefs}
\usepackage{hyperref}
\usepackage{comment}
\usepackage{tikz}
\usetikzlibrary{cd}
\usepackage{here}
\usepackage[margin=1.25in]{geometry}

\usepackage{fancybox}
\usepackage{ascmac}

\title[Fano threefolds in positive characteristic I]
{Fano threefolds in positive characteristic I} 
\author{Hiromu Tanaka} 
\address{Graduate School of Mathematical Sciences, 
The University of Tokyo, 
3-8-1 Komaba, Meguro-ku, Tokyo 153-8914, JAPAN} 
\email{tanaka@ms.u-tokyo.ac.jp}
\newcommand{\Diff}[0]{{\operatorname{Diff}}}
\newcommand{\Cl}[0]{{\operatorname{Cl}}}
\newcommand{\Sing}[0]{{\operatorname{Sing}}}

\newcommand{\red}[0]{{\operatorname{red}}}

\renewcommand{\Im}[0]{{\operatorname{Im}}}

\newcommand{\Proj}[0]{{\operatorname{Proj}}}
\newcommand{\Spec}[0]{{\operatorname{Spec}}}

\newcommand{\Bs}[0]{{\operatorname{Bs}}}
\newcommand{\Supp}[0]{{\operatorname{Supp}}}
\newcommand{\Pic}[0]{{\operatorname{Pic}}}
\newcommand{\mult}[0]{{\operatorname{mult}}}
\newcommand{\Ex}[0]{{\operatorname{Ex}}}

\newcommand{\univ}[0]{{\operatorname{univ}}}
\newcommand{\gen}[0]{{\operatorname{gen}}}

\renewcommand{\min}[0]{{\operatorname{min}}}
\newtheorem{thm}{Theorem}[section]
\newtheorem{lem}[thm]{Lemma}
\newtheorem{cor}[thm]{Corollary}
\newtheorem{prop}[thm]{Proposition}
\newtheorem{claim}[thm]{Claim}    
\newtheorem*{claim*}{Claim}         
\newtheorem{step}{Step}

\theoremstyle{definition}
\newtheorem{ex}[thm]{Example}

\newtheorem{dfn}[thm]{Definition}

\newtheorem{rem}[thm]{Remark}
      
\newtheorem{nota}[thm]{Notation}         
\newtheorem{nasi}[thm]{}

\makeatletter
  
  \@addtoreset{equation}{thm}
  \makeatother

\newcommand{\cred}{\color{black}}

\newcommand{\MO}{\mathcal{O}}

\newcommand{\Q}{\mathbb{Q}}
\newcommand{\Z}{\mathbb{Z}}
\newcommand{\F}{\mathbb{F}}

\renewcommand{\P}{\mathbb{P}}
\newcommand{\p}{\mathfrak{p}}

\newcommand{\m}{\mathfrak{m}}

\begin{document}

\maketitle


\begin{abstract}
Over an algebraically closed field of positive characteristic, we classify smooth Fano threefolds of Picard number one whose anti-canonical linear systems are not very ample. Furthermore, we also prove that an anti-canonically embedded Fano threefold of genus at least five is an intersection of quadrics. 
\end{abstract}

\tableofcontents

\section{Introduction}

One of important classes of algebraic varieties 
{\cred is the class of} Fano varieties, 
that is, smooth projective varieties whose anti-canonical divisors are ample. 
The classification of Fano varieties has been an interesting topic in algebraic geometry. 
The projective line $\P^1$ is the unique one-dimensional Fano variety. 
Two-dimensional Fano varieties are called 
del Pezzo surfaces. 
It is classically known that a del Pezzo surface is isomorphic to 
either $\P^1 \times \P^1$ or the blowup of $\P^2$ along at most general eight points. 
Study of three-dimensional Fano varieties, called Fano threefolds, 
was initiated by Fano himself. 
The classification of Fano threefolds was completed by Iskovskih, Shokurov, and Mori--Mukai 
\cite{Isk77}, \cite{Isk78}, \cite{Sho79a}, \cite{Sho79b}, \cite{MM81}, \cite{MM83}, \cite{MM03} 
(cf. \cite{IP99}, \cite{Tak89}).

It is natural to consider the classification of Fano threefolds also in positive characteristic. 
In this direction, Megyesi, {\cred Shepherd-Barron, and Saito} studied the case when $r_X=2$ (of index two), 
$\rho(X)=1$, and $\rho(X)=2$, respectively \cite{Meg98}, \cite{SB97}, 
 \cite{Sai03}. 
The purpose of this series of papers is to complete 
the classification of Fano threefolds in positive characteristic. 
This article focus on the case when $\rho(X)=1$ and $|-K_X|$ is not very ample.

\begin{thm}[{\cred Theorem \ref{t-non-bpf-main}, Section \ref{s-bpf-case}}]\label{t-main1}
Let $k$ be an algberaically closed field of characteristic $p>0$ 
and let $X$ be a Fano threefold over $k$, i.e., $X$ is a three-dimensional smooth projective variety over $k$ such that $-K_X$ is ample. 
Let $r_X$ be the index of $X$ (for its definition, see Definition \ref{d-Fano}). 
Assume that $\rho(X)=1$ and $|-K_X|$ is not very ample. 
Then  $|-K_X|$ is base point free and 
the morphism 
\[
\varphi_{|-K_X|} : X \to \P^{h^0(X, -K_X)-1}
\]
induced by $|-K_X|$ is a double cover onto its image $Y := \varphi_{|-K_X|}(X)$, i.e., 
the induced morphism $X \to Y$ is a finite surjective morphism such that the induced field extension $K(X) \supset K(Y)$ is of degree two. 
Furthermore, one and only one of the following holds. 
\begin{enumerate}
\item $r_X=1, (-K_X)^3 =2$, and $X$ is a double cover of $\P^3$. 
\item $r_X=1, (-K_X)^3 =4$, and 
$X$ is a double cover of a smooth quadric hypersurface in $\P^4$. 
\item 
$r_X=2, (-K_X)^3 = 8,$ and $X$ is isomorphic to a weighted hypersurface of $\P(1, 1, 1, 2, 3)$ of degree $6$. 
\end{enumerate}
\end{thm}


For the case when $|-K_X|$ is very ample, we shall establish the following theorem. 

\begin{thm}[Theorem \ref{t-int-of-quad}]\label{t-main2}
Let $k$ be an algberaically closed field of characteristic $p>0$ 
and let $X$ be a Fano threefold over $k$ such that $|-K_X|$ is very ample. 
Set $g := \frac{1}{2} {\cred (-K_X)^3} +1$. 
If $\rho(X)=1$ and $g \geq 5$, then the image $\varphi_{|-K_X|}(X) \subset \P^{g+1}$ is an intersection of quadrics. 
\end{thm}


Based on Theorem \ref{t-main2}, 
 we shall classify the case when $|-K_X|$ is very ample and $\rho(X)=1$ 
 in the second part \cite{Tan-Fano2}. 
Theorem \ref{t-main1} and Theorem \ref{t-main2} 
are claimed in \cite{SB97}. 
However, \cite{SB97} seems to contain several logical gaps, 
whilst \cite{SB97} introduced many ingenious techniques. 
Some proofs of this paper are base on the ideas introduced in \cite{SB97}. 

\subsection{Overview of the proofs and contents}

Given a Fano threefold $X$, 
it is easy to see that one of (I)--(III) holds (Lemma \ref{l-birat-or-hyperell}). 
\begin{enumerate}
\item[(I)] $|-K_X|$ is not base point free. 
\item[(II)] $|-K_X|$ is base point free and $\varphi_{|-K_X|}$ is birational onto its image. 
\item[(III)] $|-K_X|$ is base point free and $\varphi_{|-K_X|}$ is a double cover onto its image. 
\end{enumerate}
Following the strategy by Mori--Mukai in characteristic zero, 
it is important to classify primitive Fano threefolds, e.g. the case when $\rho(X)=1$. 
In this paper, we focus on the case when $\rho(X)=1$ and $|-K_X|$ is not very ample. 
The main part is to prove that (I) and (II) do not happen under these conditions. 
As for the case (III), we may apply almost the same argument as in characteristic zero. 



\subsubsection{Generic elephants}

In characteristic zero, if $X$ is   a Fano threefold, 
then it is known that a general element of $|-K_X|$, called a general elephant, is smooth. 
It seems to be hard to conclude 
the same conclusion in positive characteristic 
because of the failure of the Bertini theorem for base point free divisors. 
On the other hand, we shall establish the following weaker result. 

\begin{thm}[Corollary \ref{c-generic-ele}]\label{intro-elephant}
Let $k$ be an algebraically closed field of characteristic $p>0$. 
Let $X$ be a Fano threefold over $k$ with $\rho(X)=1$. 
Then the generic member $S$ of $|-K_X|$ is regular and geometrically integral. 
\end{thm}

In particular, the generic member $S$ of $|-K_X|$ is 
a regular prime divisor on $X \times_k \kappa$ 
for the purely transcendental field extension $k \subset \kappa := K(\P(H^0(X, -K_X)))$. 

\begin{rem}
The author does not know whether there exists a Fano threefold in positive characteristic such that any member of $|-K_X|$ is not smooth. 
\end{rem}

\subsubsection{K3-like surfaces} 

Let $X$ be a Fano threefold over $k$ with $\rho(X)=1$. 
By Theorem \ref{intro-elephant}, the generic member $S$ of $|-K_X|$ is regular. 
It is easy to see, by adjunction formula, 
that  $K_S \sim 0$ and $H^1(S, \MO_S)=0$, so that such a surface $S$ is called a K3-like surface. 
Although the base field $\kappa = K(\P(H^0(X, -K_X)))$ of $S$ is no longer perfect, 
we shall see that geometrically integral K3-like surfaces enjoy 
similar properties to those of K3 surfaces, e.g. the following hold: 
\begin{itemize}
\item $H^1(S, -L)=0$ for a nef and big Cartier divisor $L$ on $S$ (Theorem \ref{t-K3-vanishing}). 
\item If $L$ is a Cartier divisor on $S$, then the mobile part of $|L|$ is base point free (Corollary \ref{c-K3-fixed-part}). 
\end{itemize}
Combining with some arguments from \cite{Isk77}, 
we shall settle 
the case when $|-K_X|$ is not base point free, i.e., 
if $X$ is a Fano threefold such that $|-K_X|$ is not base point free, then $\rho(X) \geq 2$ 
(Theorem \ref{t-non-bpf-main}).

\subsubsection{Intersection of quadrics} 

The proof of Theorem \ref{t-main2} is subtler than that of characteristic zero. 
As in characteristic zero, we first consider the intersection $W$ of all the quadrics containing $X$, which can be shown to be a $4$-dimensional variety of minimal degree. 
In characteristic zero, we can easily deduce that $W$ is smooth. 
In positive characteristic, it seems to be hard to obtain the same conclusion, 
whilst we get a weaker conclusion: 
$\dim \Sing\,W \leq 1$. 
We then apply case study depending on $\dim \Sing\,W$. 
The main new ingredients are the following two results. 
\begin{enumerate}
    \item[(i)] Lefschetz hyperplane section theorem for generic members. 
    \item[(ii)] Non-existence of smooth prime divisors on a cone over the Veronese surface. 
\end{enumerate}

To explain how to use (i) and (ii), 
let us focus on the case when 
$\dim \Sing\,W=1$ and $\dim(X \cap \Sing\,W)=0$ (cf. Proposition \ref{p-tri-1dimv}). 
Let $S$ be the generic member of $|-K_X|$. 
By (i), $S \sim -K_{X \times_k \kappa}$ and $S$ is a regular projective surface with $\rho(S)=1$. 
Note that $S$ is disjoint from the blowup centre $\Sing\,W \times_k \kappa$. 
By taking the blowup $W' \to W$ along $\Sing\,W$, we obtain a surjective morphsim $S \to Z \times_k \kappa$, where $W$ is a cone over a variety $Z$ of minimal degree. 
As $\dim \Sing\,W=1$ implies $\dim Z=2$, we see that $1 = \rho(S) \geq \rho(Z \times_k \kappa) \geq \rho(Z)$, and hence $\rho(Z)=1$. 
It is easy to see that such a case only occurs when $Z$ is the Veronese surface. 
However, this contradicts (ii), because $X$ is a smooth prime divisor on 
a cone $W$ over the Veronese surface $Z$. 

\vspace{3mm}

\textbf{Acknowledgements:} 
{\cred The author would like to thank Tatsuro Kawakami 
for constructive comments and answering questions.} 
The author would like to thank 
Natsuo Saito for kindly sharing his private note on \cite{SB97}. 
{\cred 
The author also thanks the referee for reading the manuscript carefully and for suggesting several  improvements.} 
The author was funded by JSPS KAKENHI Grant numbers JP22H01112 and JP23K03028. 
\section{Preliminaries}\label{s-prelim}

\subsection{Notation}\label{ss-notation}

In this subsection, we summarise notation used in this paper. 

\begin{enumerate}
\item We will freely use the notation and terminology in \cite{Har77}. 
In particular, $D_1 \sim D_2$ means linear equivalence of Weil divisors. 
\item 
Throughout this paper, 
we work over an algebraically closed field $k$ 
of characteristic $p>0$ unless otherwise specified. 
\item For an integral scheme $X$, 
we define the {\em function field} $K(X)$ of $X$ 
as the local ring $\MO_{X, \xi}$ at the generic point $\xi$ of $X$. 
For an integral domain $A$, $K(A)$ denotes the function field of $\Spec\,A$. 

\item 
For a scheme $X$, its {\em reduced structure} $X_{\red}$ 
is the reduced closed subscheme of $X$ such that the induced closed immersion 
$X_{\red} \to X$ is surjective. 
\item For a field $\kappa$, 
we say that $X$ is a {\em variety over} $\kappa$ if 
$X$ is an integral scheme that is separated and of finite type over $\kappa$. 
We say that $X$ is a {\em curve} (resp. a {\em surface}, resp. a {\em threefold}) over $\kappa$ 
if $X$ is a variety over $\kappa$ of dimension one (resp. {\em two}, resp. {\em three}). 
\item A variety $X$ over a field is {\em regular} (resp. {\em normal}) if 
the local ring $\MO_{X, x}$ at any point $x \in X$ is regular (resp. an integrally closed domain). 
\item 
For a normal variety $X$ over a field ${\cred \kappa}$, 
we define the {\em canonical divisor} $K_X$ as a Weil divisor on $X$ such that 
$\MO_X(K_X) \simeq \omega_{X/{\kappa}}$, where $\omega_{X/{\cred \kappa}}$ denotes the dualising sheaf (cf. \cite[Section 2.3]{Tan18}). 
Canonical divisors are unique up to linear equivalence. 
Note that $\omega_{X/{\cred \kappa}} \simeq \omega_{X/{\cred \kappa'}}$ for any 
field extension ${\cred \kappa} \subset {\cred \kappa'}$ that induces 
a factorisation $X \to \Spec\,{\cred \kappa'} \to \Spec\,{\cred \kappa}$ (\cite[Lemma 2.7]{Tan18}). 
\item 
Let $X$ be a variety over a field $\kappa$ and let $L$ be an invertible sheaf on $X$. 
Given a finite-dimensional nonzero ${\cred \kappa}$-vector subspace $V \subset H^0(X, L)$, 
$|V|$ denotes the associated linear system and $\varphi_{|V|} : X \dashrightarrow \P^{\dim V-1}$ denotes the induced rational map. 
In particular, we have $|H^0(X, L)| = |L|$. 
Note that the base locus or the base scheme, denoted by $\Bs |V|$, 
is the closed subscheme whose ideal sheaf $I_{\Bs |V|}$ satisfies the following equation: 
\[
{\Im}(V \otimes_k \MO_X \hookrightarrow  H^0(X, L) \otimes_k \MO_X \to L) = I_{\Bs |V|} \cdot L. 
\]
We have the induced morphism: 
\[
\varphi_{|V|}|_{X \setminus \Bs\,|V|} : X \setminus \Bs\,|V| \to \P^{\dim V-1}. 
\]
Note that $V \neq 0$ implies $\Bs |V| \neq X$. Set 
\[
\Im\,\varphi_{|V|} := \overline{\Im (\varphi_{|V|}|_{X \setminus \Bs\,|V|})}, 
\]
where the right hand side is the closure of $\Im (\varphi_{|V|}|_{X \setminus \Bs\,|V|})$ in $\P^{\dim V-1}$. 
\item 
Given a projective variety $X$ and an invertible sheaf $L$, we set $R(X, L) := \bigoplus_{m=0}^{\infty} H^0(X, L^{\otimes m})$, which is a graded $k$-algebra. 
For a Cartier divisor $D$ on $X$, we set $R(X, D) := \bigoplus_{m=0}^{\infty} H^0(X, \MO_X(D)^{\otimes m}) = \bigoplus_{m=0}^{\infty} H^0(X, \MO_X(mD))$. 
\item 
Let $X$ be a projective normal variety over a field. 
We say that a Cartier divisor $D$ is {\em nef} if $D \cdot C \geq 0$ 
for any curve $C$ on $X$. 
Given Cartier divisors $D_1$ and $D_2$, 
the numerical equivalence $D_1 \equiv D_2$ means that $D_1 \cdot C = D_2 \cdot  C$ for any curve  $C$ on $X$. 
For a nef Cartier divisor $D$, 
its {\em numerical dimension} $\nu(X, D)$ is defined as the maximum $\nu \in \{0, 1, ..., \dim X\}$  such that 
$D^{\nu} \cdot H^{\dim X -\nu} \neq 0$ for every (equivalently, for one) ample Cartier divisor $H$. 
\item 
Let $X \subset \P^n$ be  a projective variety which is a closed subscheme of $\P^n$. 
For $s \in \Z_{>0}$, the {\em $s$-th cone} $Y \subset \P^{n+s}$ {\cred (of $X$)} is defined by the same defining ideal of $X$. 
Specifically, if 
\[
X = \Proj\,\frac{k[x_0, ..., x_n]}{(f_1, ..., f_r)} \subset 
\Proj\,k[x_0, ..., x_n] = \P^n 
\]
for homogeneous polynomials $f_1, ..., f_r \in k[x_0, ..., x_n]$, then 
\[
Y = \Proj\,\frac{k[x_0, ..., x_n, y_1, ..., y_s]}{(f_1, ..., f_r)} \subset 
\Proj\,k[x_0, ..., x_n, y_1, ..., y_s] = \P^{n+s}.  
\]
A {\em cone} of $X$ is the $s$-th cone for some $s \in \Z_{>0}$. 
\item 
{\cred Given a projective variety $X$ over a field $\kappa$, $\rho(X)$ denotes its Picard number, i.e., $\rho(X) := \dim_{\Q} ( (\Pic\,X/\equiv) \otimes_{\Z} \Q)$, 
where $\equiv$ is the numerical equivalence. 
}
\end{enumerate}

\subsubsection{Fano threefolds}

\begin{dfn}\label{d-Fano}
We say that $X$ is a {\em Fano threefold}  
if $X$ is a three-dimensional smooth projective variety over $k$ such that $-K_X$ is ample. 
The {\em index} $r_X$ of $X$ is defined as the largest positive integer $r$ such that 
there exists a Cartier divisor $H$ on $X$ {\cred satisfying} 
$-K_X \sim r H$. 
\end{dfn}

\begin{rem}\label{r-d-Fano}
Let $X$ be a Fano threefold. 
We shall frequently use the following notations. 
\begin{enumerate}
\item Set $g := \frac{(-K_X)^3}{2} +1$, which is called the {\em genus} of $X$. 
\item Set $g' := h^0(X, -K_X) -2$. 
\end{enumerate}
If $-K_X$ is very ample, then the genus of $X$ 
coincides with the genus of the smooth curve $H \cap H'$ for general members $H, H' \in |-K_X|$. 
We shall prove the following (Corollary \ref{c-kawakami}):
\[
g' = g + h^1(X, -K_X) \geq g. 
\]
\end{rem}

\subsubsection{Canonical surfaces}

For a field $\kappa$ and a normal surface $S$ over $\kappa$, 
we say that $S$ is {\em canonical} if $(S, 0)$ is canonical in the sense of 
\cite[Definition 2.8]{Kol13}. 
For the minimal resolution $\mu : T \to S$ of $S$, it is well known that the following are equivalent \cite[Theorem 2.29]{Kol13}. 
\begin{enumerate}
\item[(i)] $S$ is canonical. 
\item[(ii)] $K_S$ is Cartier and $K_T \sim \mu^*K_S$. 
\end{enumerate}
In this paper, we only need the characterisation of canonical surfaces by (ii). 
Hence (ii) may be considered as the definition of canonical surfaces. 

\subsection{Cohomologies of Fano threefolds}

Throughout this subsection, we work over an algebraically closed field $k$ of characteristic $p>0$.

{\cred 
\begin{lem}\label{l-kawakami}
Let $X$ be a Fano threefold. 
Then there exist a sequence 
\[
X=:X_0 \xrightarrow{\varphi_0} X_1 \xrightarrow{\varphi_1} \cdots \xrightarrow{\varphi_{\ell-1}} X_{\ell}  
\] 
such that the following properties hold. 
\begin{enumerate}
\item For every $i \in \{0, 1, ..., \ell\}$, $X_i$ is a Fano threefold. 
\item For every $i \in \{0, 1, ..., \ell-1\}$, $\varphi_i : X_i \to X_{i+1}$ 
is a blowup of $X_{i+1}$ along a smooth curve on $X_{i+1}$. 
\item Either 
\begin{enumerate}
    \item $\rho(X_{\ell})=1$ or 
    \item there exists a contraction  $f: X_{\ell} \to Y$ of an extremal ray 
such that $Y$ is a smooth projective rational surface and $f$ is generically smooth. 
\end{enumerate}
\end{enumerate}
\end{lem}

 This lemma is a weaker version of \cite[Lemma 3.2]{Kaw21}. 
Since \cite[Lemma 3.2]{Kaw21} uses a reference \cite{Sai03} which contains a logical gap (Remark \ref{r-Saito}(1)), we give a proof for the sake of completeness. 

\begin{proof}
We may assume that $X$ is primitive, i.e., there exists no blowup $X \to X'$ of a Fano threefold $X'$ along a smooth curve. 
Assume $\rho(X) \geq 2$. 
It is enough to show (b). 
By the same argument as in \cite[Section 8, (8.1), (8.2)]{MM83} (which is applicable by \cite[Main Theorem 1.1]{Kol91}), 
there exists a contraction  $f: X \to Y$ of an extremal ray 
such that $Y$ is a smooth projective surface. 
If $f$ is not generically smooth, then $X$ has a $\P^1$-bundle structure $f' : X \to Y'$ 
which is a contraction of another extremal ray \cite[Corollary 8 and Remark 10]{MS03}. 
Hence we may assume that $f: X \to Y$ is generically smooth. 
By \cite[Corollary 4.10(2)]{Eji19}, $-K_Y$ is big, and hence $Y$ is a smooth ruled surface. 
On the other hand, $Y$ is rationally chain connected, because so is $X$ \cite[Ch. V, Theorem 2.13]{Kol96}. 
Therefore, $Y$ is rational. 
\end{proof}
}

\begin{thm}\label{t-kawakami}
Let $X$ be a Fano threefold. 
Let $D$ be a nef Cartier divisor on $X$ with $\nu(X, D) \geq 2$. 
Then the following hold. 
\begin{enumerate}
    \item $H^j(X, -D)=0$ for any $j \leq 1$. 
    \item $H^i(X, K_X+D)=0$ for any $i \geq 2$. 
    \item $H^i(X, \MO_X)=0$ for any $i>0$. In particular, $\chi(X, \MO_X)=1$.
    \item $\Pic\,X \simeq \Z^{\oplus \rho(X)}$. 
\end{enumerate}
\end{thm}

\begin{proof}
The assertion follows from \cite[Corollary 3.6 and Corollary 3.7]{Kaw21}. 
{\cred 
Note that \cite[Corollary 3.6 and Corollary 3.7]{Kaw21} follows from \cite[Theorem 3.5]{Kaw21}, which depends on \cite[Lemma 3.2]{Kaw21}. 
As mentioned above, \cite[Lemma 3.2]{Kaw21} relies on \cite{Sai03}, which contains a logical gap. 
Although Lemma \ref{l-kawakami} is weaker than \cite[Lemma 3.2]{Kaw21}, 
Lemma \ref{l-kawakami} is enough to establish \cite[Theorem 3.5]{Kaw21}. 
}
\qedhere


\end{proof}

\begin{rem}\label{r-Saito}
\cite[the proof of Lemma 2.4]{Sai03} claims that the same argument as in \cite[Proposition 6.2]{MM83} or \cite[Corollary 4.6]{MM86} works in positive characteristic. 
However, there actually exists an example 
for which the argument does not work  \cite[Remark 7.6]{Tanb}. 
\end{rem}

Given a smooth projective threefold $X$ and a Cartier divisor $D$, 
the Riemann--Roch theorem is given as follows: 
\begin{equation}\label{e-RR3-1}
\chi(X, \MO_X(D)) = \frac{1}{12} D \cdot (D - K_X) \cdot (2D -K_X) + \frac{1}{12} D \cdot c_2(X) + \chi(X, \MO_X) 
\end{equation}
\begin{equation}\label{e-RR3-2}
\chi(X, \MO_X) = \frac{1}{24} (-K_X) \cdot c_2(X).  
\end{equation}

\begin{cor}\label{c-kawakami}
Let $X$ be a  Fano threefold. 
Let $H$ be a Cartier divisor such that the numerical equivalence $H \equiv -q K_X$ holds 
for some $q \in \Q_{\geq 0}$ (i.e., $H \cdot C = -q K_X \cdot C$ holds for any curve $C$ on $X$). 
Then the following holds 
\[
h^0(X, H) -h^1(X, H)
= \chi(X, H) =  \frac{1}{12} q(q+1)(2q+1)(-K_X)^3 + 2q + 1. 
\]
In particular, if $m \in \Z_{\geq 0}$, then 
\[
h^0(X, -mK_X) -h^1(X, -mK_X)
= \chi(X, -mK_X) =  \frac{1}{12} m(m+1)(2m+1)(-K_X)^3 + 2m + 1, 
\]
and hence $h^0(X, -K_X) -h^1(X, -K_X) =  \frac{(-K_X)^3}{2} + 3 = g+2$. 
\end{cor}

\begin{proof}
If $q=0$, then there is nothing to show (Theorem \ref{t-kawakami}). 
Assume $q >0$. 
We then have $\nu(X, H) = \nu(X, -K_X) =3$.
By (\ref{e-RR3-2}) and $H \equiv -q K_X$,  
we have 
\begin{equation}\label{e1-c-kawakami}
H \cdot c_2(X) = q (-K_X) \cdot c_2(X) = 24q \chi(X, \MO_X)=24q. 
    \end{equation}
Then the assertion holds by the following computation: 
\begin{align*}
h^0(X, H) -h^1(X, H)
&\overset{(1)}{=} \chi(X, H) \\
&\overset{(2)}{=}  \frac{1}{12} H \cdot (H - K_X) \cdot (2H -K_X) + \frac{1}{12} H \cdot c_2(X) + \chi(X, \MO_X) \\
&\overset{(3)}{=} \frac{1}{12} (-qK_X) \cdot (-qK_X -K_X) \cdot (-2q K_X -K_X) + 2q + 1\\
&= \frac{1}{12} q(q+1)(2q+1)(-K_X)^3 + 2q + 1.\\
\end{align*}
Here (1) holds by Theorem \ref{t-kawakami}, 
(2) follows from (\ref{e-RR3-1}), and 
we obtain (3) by (\ref{e1-c-kawakami}) and $\chi(X, \MO_X)=1$ (Theorem \ref{t-kawakami}). 
\end{proof}

\subsection{Bertini theorems}

\subsubsection{Generic members}

\begin{nota}\label{n-generic-member}
Let $\kappa$ be a field. 
Let $X$ be a variety over $\kappa$, 
let $L$ be a Cartier divisor on $X$, and 
let $V \subset H^0(X, L)$ be a nonzero finite-dimensional ${\cred \kappa}$-vector subspace. 
We then have the universal family $X^{\univ}_{L, V}$ that parameterises the effective Cartier divisors 
of $|V|$ whose parameter space is $\P(V)$. 
Then its generic fibre $X^{\gen}_{L, V}$ is called the {\em generic member} of $|V|$. 
To summarise, we have the following diagram in which all the squares are cartesian: 
\[
\begin{CD}
X^{\gen}_{L, V} @>>> X^{\univ}_{L, V}\\
@VVV @VVV\\
X \times_{\kappa} K(\P(V)) @>>> X \times_{\kappa} \P(V)\\
@VVV @VVV\\
\Spec\,K(\P(V)) @>>> \P(V).  
\end{CD}
\]
Note that if $X$ is normal, then we have $X^{\gen}_{L, V} \sim L_{K(\P(V))}$ for the pullback $L_{K(\P(V))}$ of $L$ to  $X  \times_{\kappa} K(\P(V))$. 
\end{nota}

\begin{rem}\label{r-generic-member}
We use Notation \ref{n-generic-member}. 
Let $X'$ be a non-empty open subset of $X$ and let 
$i: X' \hookrightarrow X$ be the induced open immersion. 
We then obtain 
\[
V \subset H^0(X, L) \hookrightarrow H^0(X', L')\qquad \text{for}\qquad 
L' := L|_{X'}. 
\]
Set $V'$ to be the image of $V$ in $H^0(X', L')$. 
We then obtain the following cartesian diagram \cite[Proposition 5.10]{Tana}: 
\[
\begin{CD}
X'^{\gen}_{L', V'} @>>> X^{\gen}_{L, V}\\
@VVV @VVV\\
X' @>i>> X. 
\end{CD}
\]
\end{rem}

\begin{thm}\label{t-Bertini-generic}
Let $\kappa$ be a field. 
Let $X$ be a variety over $\kappa$, 
let $L$ be a Cartier divisor on $X$, and 
let $V \subset H^0(X, L)$ be a nonzero finite-dimensional ${\cred \kappa}$-vector subspace. 
Set $X' := X \setminus \Bs\,|V|$ and let 
$\alpha : {\cred X^{\gen}_{L, V}}
\to X$ be the induced morphism. 
If $X'$ is regular, then $\alpha^{-1}(X')$ is regular. 
\end{thm}

\begin{proof}
After replacing $X$ by $X'$, we may assume that $|V|$ is base point free (Remark \ref{r-generic-member}). 
Then the assertion follows from  \cite[Theorem 4.9 and Remark 5.8]{Tana}. 
\end{proof}

\subsubsection{General members}

\begin{prop}\label{p-Bertini-irre}
We work over an algebraically closed field $k$. 
Let $f: X \to \P^N$ be a morphism from a variety $X$. 
If $\dim \overline{f(X)} \geq 2$, then $f^{-1}(H)$ is irreducible 
for a general hyperplane $H \subset \P^N$. 
\end{prop}

\begin{proof}
See \cite[4) of Theoreme 6.3]{Jou83}. 
\end{proof}

The following proposition is a known result \cite[Theorem 2.7]{Fuj90}. 
However, we here give a proof for the sake of compleness, 
because \cite[Theorem 2.7]{Fuj90} depends on \cite{Wei62} which is written 
in a classical language of algebraic geometry.

\begin{prop}\label{p-Bertini-integral}
We work over an algebraically closed field $k$. 
Let $X$ be a projective normal variety and 
let $L$ be a Cartier divisor on $X$ 
such that $\dim \Bs\,|L| \leq \dim X-2$ and 
$\dim(\Im\,\varphi_{|L|}) \geq 2$. 
Then general members of $|L|$ are prime divisors. 
In particular, the generic member of $|L|$ is geometrically integral. 
\end{prop}

\begin{proof}
We first reduce the problem to the case when $|L|$ is base point free. 
Let 
\[
\mu : X' \to X 
\] 
be the normalisation of the resolution of the indeterminacies of $\varphi_{|L|} : X \dashrightarrow \P^{h^0(X, L) -1}$. 
We then have 
\[
\mu^*L = M + F, 
\]
where $|M|$ is base point free and $F$ is the fixed part of $|\mu^*L|$. 
Then we obtain the induced morphisms: 
\[
\varphi_{|M|} : X' \xrightarrow{\mu} X \overset{\varphi_{|L|}}{\dashrightarrow} \P^{h^0(X, L) -1}. 
\]
By construction, we have 
\[
H^0(X', M) \xrightarrow{\simeq} H^0(X', M+F) = H^0(X', \mu^*L) \simeq H^0(X, L). 
\]
Furthermore, this composite isomorphism induces the bijection: $|M| \to |L|, D \mapsto \mu_*D$. 
Therefore, if a general member $D$ of $|M|$ is a prime divisor, then 
also its push-forward $\mu_*D$ is a prime divisor. 
Hence we may assume that $|D|$ is base point free 
after replacing $(X, L)$ by $(X', M)$.

We have the following induced morphisms: 
\[
\varphi_{|L|} : X \xrightarrow{\varphi} Y \hookrightarrow \P^N, \qquad Y := \varphi_{|L|}(X), \qquad N := h^0(X, L) -1. 
\]
By construction, it holds that 
\[
H^0(X, L) \simeq H^0(Y, \MO_{\P^N}(1)|_Y) \simeq H^0(\P^N, \MO_{\P^N}(1)). 
\]
In particular, a general member $H_Y \in |\MO_{\P^N}(1)|_Y|$ (a general hyperplane section) 
is an integral scheme. 

For a suitable non-empty open subset $Y^{\circ} \subset Y$, which is normal, and 
its inverse image $X^{\circ} := \varphi^{-1}(Y^{\circ}) \subset X$, 
the Noether normalisation theorem, applied for the generic fibre $X \times_Y \Spec\,K(Y)$, 
induces following factorisation: 
\[
X^{\circ} \xrightarrow{\alpha} Z := Y^{\circ} \times \P^r \xrightarrow{{\rm pr}_1} Y^{\circ}, 
\]
where $\alpha : X^{\circ} \to Z = Y^{\circ} \times \P^r$ is a finite surjective morphism 
of normal varieties and the latter morphism 
{\cred ${\rm pr}_1 : Y^{\circ} \times \P^r \to Y^{\circ}$} is the projection (cf. \cite[Lemma 2.14]{Tan20}). 
For a general hyperplane section $H_Y \subset Y$, 
its inverse image $H_Y|_Z$ to $Z$ is an integral scheme, 
because  it can be written as $(H_Y \cap Y^{\circ}) \times \P^r$. 
Furthermore, we take the decomposition via the separable closure: 
\[
X^{\circ} \xrightarrow{\beta} W \xrightarrow{\gamma} Z = Y^{\circ} \times \P^r \xrightarrow{{\rm pr}_1} Y^{\circ}, 
\]
where 
\begin{itemize}
\item $W$ is a normal variety, 
\item $\beta : X^{\circ} \to W$ is a finite surjective purely inseaprable morphism, and 
\item $\gamma : W \to Z$ is a finite surjective separable morphism. 
\end{itemize}
Then also the pullback  $H_Y|_W$  of $H_Y$ to $W$ is integral (i.e., a prime divisor). 
Indeed, $H_Y|_X$ is irreducible (Proposition \ref{p-Bertini-irre}), so it suffices to prove that $H_Y|_W$ is reduced. 
It is $S_1$, and hence it suffices to prove that $H_Y|_W$ is $R_0$. 
This holds because we may assume that $H_Y|_Z$ intersects with the \'etale locus of $\gamma : W \to Z$. 

Since ${\cred \beta :} X^{\circ} \to W$ is a finite purely inseparable morphism, we have the following factorisation for some $e \in \Z_{>0}$: 
\[
F^e : W \to X^{\circ} \xrightarrow{{\cred \beta}} W, 
\]
{\cred where $F^e : W \to W$ denotes the $e$-th iterated absolute Frobenius morphism.} 
By $(F^e)^*(H_Y|_W) = p^e H_Y|_W$, 
it holds that $H_Y|_{X^{\circ}} = p^d (H_Y|_{X^{\circ}})_{\red}$ for some $0 \leq d \leq e$. 
A general member of $|L|$ is nothing but the pullback $H_Y|_X$ of 
a general member $H_Y$ of $|\MO_{\P^N}(1)|_Y|$. 
In particular, 
\[
(H_Y|_X)|_{X^{\circ}} = H_Y|_{X^{\circ}} =p^d (H_Y|_{X^{\circ}})_{\red}. 
\]
Suppose that a general member of $|L|$ is not ingeral. 
Since $H_Y|_X$ is irreducible  (Proposition \ref{p-Bertini-irre}), 
we have $H_Y|_X = pD$ for some effective Weil divisor $D$. 

Therefore, a general member of $|L|$ is of the form $pD$. 
Replacing $L$ by $pD$, we may assume that $L = pL'$ for some Weil divisor $L'$. 
This implies that the image of 
\[
F : H^0(X, L') \to H^0(X, L), \qquad s \mapsto s^p
\]
is dense with respect to the Zariski topology of the affine space $H^0(X, L)$. 
As the image $\Im (F)$ is a closed subset in $H^0(X, L)$, 
$F$ is surjective. 
Therefore, any member of $|L|$ can be written as $pD$ for some Weil divisor $D$.

We have an isomorphism: 
\[
F: H^0(X, L') \xrightarrow{\simeq} H^0(X, L), \qquad t \mapsto t^p. 
\]
Fix a $k$-linear basis of $H^0(X, L')$:  
\[
t_0, ..., t_N \in H^0(X, L')
\]
and their images to $H^0(X, L)$: 
\[
t_0^p, ..., t_N^p \in H^0(X, L),  
\]
which form a $k$-linear basis of $H^0(X, L)$. 
We then obtain 
\[
\varphi_{|L|} : X \to \P^N, \qquad x \mapsto [ (t_0(x))^p : \cdots : (t_N(x))^p]. 
\]
Therefore, we get the following factorisation: 
\[
\varphi_{|L|} : X \xrightarrow{\varphi_{|L'|} } \P^N \xrightarrow{F} \P^N, \qquad 
x \mapsto [t_0(x) : \cdots : t_N(x)] \mapsto [ (t_0(x))^p : \cdots : (t_N(x))^p].
\]
By $Y = \varphi_{|L|}(X)$, we obtain the factorisation 
\[
X \to Y \xrightarrow{F} Y, 
\]
which induces the following isomorphisms 
\[
H^0(Y, \MO_{\P^N}(1)|_Y) \xrightarrow{\simeq} H^0(Y, {\cred \MO_{\P^N}(p)|_Y}) \xrightarrow{\simeq}H^0(X, L), 
\]
because the composition is an isomorphism and each map is injective. 
This implies that any member of $|\MO_{\P^N}(p)|_Y|$ is non-reduced. 
However, this is a contradiction, because 
$H \cap Y \in |\MO_{\P^N}(p)|_Y|$ is an integral scheme for a general member $H \in |\MO_{\P^N}(p)|$. 
\end{proof}

\begin{lem}\label{l-Bertini-birat}
We work over an algebraically closed field $k$. 
Let $X$ be a projective CM variety with $\dim X \geq 2$ and 
let 
\[
\psi : X \to Y 
\]
be a birational morphism to a projective variety $Y$. 
Fix a closed embedding $Y \subset \P^N$. 
If $H$ is a general hyperplane of $\P^N$, then 
\[
\psi^{-1}(Y \cap H) \to Y \cap H
\]
is a birational morphism 
from a projective CM variety $\psi^{-1}(Y \cap H)$ 
to a projective variety $Y \cap H$. 
In particular, if $H_1, ..., H_r \subset \P^N$ are general hyperplanes with $1 \leq r \leq \dim X-1$, then 
\[
\psi^{-1}(Y \cap H_1 \cap \cdots \cap H_{r})
\to Y \cap H_1 \cap \cdots \cap H_{r}
\]
is a birational morphism from a projective CM variety 
$\psi^{-1}(Y \cap H_1 \cap \cdots \cap H_{r})$ to a projective variety 
$Y \cap H_1 \cap \cdots \cap H_{r}$. 
\end{lem}

\begin{proof}
Fix a general hyperplane $H$. 
By the classical Bertini theorem, $Y \cap H$ is an integral scheme. 
It follows from Proposition \ref{p-Bertini-irre} that $\psi^{-1}(Y \cap H)$ is irreducible. 
Since $\psi^{-1}(Y \cap H)$ is an effective Cartier divisor on $X$, 
also $\psi^{-1}(Y \cap H)$ is CM. 
Since $\psi^{-1}(Y \cap H) \to Y \cap H$ is isomorphic over the generic point of $Y \cap H$, 
$\psi^{-1}(Y \cap H)$ is generically reduced. 
Then $\psi^{-1}(Y \cap H)$ is reduced, because $\psi^{-1}(Y \cap H)$ is $R_0$ and $S_1$.  
Therefore, $\psi^{-1}(Y \cap H)$ is an integral scheme. 
\end{proof}

\subsection{$\Delta$-genera}

Throughout this subsection, we work over an algebraically closed field $k$ of characteristic $p>0$. 


\begin{dfn}\label{d-Delta}
We say that  $(X, L)$ is a {\em polarised variety} 
if $X$ is a projective variety and $L$ is an ample invertible sheaf or an ample Cartier divisor.  
Set 
\[
\Delta(X, L) := \dim X + L^{\dim X} -h^0(X, L). 
\]
We say that polarised varieties $(X, L)$ and $(X', L')$ are {\em isomorphic}, 
denoted by $(X, L) \simeq (X', L')$, 
if there exists a $k$-isomorphism $\theta : X \xrightarrow{\simeq} X'$ such that $L \simeq \theta^*L'$. 
A polarised variety $(X, L)$ is called {\em smooth} if $X$ is smooth. 
\end{dfn}

\begin{thm}\label{t-Delta-geq0}
Let $(X, L)$ be a polarised variety. 
Then the following hold. 
\begin{enumerate}
    \item $\Delta(X, L) \geq 0$, i.e., 
$h^0(X, L) 
\leq \dim X + L^{\dim X}$. 
\item If $\Delta(X, L)=0$, then $|L|$ is very ample. 
\end{enumerate}
\end{thm}

\begin{proof}
See \cite[Theorem (2.1) and Theorem (4.2)]{Fuj82b}. 
\end{proof}

\begin{thm}\label{t-Delta0-sm}
Let $(X, L)$ be a smooth polarised variety such that $\Delta(X, L)=0$. 
Set $n := \dim X$. 
Then one of the following holds. 
\begin{enumerate}
\item $L^n =1$ and $(X, L) \simeq (\P^n, \MO_{\P^n}(1))$. 
\item $L^n =2$ and $(X, L) \simeq (Q^n, \MO_{\P^{n+1}}(1)|_{Q^n})$, 
where $Q^n \subset \P^{n+1}$ is a smooth quadric hypersurface. 
\item $L^n \geq 3$ and $X \simeq \P_{\P^1}(E)$, where $E$ is a locally free sheaf on $\P^1$ of rank $n$. 
\item $(X, L) \simeq (\P^2, \MO_{\P^2}(2))$ and $L^2 =4$. 
\end{enumerate}
\end{thm}

\begin{proof}
See \cite[Corollary (4.3) and Theorem (4.9)]{Fuj82b}. 
\end{proof}

\begin{rem}\label{r-Delta0-sing}
Let $(X, L)$ be a polarised variety such that $X$ is not smooth and  $\Delta(X, L)=0$. 
Then $|L|$ is very ample (Theorem \ref{t-Delta-geq0}(2)). 
Set $N := h^0(X, L) -1$ and we fix a closed embedding $X \subset \P^N$ induced by $|L|$. 
Then $\Delta(X, L) = 0$ implies $\deg X = 1 + {\rm codim}\,X$. 
Such a variety is called a {\em variety of minimal degree}, 
since the inequality $\deg X \geq 1 + {\rm codim}\,X$ holds in general. 
By \cite[Theorem (4.11)]{Fuj82b} or \cite[Theorem 1]{EH87}, 
it is known that $(X, L)$ is a cone over a smooth polarised variety $(Z, L_Z)$ with $\Delta(Z, L_Z)=0$. 
More specifically, it holds by \cite[page 4]{EH87} that 
\begin{itemize}
\item $Z  = \Proj\,k[x_0, ..., x_{N-r}] /(f_1, ..., f_s) \subset \P^{N-r}$ 
for $r := \dim X - \dim Z$ and some $f_1, ..., f_s \in k[x_0, ..., x_{N-r}]$, 
\item $X = \Proj\,k[x_0, ..., x_{N}] /(f_1, ..., f_s) \subset \P^{N}$, i.e., 
the same equations as those of $Z$ define $X \subset \P^N$. 
\end{itemize}
\end{rem}

\begin{rem}\label{r-Delta-geom-int}
{\cred 
Let $\kappa$ be a (not necessarily algebraically closed) field. 
For a geometrically integral projective variety $X$ over $\kappa$ and an ample divisor $L$, we set 
\[
\Delta(X, L) := \dim X +L^{\dim X} -h^0(X, L). 
\]
By definition, we have 
\[
\Delta(X, L) = \Delta(X \times_{\kappa} \overline{\kappa}, \alpha^*L) \geq 0, 
\]
where $\overline \kappa$ denotes the algebraic closure of $\kappa$ and 
$\alpha : X \times_{\kappa} \overline{\kappa} \to X$ is the induced morphism.}
\end{rem}
\subsection{The case of index $\geq 2$}

Throughout this subsection, we work over an algebraically closed field $k$ of characteristic $p>0$.

\begin{thm}\label{t-index-34}
Let $X$ be a Fano threefold. 
Then the following hold. 
\begin{enumerate}
\item $1 \leq r_X \leq 4$. 
\item $r_X=4$ if and only if $X \simeq \P^3$. 
\item $r_X=3$ if and only if $X$ is isomorphic to a smooth quadric hypersurface in $\P^4$. 
\end{enumerate}
\end{thm}

\begin{proof}
See \cite[Proposition 4]{Meg98}. 
\end{proof}

\begin{prop}\label{p-wt-lift}
Let $(X, L)$ be a polarised variety. 
Fix $a \in \Z_{>0}$ and $\delta \in H^0(X, aL)$. 
Let $D$ be the member of $|aL|$ corresponding to $\delta$. 
Let 
\[
\rho : \bigoplus_{m=0}^{\infty} H^0(X, mL) \to 
\bigoplus_{m=0}^{\infty} H^0(D, m(L|_D)) 
\]
be the restriction map. 
Let $\xi_1, ..., \xi_r$ be homogeneous elements of $\bigoplus_{m=0}^{\infty} H^0(X, mL)$. 
If $\bigoplus_{m=0}^{\infty} H^0(D, m(L|_D))$ is generated by $\rho(\xi_1), ..., \rho(\xi_r)$ 
as a $k$-algebra, 
then $\bigoplus_{m=0}^{\infty} H^0(X, mL)$ is generated by $\delta, \xi_1, ..., \xi_r$. 
\end{prop}

\begin{proof}
See \cite[Theorem (2.3)]{Fuj90}. 
\end{proof}

{\cred 
\begin{dfn}[(1.5) and (1.6) of \cite{Fuj82b}]
Let $(X, L)$ be a polarised variety and set $n:= \dim X$. 
\begin{enumerate}
\item 
We say that 
$X =:X_n \supset X_{n-1} \supset 
\cdots \supset X_1$ 
is a {\em ladder} (of $(X, L)$) if, for every $i \in \{n, n-1, ..., 2\}$, 
$X_{i-1}$ is a member of $|L|_{X_i}|$ which is an integral scheme. 
By definition, we have $\dim X_i = i$ for each $i \in \{n, n-1, ..., 1\}$. 
\item 
We say that 
$X =:X_n \supset X_{n-1} \supset 
\cdots \supset X_1$ is a {\em regular ladder} of $(X, L)$ if 
this is a ladder and the restriction map 
\[
H^0(X_i, L|_{X_i}) \to H^0(X_{i-1}, L|_{X_{i-1}})
\]
is surjective for every $i \in \{n, n-1, ..., 2\}$. 
We say that $(X, L)$ {\em has a regular ladder} if there exists a regular ladder of $(X, L)$. 
\end{enumerate}
\end{dfn}
}

\begin{thm}\label{t-index2-Delta}
Let $X$ be a Fano threefold with $r_X=2$. 
Let $H$ be a Cartier divisor such that $-K_X \sim 2H$. 
Then the following hold. 
\begin{enumerate}
\item $\Delta(X, H)=1$. 
\item ${\cred (X, H)}$ has a regular ladder. 
\item If $H^3 \geq 2$, then $|H|$ is base point free. 
\item If $H^3 \geq 3$, then $|H|$ is very ample. 
\end{enumerate}
\end{thm}

\begin{proof}
The assertion (1) follows from \cite[Lemma 5]{Meg98}. 
Then (2)--(4) hold by \cite[Corollary (5.5)]{Fuj82b}. 
\end{proof}

\begin{lem}\label{l-ell-wt}
Let $C$ be a projective Gorenstein curve such that $\omega_C \simeq \MO_C$. 
Let $L$ be a Cartier divisor on $C$. 
Then the following hold. 
\begin{enumerate}
\item 
If $\deg L \geq 3$, then $R(X, L)$ is generated by $H^0(X, L)$ as a $k$-algebra. 
\item 
If $\deg L = 2$, then $R(X, L)$ is generated by 
$\eta_1, \eta_2 \in H^0(X, L)$ and $\eta_3 \in H^0(X, 2L)$ as a $k$-algebra. 
\item 
If $\deg L=1$, then 
$R(X, L)$ is generated by 
$\eta_1 \in H^0(X, L), \eta_2 \in H^0(X, 2L), \eta_3 \in H^0(X, 3L)$
 as a $k$-algebra. 
\end{enumerate}
\end{lem}

\begin{proof}
By \cite[Proposition 11.11]{Tan21}, the following hold. 
\begin{enumerate}
\item[(i)] If $\deg L \geq 3$, then $R(X, L)$ is generated by $H^0(X, L)$.
\item[(ii)] If $\deg L = 2$, then $R(X, L)$ is generated by $H^0(X, L) \oplus H^0(X, 2L)$. 
\item[(iii)] If $\deg L = 1$, then $R(X, L)$ is generated by 
$H^0(X, L) \oplus H^0(X, 2L) \oplus H^0(X, 3L)$. 
\end{enumerate}
Then (i) implies (1). 

We only prove (3), as the proof of (2) is similar. 
By the Riemann--Roch theorem, we have $h^0(X, mL) =m$ for any $m \geq 1$. 
Fix a nonzero element $\eta_1 \in H^0(X, L)$, so that we have 
$H^0(X, L) = k \eta_1$. 
Then $\eta_1^2 \in H^0(X, 2L)$ is nonzero, and hence we can find $\eta_2 \in H^0(X, 2L)$ such that 
$H^0(X, 2L) = k \eta_1^2 \oplus k \eta_2$. 
Since 
\[
\times \eta_1 : H^0(X, 2L) \to H^0(X, 3L), \qquad s \mapsto \eta_1 s
\]
is injective, $\eta_1^3, \eta_1\eta_2 \in H^0(X, 3L)$ are linearly independent over $k$. 
Therefore, we can find $\eta_3 \in H^0(X, 3L)$ such that 
\[
H^0(X, 3L) = k \eta_1^3 \oplus k \eta_1 \eta_2 \oplus k \eta_3. 
\]
Hence $\eta_1, \eta_2, \eta_3$ generate $R(X, L)$ as a $k$-algebra. 
Thus (3) holds. 
\end{proof}

\begin{thm}\label{t-index-2}
Let $X$ be a Fano threefold with $r_X=2$. 
Let $H$ be a Cartier divisor such that $-K_X \sim 2H$. 
Then $1 \leq H^3 \leq 7$. 
Furthermore, the following hold. 
\begin{enumerate}
\item If $H^3 =1$, then $X$ is isomorphic to a weighted hypersurface in $\P(1, 1, 1, 2, 3)$ of degree $6$. 
\item If $H^3 =2$, then $X$ is isomorphic to a weighted hypersurface in $\P(1, 1, 1, 1, 2)$ of degree $4$. 
\item If $H^3 = 3$, then $X$ is isomorphic to a cubic hypersurface in $\P^4$. 
\item If $H^3= 4$, then $X$ is isomorphic to a complete intersection of two quadrics in $\P^5$. 
\item If $H^3 =5$, then $X$ is isomorphic to ${\rm Gr}(2, 5) \cap H_1 \cap H_2 \cap H_3$ in $\P^9$, 
where the inclusion ${\rm Gr}(2, 5) \subset \P^9$ is the Pl\"ucker embedding and 
$H_1, H_2, H_3$ are general hyperplanes of $\P^9$. 
\item If $H^3=6$, then $X$ is isomorphic to either $\P^1 \times \P^1 \times \P^1$ or 
a member of $|\MO_{\P^2 \times \P^2}(1, 1)|$ on $\P^2 \times \P^2$. 
\item If $H^3 =7$, then $X$ is isomorphic to a blowup of $\P^3$ at a point. 
\end{enumerate}
\end{thm}

\begin{proof}
It follows from \cite[Theorem 6]{Meg98} that (3)--(7) holds. 
Let us show (1) and (2). 
By Theorem \ref{t-index2-Delta}, there exists a regular ladder {\cred of $(X, H)$}: 
\[
X \supset S \supset C, 
\]
and hence the restriction map 
\[
H^0(X, H) \to H^0(C, H_C), \qquad H_C := H|_C
\]
is {\cred surjective}. 

We now prove (1). 
By Proposition \ref{p-wt-lift}, it is enough to find 
\[
\eta_1 \in H^0(C, H_C),\quad 
\eta_2 \in H^0(C, 2H_C), \quad 
\eta_3 \in H^0(C, 3H_C)
\]
such that $R(C, H_C)$ is generated by $\eta_1, \eta_2, \eta_3$. 
This follows from Lemma \ref{l-ell-wt}(3). 
Thus (1) holds. 
By Proposition \ref{p-wt-lift} and Lemma \ref{l-ell-wt}(2), 
we can apply a similar argument for (2) to that of (1). 
\end{proof}

\begin{cor}\label{c-index2-va}
Let $X$ be a Fano threefold with $r_X\geq 2$. 
Then the following hold. 
\begin{enumerate}
\item If $(r_X, (-K_X/2)^3)  \neq (2, 1)$, then $R(X, -K_X)$ is generated by $H^0(X, -K_X)$, 
and hence $|-K_X|$ is very ample. 
\item If $(r_X, (-K_X/2)^3)  = (2, 1)$, then 
$|-K_X|$ is base point free and $\varphi_{|-K_X|}$ is a double cover onto its image 
(in particular, $|-K_X|$ is not very ample). 
\end{enumerate}
\end{cor}

\begin{proof}
Let us show (1). 
If $r_X \geq 3$, then the assertion holds by Theorem \ref{t-index-34}. 
We may assume that $r_X=2$. 
Fix an ample Cartier divisor $H$ on $X$ with $-K_X \sim 2H$. 
If $H^3 \geq 3$, then $|H|$ is very ample, and hence also $|-K_X| = |2H|$ is very ample. 
Thus we may assume that $(r_X, H^3)= (2, 2)$. 
In this case, $|H|$ is base point free. 
By the proof of Theorem \ref{t-index-2}, 
$R(X, H)$ is generated by $x_0, x_1, x_2, x_3 \in H^0(X, H)$ and $y \in H^0(X, 2H)$ as a $k$-algebra. 
This immediately deduces that $R(X, 2H)$ is generated by $H^0(X, 2H)$, 
because if a monomial $x_0^{a_0}x_1^{a_1}x_2^{a_2}x_3^{a_3}y^b$ is of even degree, i.e., 
$a_0 + a_1 + a_2+a_3 + 2 b \in 2\Z_{>0}$, then we get $a_0 + a_1 + a_2+a_3 \in 2\Z_{\geq 0}$, and hence 
$x_0^{a_0}x_1^{a_1}x_2^{a_2}x_3^{a_3}y^b$ can be written as a multiple of elements of $H^0(X, 2H)$. 
Therefore, $|-K_X| = |2H|$ is very ample. 

The assertion (2) follows from \cite[Theorem 8]{Meg98}. 
\end{proof}

\section{K3-like surfaces}\label{s-K3-like}

Throughout this section, we work over a field $\kappa$ of characteristic $p>0$. 
In our applications, we have $\kappa = k (t_1, ..., t_N)$, which is a purely transcendental field extension over an algebraically closed field $k$ of characteristic $p>0$. 
In particular, $\kappa$ can not be assumed to be perfect. 

\begin{dfn}
We work over a field $\kappa$ of characteristic $p>0$.  
We say that $S$ is a {\em K3-like surface} 
if $S$ is a projective normal surface such that $H^1(S, \MO_S)=0$ and $K_S \sim 0$. 
\end{dfn}

We are mainly interested in the following two cases: 
geometrically integral regular K3-like surfaces and 
geometrically integral canonical K3-like surfaces, i.e., $S$ has at worst canonical singularities. 

\subsection{Vanishing theorems}


The purpose of this subsection is to establish the Kodaira (or Ramnujam) vanishing theorem for geometrically integral K3-like surfaces (Theorem \ref{t-K3-vanishing}). 
To this end, we shall establish a vanishing of Mumford type for geometrically integral surfaces (Theorem \ref{t-Mumford-vanishing}). 
We start with the following auxiliary result. 

\begin{prop}\label{p-Enokizono}
We work over a field $\kappa$ of characteristic $p>0$. 
Let $S$ be a projective normal surface and let $D$ be a nef and big effective $\Q$-Cartier $\Z$-divisor. 
Then $H^0(D, \MO_D)$ is a field and 
\[
H^0({\cred S}, \MO_{\cred S}) \hookrightarrow H^0(D, \MO_D)
\]
is a finite purely inseparable extension. 
\end{prop}

\begin{proof}
{\cred 
For $\kappa' := H^0(S, \MO_S)$,  we have the following induced morphisms: 
\[
S \to \Spec\,\kappa' \to \Spec\,\kappa. 
\]
Recall that $\kappa'$ is a field which is a finitely generated $\kappa$-module \cite[Ch. III, Theorem 8.8(b)]{Har77}, i.e., 
$\kappa \subset \kappa'$ is a field extension of finite degree.}
Replacing $\kappa$ by ${\cred \kappa' =} H^0(S, \MO_S)$, we may assume that $H^0(S, \MO_S)=\kappa$. 
Then the assertion follows from \cite[Collorary 3.13 and Corollary 3.17]{Eno}. 
\end{proof}

\begin{thm}\label{t-Mumford-vanishing}
We work over a field $\kappa$ of characteristic $p>0$. 
Let $S$ be a geometrically integral projective normal surface. 
Let $D$ be a nef and big effective $\Q$-Cartier $\Z$-divisor such that 
$h^0(S, \MO_S(D)) \geq 2$. 
Then the following hold. 
\begin{enumerate}
\item The induced map $H^0(S, \MO_S) \to H^0(D, \MO_D)$ is an isomorphism. 
\item The induced map $H^1(S, \MO_S(-D)) \to H^1(S, \MO_S)$ is injective. 
\end{enumerate}
\end{thm}

\begin{proof}
{\cred Since $S$ is geometrically integral, we have $H^0(S, \MO_S) = \kappa$.} 
As the assertion (2) immediately holds by (1), it suffices to show (1). 
Taking the base change to the separable closure of $\kappa$, 
we may assume that $\kappa$ is separably closed. 
Since $S$ is geometrically integral, there exists a ${\cred \kappa}$-rational point $P \in S$ around which $S$ is smooth. 
In particular, $\MO_S(D)$ is invertible around $P \in S$. 
We then obtain a short exact sequence  
\[
0 \to \MO_S(D) \otimes \m_P \to \MO_S(D) \to \MO_S(D) \otimes \MO_{P} \to 0, 
\]
which induces another exact sequence: 
\[
0 \to H^0(S, \MO_S(D) \otimes \m_P) \to H^0(S, \MO_S(D)) \to H^0(P, \MO_S(D)|_P) \simeq {\cred \kappa}. 
\]
By $h^0(S, \MO_S(D)) \geq 2> 1=h^0(P, \MO_S(D)|_P)$, we get $H^0(S, \MO_S(D) \otimes \m_P) \neq 0$. 
Hence we may assume that $P \in \Supp\,D$. 
By Proposition \ref{p-Enokizono}, we have field extensions 
\[
{\cred \kappa}= H^0(S, \MO_S) \hookrightarrow H^0(D, \MO_D) \hookrightarrow H^0(P, \MO_P) = {\cred \kappa}
\]
whose composite is identity. 
Therefore, $H^0(S, \MO_S) \to H^0(D, \MO_D)$ is an isomorphism. 
Thus (1) holds. 
\end{proof}

\begin{thm}\label{t-K3-vanishing}
We work over a field $\kappa$ of characteristic $p>0$. 
Let $S$ be a geometrically integral regular K3-like surface. 
Let $L$ be a nef and big Cartier divisor.  
Then the following hold. 
\begin{enumerate}
\item $H^j(S, \MO_S(-L))=0$ for any $j <2$. 
\item $H^i(S, \MO_S(L)) =0$ for any $i>0$. 
\item $h^0(S, \MO_S(L)) = 2 + \frac{1}{2}L^2$. 
\end{enumerate}
\end{thm}

\begin{proof}
By Serre duality, we have $H^2(S, L)=0$. 
It follows from the Riemann--Roch theorem that 
\[
h^0(S, L) \geq 
h^0(S, L) -h^1(S, L) = \chi(S, L) = \chi(S, \MO_S) + \frac{1}{2} L \cdot (L - K_S) = 2 + \frac{1}{2} L^2 >2. 
\]
We then get $H^1(S, -L)=0$ (Theorem \ref{t-Mumford-vanishing}), 
which implies $H^1(S, L)=0$ by Serre duality. 
Thus (1) and (2) hold, which implies (3). 
\end{proof}

\subsection{Linear systems}

In this subsection, we study base locus of divisors on geometrically integral regular K3-like surfaces (Theorem \ref{t-K3-Bs-reg}). 
The hardest part is to prove that nef and big divisors without fixed components is base point free (Proposition \ref{p-K3-big-key}). 
Many arguments are based on the proofs of the corresponding results over algebraically closed fields \cite[Chapter 2]{Huy16}.

\begin{prop}\label{p-curve-Kbpf}
We work over a field $\kappa$ of characteristic $p>0$. 
Let $C$ be a projective regular curve such that $\deg K_C \geq 0$. 
Then $|K_C|$ is base point free. 
\end{prop}

\begin{proof}
Replacing $\kappa$ by $H^0(C, \MO_C)$, 
the problem is reduced to the case when  $H^0(C, \MO_C)  = \kappa$. 
Set 
\[
g :=h^1(C, \MO_C) = h^0(C, K_C). 
\]
By the Riemann--Roch theorem and Serre duality, we obtain 
\[
-\chi(C, \MO_C) = \chi(C, K_C) = \deg K_C + \chi(C, \MO_C), 
\]
which implies 
\[
\deg K_C = 2g -2. 
\]
If $\deg K_C =0$, i.e., $g =h^1(C, \MO_C)=1$, then $H^0(C, K_C) \neq 0$, which implies $K_C \sim 0$. 

We may assume that $2g -2 =\deg K_C >0$, i.e., $g \geq 2$. 
Then $K_C$ is ample. 
By $h^0(C, K_C) = g >0$, 
there exists an effective Cartier divisor $D$ on $C$ with $K_C \sim D$. 

Fix a closed point $P \in C$. 
It suffices to show that $P \not\in |K_C|$. 
We may assume that $P \in \Supp\,D$, i.e., we can write 
\[
D = aP + D'
\]
for some $a \in \Z_{>0}$ and an effective Cartier divisor $D'$ with $P \not\in \Supp\,D'$. 
We have the following exact sequence: 
\[
H^0(C, \MO_C(K_C)) \to H^0(P, \MO_C(K_C)|_P) \to H^1(C, \MO_C(K_C-P)) \to H^1(C, \MO_{\cred C}(K_C)) \to 0. 
\]
If $\dim_{\kappa} H^1(C, K_C-P)=1$, then  $H^1(C, K_C-P) \xrightarrow{\simeq} H^1(C, K_C)$, 
and hence we are done. 
Therefore, we may assume that $\dim_{\kappa} H^1(C, K_C-P) \geq 2$, i.e.,  $\dim_{{\cred \kappa}} H^0(C, \MO_C(P)) \geq 2$. 
Then $P$ is linearly equivalent to an effective Cartier divisor $E$ with $E \neq P$. 
Therefore, we obtain $P \not\in E$. 
This implies that $|P|$ is base point free. 
Hence $P$ is not a base point of $|K_C|$. 
\end{proof}

\begin{prop}\label{p-K3-non-big}
We work over a field $\kappa$ of characteristic $p>0$. 
Let $S$ be a 
regular K3-like surface. 
Let $L$ be a Cartier divisor on $S$. 
If $L$ is nef and $L^2=0$, then $|L|$ is base point free. 
\end{prop}

\begin{proof}
Replacing $\kappa$ by $H^0(S, \MO_S)$, we may assume that $H^0(S, \MO_S) = \kappa$. 

We first treat the case when $L \equiv 0$. 
It suffices to show that $L \sim 0$, i.e.,  $H^0(S, L) \neq 0$. 
Otherwise, we have that $H^0(S, L)= H^0(S, -L)=0$, and hence $H^2(S, L)=0$ by Serre duality. 
Therefore, the Riemann--Roch theorem implies that 
\[
0 \geq -h^1(S, L) = \chi(S, L) =  \frac{1}{2} L^2 +2 = 2 >0,  
\]
which is absurd.

Let us handle the remaining case: $L \not\equiv 0$. 
By Serre duality, we obtain $H^2(S, L) =0$. 
Then the Riemann--Roch theorem implies that 
\[
h^0(S, L) \geq h^0(S, L) - h^1(S, L) = \chi(S, L) =\frac{1}{2} L^2 +2 = 2. 
\]
Take the decomposition $L=M+F$, where $M$ is the mobile part of $|L|$ 
and $F$ is the fixed part of $|L|$. 
Note that $M$ is nef. 
We have 
\[
0 = L^2 = L \cdot (M+F) = L \cdot M + L \cdot F. 
\]
It follows from $L \cdot M \geq 0$ and $L \cdot F \geq 0$ that
\[
L \cdot M = L \cdot F=0.
\]
By $M^2 \geq 0, F \cdot M \geq 0$, and $(M+F) \cdot M = L \cdot M =0$, we get 
\[
M^2 = M \cdot F=0.
\]
Therefore, we obtain 
\[
F^2  = (L-M) \cdot F = 0. 
\]
If $F \neq 0$, then 
\[
h^0(S, F) \geq h^0(S, F) -h^1(S, F) =\chi(S, F) = \frac{1}{2}F^2 +2 =2, 
\]
which is a contradiction. 
Therefore, we have $F=0$, i.e., the base locus ${\rm Bs}\,|L|$ is zero-dimensional (if non-empty). 
However, if ${\rm Bs}\,|L|$ is non-empty and zero-dimensional, 
then we again obtain a contradiction by $L^2 =0$. 
Thus $|L|$ is base point free. 
\end{proof}

\begin{lem}\label{l-K3-big-curve}
We work over a field $\kappa$ of characteristic $p>0$. 
Let $S$ be a 
regular K3-like surface. 
Let $C$ be a prime divisor on $S$ with $C^2 \geq 0$. 
Then the following hold. 
\begin{enumerate}
\item 
$\dim \Bs\,|C| \leq 0$. 
\item 
If $C$ is regular, then $|C|$ is base point free. 
\end{enumerate}
\end{lem}

\begin{proof}
If $C^2 =0$, then $|C|$ is base point free by Proposition \ref{p-K3-non-big}. 
Hence we may assume that $C^2 >0$. 
We then obtain $h^0(S, C) \geq 2$ by the Riemann--Roch theorem, which implies (1). 

Let us show (2). 
Consider the exact sequence: 
\[
H^0(S, \MO_S(C)) \to H^0(C, \MO_S(C)|_C) \simeq H^0(C, K_C) \to H^1(S, \MO_S)=0. 
\]
As $|K_C|$ is base point free (Proposition \ref{p-curve-Kbpf}), 
$|C|$ is base point free. 
\end{proof}

\begin{lem}\label{l-conic-nef}
We work over a field $\kappa$ of characteristic $p>0$. 
Let $C$ be a projective Gorenstein curve such that $\omega_C^{-1}$ is ample. 
Let $L$ be a nef line bundle on $C$. 
Then $|L|$ is base point free. 
\end{lem}

\begin{proof}
{\cred Replacing $\kappa$ by $H^0(C, \MO_C)$, we may assume that $H^0(C, \MO_C) = \kappa$.} 
By \cite[Lemma 10.6]{Kol13}, 
there exists an ample invertible sheaf $H$ on $C$ such that $\Pic\,C \simeq \Z[H]$, 
where $[H]$ denotes the isomorphism class of $H$. 
Therefore, it suffices to show that $|H|$ is base point free. 
Again by \cite[Lemma 10.6]{Kol13}, one of the following holds. 
\begin{enumerate}
\item $C \simeq \P^1_{\kappa}$. 
\item $C$ is a conic on $\P^2_{\kappa}$ and $H \simeq \omega_C^{-1}$. 
\end{enumerate}
Hence we may assume that (2) holds. 
For the embedding $C \subset \P^2_{\kappa}$, we obtain 
\[
\MO_{\P^2_{\kappa}}(-1)|_C \simeq (\omega_{\P^2_{\kappa}} \otimes \MO_{\P^2}(C))|_C \simeq \omega_C. 
\]
Therefore, $H \simeq \omega_C^{-1} \simeq \MO_{\P^2_k}(1)|_C$, which implies that $|H|$ is very ample. 
\end{proof}

\begin{lem}\label{l-1-conne}
We work over a field $\kappa$ of characteristic $p>0$. 
Let $S$ be a regular projective surface. 
Let $D$ be an effective $\Z$-divisor. 
{\cred Then the following hold.} 
\begin{enumerate}
\item {\cred If $D$ is nef and big, then 
$D$ is $1$-connected, i.e., 
$C_1 \cdot C_2 >0$ holds for all nonzero effective $\Z$-divisors $C_1, C_2$
with $D = C_1 +C_2$.} 
\item If $D$ is $1$-connected, then $H^0(D, \MO_D)$ is a field. 
\end{enumerate}
\end{lem}

\begin{proof}
{\cred 
Let us show (1). 
By the same argument as in \cite[Ch. V, Theorem 1.9]{Har77}, 
the Hodge index theorem holds even over an imperfect field, i.e., 
the signature of the intersection product is $(1, \rho(S)-1)$. 
Then the argument as in \cite[Remark 1.7]{Huy16} works without any changes. 
Thus (1) holds.}

{\cred Let us show (2).} 
Take a maximal effective $\Z$-divisor $D'$ 
such that $0 \leq D' \leq D$ and $H^0(D', \MO_{D'})$ is a field (such $D'$ exists, since a prime divisor satisfies this condition). 
Assume $D \neq D'$. Then we have $D' \cdot (D -D')>0$. 
Then $D' \cdot C >0$ for some prime divisor $C \subset \Supp\,(D -D')$. 
We have an exact sequence 
\[
0 \to \MO_S(-D')|_C \to \MO_{D'+C} \to \MO_{D'} \to 0, 
\]
which induces the following injection 
\[
H^0(D'+C, \MO_{D'+C}) \hookrightarrow H^0(D', \MO_{D'}). 
\]
Then $H^0(D'+C, \MO_{D'+C})$ is an intermediate ring between $\kappa$ and $H^0(D', \MO_{D'})$. 
Since the field extension $\kappa \subset H^0(D', \MO_{D'})$ is of finite degree, 
$H^0(D' {\cred +C}, \MO_{D'+C})$ is a field. 
This contradicts the maximality of $D'$. 
{\cred Thus (2) holds.} 
\end{proof}

\begin{prop}\label{p-K3-big-key}
We work over a field $\kappa$ of characteristic $p>0$. 
Let $S$ be a geometrically integral regular K3-like surface. 
Let $M$ be a nef and big Cartier divisor on $S$. 
If $\dim {\rm Bs}|M| \leq 0$, then $|M|$ is base point free. 
\end{prop}


\begin{proof}
Taking the base change to the separable closure of $\kappa$, 
we may assume that $\kappa$ is separably closed. 
We prove the assertion by induction on $M^2$. 
We may start with the case when $M^2=1$ as the base of this induction. 
By $M^2 \in 2\Z$, there is nothing to show. 
In what follows, we assume that $|{\cred \widetilde M}|$ is base point free 
if ${\cred \widetilde M}$ is a nef and big Cartier divisor on $S$ such that  
${\cred \widetilde M}^2 < M^2$ and $\dim {\rm Bs}\,|{\cred \widetilde M}| \leq 0$.

Suppose that $|M|$ is not base point free. 
Let us derive a contradiction. 
Note that ${\rm Bs}\,|M|$ (set-theoretically) consists of finitely many closed points $Q:=Q_1, Q_2, ..., Q_r$. 
Set $d(Q):=[k(Q):{\cred \kappa}]$.

\begin{step}\label{s1-K3-big-key}
There exists an index $1 \leq i \leq r$ such that 
the generic member $M'$ of $|M|$ is not regular at $Q_i$. 
\end{step}

\begin{proof}[Proof of Step \ref{s1-K3-big-key}]
Suppose the contrary, i.e., 
suppose that $M'$ is regular at each of $Q_1, ..., Q_r$. 
By Theorem \ref{t-Bertini-generic}, 
$M'$ is regular outside $Q_1, ..., Q_r$. 
In particular, $M'$ is a regular curve on a regular K3-like surface. 
Then $|M \otimes_{\kappa} \kappa'| = |M'|$ is base point free (Lemma \ref{l-K3-big-curve}(2)), where $\kappa' := K(\P(H^0(S, M)))$. 
Hence also $|M|$ is base point free, which is absurd. 
This completes the proof of Step \ref{s1-K3-big-key}. 
\end{proof}

After possibly permuting the indices, we may assume that $i=1$. 
In particular, any member $N \in |M|$ is not regular at $Q$. 
Set
\[
\mu := \min_{D \in |M|} \{ \mult_Q\,D\}.
\]
We have $\mu \geq 2$. 
Let $f: T \to S$ be the blowup at $Q$. Set $E:= \Ex(f)$ and $N := f^*M -2E$.

\begin{step}\label{s2-K3-big-key}
It holds that $H^1(T, -N) \neq 0$. 
\end{step}

\begin{proof}[Proof of Step \ref{s2-K3-big-key}]
By the exact sequence
\[
H^0(T, f^*M) \xrightarrow{0} H^0(E, f^*M) \to H^1(T, f^*M-E)
\]
and $H^0(E, f^*M) \neq 0$, we obtain 
$H^1(T, f^*M-E) \neq 0$. 
By Serre duality and $K_T \sim f^*K_S+E\sim E$, we get 
$H^1(T, -(f^*M-2E)) \neq 0$. 
This completes the proof of Step \ref{s2-K3-big-key}. 
\end{proof}

\begin{step}\label{s3-K3-big-key}
$f^*M-\mu' E$ is nef for any $0 \leq \mu' \leq \mu$. 
In particular, $N$ is nef. 
\end{step}

\begin{proof}[Proof of Step \ref{s3-K3-big-key}]
Since $f^*M$ is nef, it suffices to show that $f^*M -\mu E$ is nef. 
Pick two general members $M_1, M_2 \in |M|$. 
We then obtain 
\[
\mult_{Q} M_1 = \mult_Q M_2 = \mu, \qquad \dim (M_1 \cap M_2) =0. 
\]
Therefore, 
$f^*M_1 - \mu E$ and $f^*M_2 - \mu E$ share no irreducible components. 
In particular,  $f^*M - \mu E$ is nef. 
This completes the proof of Step \ref{s3-K3-big-key}. 
\end{proof}

\begin{step}\label{s4-K3-big-key}
The following hold. 
\begin{enumerate}
\item $N^2=0$. 
\item $\mu=2$. 
\item $r=1$, i.e., the set-theoretic equation ${\rm Bs}|M| = \{Q\}$ holds. 
\item $M^2 = 4d(Q)$. 
\item For any $M_1 \in |M|$, we have $\mult_Q M_1 =2$. 
\item $|N|$ is base point free. 
\end{enumerate}
\end{step}

\begin{proof}[Proof of Step \ref{s4-K3-big-key}]
Note that $N$ is nef by Step \ref{s3-K3-big-key}. 

Let us show (1). 
If $N^2>0$, then $N$ is nef and big. 
Then also $f_*N$ is nef and big, and hence 
$h^0(T, N) \geq h^0(S, f_*N) \geq 2$. 
Since $T$ is geometrically integral and $H^1(T, \MO_T)=0$, 
we get $H^1(T, -N)=0$ (Theorem \ref{t-Mumford-vanishing}), 
which contradicts Step \ref{s2-K3-big-key}. 
Thus (1) holds. 

Let us show (2). 
Suppose the contrary, i.e., $\mu \geq 3$. 
Note that $f^*M -\mu'E$ is nef for any $0 \leq \mu' \leq 3$. 
By $(f^*M)^2>0$ and $(f^*M -3E)^2 \geq 0$, we have that $N^2 = (f^*M -2E)^2 >0$, 
which contradicts (1). 
Thus (2) holds. 

Let us show (3). 
Suppose $r \geq 2$. 
Pick two general members $M_1, M_2 \in |M|$. 
Then we have that 
\[
\mult_{Q} M_1 = \mult_Q M_2 = 2, \qquad \dim (M_1 \cap M_2) =0. 
\]
Then  
\[
(f^*M_1 -2E) \cap 
(f^*M_2 -2E)
\]
is zero-dimensional. 
By $r \geq 2$, 
\[
(f^*M_1 -2E) \cap 
(f^*M_2 -2E) \neq \emptyset, 
\]
i.e., this set contains the point $f^{-1}(Q_2)$. 
Hence, we obtain $N^2 =(f^*M_1 -2E) \cdot (f^*M_2 -2E) >0$, which contradicts (1). 
Thus (3) holds. 

Let us show (4). 
We have that $0 = N^2 = (f^*M -2E)^2 = M^2 +4E^2 = M^2 -4d(Q)$, 
and hence (4) holds. 

Let us show (5). 
Suppose the contrary, i.e., $\mult_Q M_1 \geq 3$ for some $M_1 \in |M|$. 
Pick a general member $M_2 \in |M|$. 
We have that $\mult_Q M_2 =2$, and $\dim (M_1 \cap M_2) \leq 0$. 
It holds that $N^2 = (f^*M_1 -2E) \cdot (f^*M_2 -2E) >0$, which contradicts (1). 
Thus, (5) holds. 

The assertion (6) follows from {\cred the proof of (3)}. 
This completes the proof of Step \ref{s4-K3-big-key}. 
\end{proof}

\begin{step}\label{s5-K3-big-key}
A general member of $|N|$ is a disjoint union of irreducible effective $\Z$-divisors. 
\end{step}

\begin{proof}[Proof of Step \ref{s5-K3-big-key}]
Let 
\[
\pi :=\varphi_{|N|} :T \to B \subset \mathbb P_{\cred \kappa}^{h^0(T, N) -1}. 
\]
be the morphism induced by $|N|$, where $B$ denotes its image. 
Since $|N|$ is base point free with $N \not\equiv 0$ (Step \ref{s2-K3-big-key}), 
$B$ is a (possibly non-regular) curve. 
Let 
\[
\pi: T \xrightarrow{\pi'} B' \xrightarrow{\theta } B  
\]
be the Stein factorisation. 
Note that each of $B$ and $B'$ has infinitely many $\kappa$-rational points. 
Pick a general closed point $b \in B$. 
Then $\theta^{-1}(b)$ consists of finitely many general closed points of $B'$. 
Since general fibres of $T \to B'$ are (geometrically) irreducible, 
$\pi^*(b)$ is a disjoint union of irreducible divisors. 
This implies that a general member of $|N|$ is a disjoint union of irreducible divisors. 
This completes the proof of Step \ref{s5-K3-big-key}. 
\end{proof}

\begin{step}\label{s6-K3-big-key}
Fix a non-empty open subset $T' \subset T$ on which $f$ is an isomorphism. 
Then there exist a ${\cred \kappa}$-rational point $R \in T'$ and $N' \in |N|$ such that 
$R \in N'$ and $N'$ is a disjoint union of irreducible divisors. 
Set $M' :=f_*N'$. In particular, $M'$ passes through the  $\kappa$-rational point $R_S :=f(R) \in S$. 
\end{step}

\begin{proof}[Proof of Step \ref{s6-K3-big-key}]
Fix two  general members $N_1, N_2 \in |N|$, 
so that  each of $N_1$ and $N_2$ is a disjoint union of irreducible effective $\Z$-divisors 
(Step \ref{s6-K3-big-key}). 
Let $V \subset H^0(T, N)$ be the two-dimensional ${\cred \kappa}$-vector subspace corresponding to the pencil generated by $N_1, N_2$. 
For the pencil $|V|$ corresponding to $V$, 
a general member $N' \in |V|$ is a disjoint union of irreducible divisors. 
Therefore, except for finitely many members $\widetilde N_1, ..., \widetilde N_t$, 
any member of $|V|$ is a disjoint union of irreducible divisors. 
Pick a $\kappa$-rational point $R \in T'$ such that $R \not\in \bigcup_{j=1}^t \widetilde N_j$. 
Consider an exact sequence: 
\[
0 \to H^0(T, N \otimes \m_R) \to H^0(T, N) \to H^0(R, N|_R). 
\]
By $\dim V =2$ and  $\dim H^0(R, N|_R) =1$, 
we have that $H^0(T, N \otimes \m_R) \cap V \neq \{0\}$. 
Then there is a member $N' \in |N|$ such that $R \in N'$. 
By the choice of $R$, we have that $N' \neq \widetilde N_1, ..., N' \neq \widetilde N_t$. 
Then $N'$ is a  disjoint union of irreducible divisors, as required. 
This completes the proof of Step \ref{s6-K3-big-key}. 
\end{proof}

\begin{step}\label{s7-K3-big-key}
$N'$ is not $1$-connected, i.e., there exists a decomposition 
\[
N' = C_{1, T} +C_{2, T}
\]
for some nonzero effective $\Z$-divisors $C_{1, T}, C_{2, T}$ 
such that $C_{1, T} \cdot C_{2, T} \leq 0$. 
Set $C_i := f_*C_{i, T}$. 
\end{step}

\begin{proof}[Proof of Step \ref{s7-K3-big-key}]
If $N'$ is $1$-connected, then $H^0(N', \MO_{N'})$ is a field (Lemma \ref{l-1-conne}). 
Since $N'$ contains a ${\cred \kappa}$-rational point $R$ (Step \ref{s6-K3-big-key}), it holds that $H^0(N', \MO_{N'}) =\kappa$. 
Then we have an exact sequence 
\[
H^0(T, \MO_T) \xrightarrow{\simeq} H^0(N', \MO_{N'}) \to H^1(T, -N') \to H^1(T, \MO_T)=0. 
\]
Thereofore, we obtain 
\[
0 = H^1(T, -N') \simeq H^1(T, -N), 
\]
which contradicts Step \ref{s2-K3-big-key}. 
This completes the proof of Step \ref{s7-K3-big-key}. 
\end{proof}

\begin{step}\label{s8-K3-big-key}
The following hold. 
\begin{enumerate}
\renewcommand{\labelenumi}{(\roman{enumi})}
\item $C_1 \neq 0$ and $C_2 \neq 0$. 
\item For each $i \in \{1, 2\}$, we have $Q \in C_i$. 
\item For each $i \in \{1, 2\}$, we have $\mult_{Q} C_i=1$ and $C_{i, T} = f^*C_i -E$. 
\item $C_{1, T} \cdot C_{2, T}  = C_{1} \cdot C_{2} -  d(Q)$. 
\end{enumerate}
\end{step}

\begin{proof}[Proof of Step \ref{s8-K3-big-key}]
The assertion (i) follows from the fact that $N'$ does not contain $E$. 

Let us show (ii). 
If one of $C_1$ and $C_2$, say $C_1$, does not contain $Q$, then 
$C_{1, T} \cdot E=0$, which implies 
\[
C_{1, T} \cdot C_{2, T} = C_{1, T} \cdot f^*f_*C_{2, T} = C_1 \cdot C_2 >0.  
\]
{\cred Here $C_1 \cdot C_2 >0$ follows from $M \sim f_*N' =C_1+C_2$ 
and $M$ is $1$-connected (Lemma \ref{l-1-conne}(1)).}  
This contradicts Step \ref{s7-K3-big-key}. Thus (ii) holds. 

Let us show (iii). 
It follows from (ii) that $\mult_Q C_1 \geq 1$ and $\mult_Q C_2 \geq 1$. 
By $\mult_Q (M') = 2$ and $M' = C_1+C_2$, 
we get $\mult_Q C_1 = \mult_Q C_2 = 1$. 
Since $C_{i, T}$ does not contain $E$, we have that $C_{i, T} = f^*C_i -E$. 
Thus (iii) holds. 

The assertion (iv) follows from 
\[
C_{1, T} \cdot C_{2, T} = (f^*C_{1} - E) \cdot (f^*C_{2} - E) 
=C_{1} \cdot C_{2} -  d(Q). 
\]
This completes the proof of Step \ref{s8-K3-big-key}. 
\end{proof}

\begin{step}\label{s9-K3-big-key}
Each of $C_1$ and $C_2$ is a prime divisor. 
Furthermore, one of the following holds. 
\begin{enumerate}
\item[(I)] 
$M'$ and $N'$ are irreducible. Furthermore, $M' =2C$ and $N' =2C_T$ 
for $C:=C_1=C_2$ and $C_T := C_{1, T} = C_{2, T}$. 
\item[(II)] 
$C_1 \neq C_2$, $C_{1, T} \cdot C_{2, T}=0$, and $C_1 \cdot C_2 = d(Q)$. 
Moreover, $C_{1, T} \cap C_{2, T} =\emptyset$ and 
$M'=C_1+C_2$ is simple normal crossing at $Q$. 
\end{enumerate}
\end{step}

\begin{proof}[Proof of Step \ref{s9-K3-big-key}]
Recall that $N'$ is a disjoint union of irreducible divisors and $M'=f_*N'$ is connected. 
We can write 
\[
N' = N'_1 + \cdots {\cred +} N'_s,
\]
where each $N'_i$ is an irreducible effective $\Z$-divisor and $N'_{i_1} \cap N'_{i_2} = \emptyset$ for any pair $(i_1, i_2)$ with $i_1 \neq i_2$. 
Note that $E \not\subset {\cred N'}$, and hence $M'_i := f_*N'_i \neq 0$ for any $i$. 

We now show that each $M'_i$ passes through $Q$. 
Suppose the contrary, i.e., $N'_i \cap E =\emptyset$. 
Then $f:T \to S$ is an isomorphism around $N'_i$, 
and hence 
$M'_i$ is disjoint from $\bigcup_{j \neq i} M'_j$. 
This is a contradiction. 

By $\mult_Q M' = 2$, we have $s \leq 2$. 
We first treat the case when $s=1$. 
In this case, $N'$ and $M'$ are irreducible. 
Hence (I) holds. 

We may assume that $s=2$. 
By $C_1 \cdot C_2 \geq d(Q)$ and $C_{1, T} \cdot C_{2, T} \leq 0$, 
we obtain $C_1 \cdot C_2 =d(Q)$ and $C_{1, T} \cdot C_{2, T} = 0$. 
In this case, we obtain (II). 
This completes the proof of Step \ref{s9-K3-big-key}. 
\end{proof}

\begin{step}\label{s10-K3-big-key}
(I) does not hold. 
\end{step}

\begin{proof}[Proof of Step \ref{s10-K3-big-key}]
Suppose that (I) holds. 
By $M^2 >0$, we have that 
\[
M^2 =M'^2 = (2C)^2 =4C^2 > C^2 >0.  
\]
It holds that $\Bs\,|C| \leq 0$ {\cred (Lemma \ref{l-K3-big-curve}(1))}, 
and hence the induction hypothesis can be applied for $|C|$, so that $|C|$ is base point free. 
Then also $|M| = |2C|$ is base point free, which is a contradiction. 
This completes the proof of Step \ref{s10-K3-big-key}. 
\end{proof}

\begin{step}\label{s11-K3-big-key}
(II) does not hold. 
\end{step}

\begin{proof}[Proof of Step \ref{s11-K3-big-key}]
Suppose that (II) holds. 
By symmetry, we may assume that $C_1^2 \geq C_2^2$. 
By $C_1 \cdot C_2 =d(Q) >0$, we get 
\[
M^2 = (C_1+C_2)^2 = C_1^2 + 2 C_1 \cdot C_2 +C_2^2 = C_1^2 + C_2^2 +2d(Q) > C_1^2 + C_2^2. 
\]
By $M^2 =4d(Q)$, 
it holds that  
\[
C_1^2 + C_2^2 = 2d(Q) >0. 
\]
In particular, we get $C_1^2 >0$. 

If $C_2^2 \geq 0$, we obtain 
\[
M^2 >C_1^2 +C_2^2 \geq C_i^2 
\]
for any $i \in \{1, 2\}$. 
Hence $|C_1|$ and $|C_2|$ are base point free by the induction hypothesis, 
{\cred Proposition \ref{p-K3-non-big}, and Lemma \ref{l-K3-big-curve}(1)}. 
Then also $|M| = |C_1+C_2|$ is base point free, which is a contradiction. 

We may assume that $C_2^2 <0$. 
In this case, $|\MO_S(M)|_{C_2}|$ is base point free (Lemma \ref{l-conic-nef}). 
As $C_1$ is nef and big, we have the following exact sequence (Theorem \ref{t-K3-vanishing}): 
\[
H^0(S, \MO_S(M)) \to H^0(C_2, \MO_S(M)|_{C_2}) \to H^1(S, \MO_S(C_1)) =0. 
\]
Therefore, ${\rm Bs}|M|$ is disjoint from $C_2$. 
This contradicts $Q \in \Bs\,|M|$ and $Q \in C_2$. 
This completes the proof of Step \ref{s11-K3-big-key}. 
\end{proof}
Step \ref{s9-K3-big-key}, 
Step \ref{s10-K3-big-key}, and Step \ref{s11-K3-big-key}
complete the proof of Proposition \ref{p-K3-big-key}. 
\end{proof}

\begin{cor}\label{c-K3-curve-bpf}
We work over a field $\kappa$ of characteristic $p>0$. 
Let $S$ be a geometrically integral regular K3-like surface. 
Let $C$ be a prime divisor on $S$. 
If $C^2 \geq 0$, then $|C|$ is base point free. 
\end{cor}

\begin{proof}
If $C^2=0$, then apply Proposition \ref{p-K3-non-big}. 
Otherwise, $C$ is nef and big. 
In this case, $|C|$ is base point free by Lemma \ref{l-K3-big-curve}(1) and Proposition \ref{p-K3-big-key}. 
\end{proof}

\begin{cor}\label{c-K3-fixed-part}
We work over a field $\kappa$ of characteristic $p>0$. 
Let $S$ be a geometrically integral regular K3-like surface. 
Let $L$ be a Cartier divisor on $S$. 
Then ${\rm Bs}\,|L| = F$ for the fixed part $F$ of $L$. 
\end{cor}

Note that the scheme-theoretic equation $\Bs\,|L| = F$ holds, i.e., 
the equation $I_{\Bs\,|L|} = I_F$ holds for the corresponding ideal sheaves 
$I_{\Bs\,|L|}, I_F$ on $S$. 

\begin{proof}
Let $M$ be the mobile part of $|L|$. 
If $M=0$, then there is nothing to show. 
If $M \neq 0$, then $M$ is a nef Cartier divisor with $\dim \Bs |M| \leq 0$. 
Then $M$ is base point free by Proposition \ref{p-K3-non-big} and Proposition \ref{p-K3-big-key}. 
\end{proof}

\subsection{Nef divisors with base points}

\begin{lem}\label{l-H^0-image}
We work over a field $\kappa$. 
Let $X$ be a projective variety. 
For a  Cartier divisor $L$ such that $|L|$ is base point free, 
let $\varphi_{|L|} : X \to \P^{h^0(X, L)-1}$ be the induced morphism. 
Assume that we have morphisms whose composition is $\varphi_{|L|}$: 
\[
\varphi_{|L|} : X \xrightarrow{\psi} Y \xrightarrow{\theta} \P^{N}
\]
with $N := h^0(X, L)-1$ such that $Y$ is a projective variety and $\psi : X \to Y$ is surjective 
(e.g., $Y$ is the image of $\varphi_{|L|}$ or the Stein factorisation of $\varphi_{|L|}$). 
Then the following induced maps are isomorphisms: 
\[
\varphi_{|L|}^* : 
H^0(\P^N, \MO_{\P^N}(1))
 \xrightarrow{\theta^*, \simeq} H^0(Y, \theta^*\MO_{\P^N}(1)) 
 \xrightarrow{\psi^*, \simeq} H^0(X, L). 
\]
\end{lem}

\begin{proof}
By construction, $\varphi_{|L|}^* = \psi^* \circ \theta^*$ is an isomorphism. 
In particular, $\psi^*$ is surjective. 
Since $\psi : X \to Y$ is surjective, $\MO_Y \to \psi_*\MO_X$ is injective, and hence $\psi^*$ is injective. 
Therefore, $\psi^*$ is an isomorphism. Thus the remaining map $\theta^*$ is an isomorphism. 
\end{proof}

\begin{lem}\label{l-K3-to-curve}
We work over a field $\kappa$ of characteristic $p>0$. 
Let $S$ be a geometrically integral regular K3-like surface. 
Let $\pi:S \to B$ be a morphism to a regular projective curve $B$ with $\pi_*\MO_S =\MO_B$. 
Then the following hold. 
\begin{enumerate}
\item $B$ is smooth and $-K_B$ is ample. 
In particular, $B \simeq \mathbb P^1_{\kappa}$ if and only if $B$ has a $\kappa$-rational point. 
\item 
For a nef Cartier divisor $N_B$ on $B$, it hold that 
\[
 \dim_{\kappa} H^0(S, \pi^*N_B) = \deg_{\kappa} N_B +1. 
\]
\end{enumerate}
\end{lem}

\begin{proof}
We have the injection $H^1(B, \MO_B) \hookrightarrow H^1(S, \MO_S)=0$ 
induced by the corresponding Leray spectral sequence. 
It follows from $\deg_{\kappa} K_B = 2h^1(B, \MO_B) -2 =-2$ that $-K_B$ is ample. 
Then $B$ is a geometrically integral conic on $\P^2_{\kappa}$ \cite[Lemma 10.6]{Kol13}, and hence $B$ is smooth. 
Since $B \times_{\kappa} \overline{\kappa} \simeq \P^1_{\overline{\kappa}}$, 
it holds by Ch\^atelet's theorem \cite[Theorem 5.1.3]{GS17} that $B \simeq \mathbb P^1_{\kappa}$ if and only if $B$ has a $\kappa$-rational point.
Thus (1) holds. 
By $\pi_*\MO_S = \MO_B$ and the Riemann--Roch theorem for regular projective curves, we get 
\[
 \dim_{\kappa} H^0(S, \pi^*N_B) = \dim_{\kappa} H^0(B, N_B) = \chi (B, N_B) = \deg_{\cred \kappa} N_B +1, 
\]
where $H^1(B, N_B) =0$ follows from Serre duality. 
Thus (2) holds.  
\end{proof}

\begin{lem}\label{l-K3-Bs-reg}
We work over a field $\kappa$ of characteristic $p>0$. 
Let $S$ be a geometrically integral regular K3-like surface. 
For an effective Cartier divisor $F$ and nef and big Cartier divisors $L$ and $M$, 
assume that 
\[
L = M + F \qquad \text{and} \qquad h^0(S, L) = h^0(S, M). 
\]
Then $F=0$. 
\end{lem}

\begin{proof}
Suppose $F \neq 0$. 
Let us derive a contradiction. 
We can write $F = \sum_{i=1}^r F_i$, where $r \geq 1$ and each $F_i$ is a prime divisor 
($F_i =F_j$ might hold even if $i \neq j$). 
Note that $h^0(S, L) = h^0(S, M)$ implies that the following induced injection is an isomorphism: 
\[
H^0(S, M) \xrightarrow{\simeq} 
 H^0(S, M + F) = H^0(S, L). 
\]
We have $H^i(S, L) = H^i(S, M)=0$ for $i>0$ (Theorem \ref{t-K3-vanishing}), 
and hence $\chi(S, L) = h^0(S, L) = h^0(S, M) = \chi(S, M)$. 
By the  Riemann--Roch theorem, we conclude $L^2 = M^2$. 
We then get $M^2 = L^2 = (M+F)^2 = M^2 + 2 M \cdot F + F^2$, which implies 
\[
M \cdot F + L \cdot F = 2M \cdot F + F^2=0. 
\]
By $M \cdot F \geq 0$ and $L \cdot F \geq 0$, 
we have that $M \cdot F = L \cdot F=0$. 
Since $M$ and $L$ are nef, 
we get $M \cdot F_i = L \cdot F_i =0$ and hence $F \cdot F_i =0$ for all $1 \leq i \leq r$. 
Therefore, $F$ is nef with $F^2 = \sum_{i=1}^r F \cdot F_i =0$. 
Then $|F|$ is base point free by Proposition \ref{p-K3-non-big}. 
This contradicts the fact that $F$ is contained in the fixed part of $|L|$. 
\end{proof}

\begin{thm}\label{t-K3-Bs-reg}
We work over a field $\kappa$ of characteristic $p>0$. 
Let $S$ be a geometrically integral regular K3-like surface. 
Let $L$ be a nef and big Cartier divisor on $S$. 
Set $g :=\frac{1}{2}L^2 +1$. 
Let 
\[
L = M+F
\]
be the decomposition into the mobile part $M$ of $|L|$ and the fixed part $F$ of $|L|$. 
Assume that $|L|$ is not base point free. 
Then the following hold. 
\begin{enumerate}
\item $h^0(S, L) = \frac{1}{2}L^2 +2 = g+1$ and $g \geq 2$. 
\item $\Bs\,|L| = F$ as closed subschemes of $S$. 
\item $|M|$ is base point free and $M^2 =0$. 
\end{enumerate}
For $B := \varphi_{|M|}(S)$, we have the induced morphisms: 
\[
\varphi_{|M|} : S \xrightarrow{\pi} B \hookrightarrow \P^g. 
\] 
\begin{enumerate}
\setcounter{enumi}{3}
\item $\Delta(B, \MO_{\P^g}(1)|_B)=0$, $\pi_*\MO_S = \MO_B$, $\deg (\MO_{\P^g}(1)|_B) =g$, 
and $B$ is a smooth projective curve with $H^1(B, \MO_B)=0$. 
\item $F$ is a prime divisor which is a section of $\pi : S \to B$, i.e., $\pi|_F : F \xrightarrow{\simeq} B$. 
Moreover, $H^0(F, \MO_F)=\kappa$, $F^2 = -2$, and $L \cdot F = g-2$. 
\item 
If $L \cdot F=0$, then $L^2=2$ and $g=2$. 
\item 
If $L \cdot F>0$, then $L \cdot C >0$  
for any curve $C {\cred \subset} S$ with $C \cap F \neq \emptyset$. 
\end{enumerate}
Furthermore, if $S$ has a $\kappa$-rational point, then $F \simeq B \simeq \P^1_{\kappa}$. 
\end{thm}

\begin{proof}
Taking the base change to the separable closure, we may assume that $\kappa$ is separably closed. 
In particular, $S$ has a $\kappa$-rational point and it suffices to show (1)--(7).

The assertions (1) and (2) hold by Theorem \ref{t-K3-vanishing} and Corollary \ref{c-K3-fixed-part}, respectively. 
By $F \neq 0$, we get $M^2=0$ (Lemma \ref{l-K3-Bs-reg}). 
Since $M$ is nef, $|M|$ is base point free (Proposition \ref{p-K3-non-big}). 
Thus (3) holds.

Let us show (4). 
Let 
\[
\pi : S \to B' \to B
\]
be the Stein factorisation of $\pi : S \to B$. 
{\cred By Remark \ref{r-Delta-geom-int},} we have 
\[
0 \leq \Delta(B, \MO_{\P^g}(1)|_B) = \dim B + \deg(\MO_{\P^g}(1)|_B) -h^0(B, \MO_{\P^g}(1)|_B). 
\]
For $\ell := \deg (B' \to B)$, it follows from $B' \simeq \P^1_{\kappa}$ and 
\[
H^0(S, L) \simeq H^0(S, M) \simeq H^0(B', \MO_{\P^g}(1)|_{B'}) \simeq H^0(B,  \MO_{\P^g}(1)|_{B})
\qquad 
(\text{Lemma \ref{l-H^0-image}})
\]
that 
\begin{itemize}
\item $h^0(B, \MO_{\P^g}(1)|_B) = g+1$ and 
\item $\deg(\MO_{\P^g}(1)|_B) = \frac{1}{\ell} \deg(\MO_{\P^g}(1)|_{B'}) = \frac{1}{\ell} (h^0(B', \MO_{\P^g}(1)|_{B'})-1) = \frac{g}{\ell}$. 
\end{itemize}
Hence we obtain 
\[
0 \leq \Delta(B, \MO_{\P^g}(1)|_B) = 
{\cred 1 + \frac{g}{\ell} -(g+1) = \frac{g}{\ell}(1-\ell).}
\]
Therefore, $\Delta(B, \MO_{\P^g}(1)|_B) =0$ and $\ell =1$. 
Since $B$ has a $\kappa$-rational point, we get $B \simeq \P^1_{\kappa}$. 
Then $B' \to B$ is a finite birational morphism to a normal variety, and hence an isomorphism. 
Hence $\deg(\MO_{\P^g}(1)|_B) = g$. 
Thus (4) holds.

Let us show (5). 
We first prove that 
\begin{enumerate}
\item[(5a)] $F$ is a prime divisor, and 
\item[(5b)] $\pi(F) = B$.
\end{enumerate}
Since $L = M+F$ is big, there exists a prime divisor $C \subset \Supp\,F$ such that $\pi(C) =B$. 
We can write $F = C+F'$ for some effective $\Z$-divisor $F'$. 
Suppose $F' \neq 0$. 
It suffices to derive a contradiction. 
By Lemma \ref{l-K3-Bs-reg}, it is enough to prove that $M+C$ is nef and big. 
For a $\kappa$-rational point $Q \in B = \P^1_{\kappa}$ 
and $\kappa_C:=H^0(C, \MO_C)$, 
we have $M \sim g \pi^*Q$ and  
\[
\deg_{\kappa_C} (M+C)|_C  \geq \deg_{\kappa_C} g (\pi^*Q|_C) + (-2) \geq g-2 \geq 0.  
\]
Hence $M+C$ is nef. {\cred Since also $M$ is nef, 
$M+a C$ is nef for every $0 \leq a \leq 1$. 
For some $0< \epsilon \ll 1$, 
it follows from $M^2 =0$ and $M \cdot C >0$ that 
\[
(M+\epsilon C)^2 = M^2 +2\epsilon M \cdot C + \epsilon^2 C^2 = 
\epsilon ( 2M \cdot C +\epsilon C^2) >0. 
\]
Hence $M+\epsilon C$ is big. Therefore, $M+C$ is nef and big.} 
Thus (5a) and (5b) hold.

Set $\delta := \deg(\pi|_F :  F \to B) =[K(F):K(B)]$, $\kappa_F :=H^0(F, \MO_F)$, 
and $d_F := [\kappa_F :\kappa]$. 
By the following commutative diagram 
\[
\begin{CD}
F @>>> \P^1_{\kappa} =B\\
@VVV @VVV\\
\Spec\,\kappa_F @>>> \Spec\,\kappa, 
\end{CD}
\]
we get the factorisation: 
\[
\pi|_F : F \xrightarrow{\alpha} \P^1_{\kappa_F} \xrightarrow{\beta} \P^1_{\kappa} =B. 
\]
We then obtain 
\[
\delta = \deg(\pi|_F : F \to \P^1_{\kappa}) = \deg(\alpha) \cdot \deg(\beta) = \deg(\alpha) [\kappa_F : \kappa] 
\in d_F \Z_{>0}. 
\]
It holds that 
\[
2g - 2 = L^2 = (M+F)^2 = M^2 + 2M \cdot F + F^2 = 0 + 2\delta g +F^2
\geq 2d_Fg -2d_F, 
\]
where the last inequality follows from $\delta \in d_F \Z_{>0}$ and 
\[
F^2 = (K_S+F) \cdot F = \deg_{\kappa} \omega_F = 
2\dim_{\kappa} H^1(F, \MO_F) - 2\dim_{\kappa} H^0(F, \MO_F) \geq  -2 d_F. 
\]
We then obtain $d_F=1$ and the above inequalities are equalities: 
$\delta =1$ and $F^2 = -2$. 
Since $\pi|_F : F \to B = \P^1_{\kappa}$ is birational, this is an isomorphism. 
Finally, we get $L \cdot F = (M+F) \cdot F = g-2$. 
Thus (5) holds.

Let us show (6). 
By $g -2 =L \cdot F =0$, we have $g=2$. Hence $L^2 = 2g -2 = 2$. 
Thus (6) holds.

Let us show (7). 
Assume $L \cdot F >0$. 
Fix a curve $C$ on $S$ with $C \cap F \neq \emptyset$. 
It suffices to show $L \cdot C>0$. 
If $C = F$, then $L \cdot C = L \cdot F >0$ by our assumption. 
If $C \neq F$, then we get $F \cdot C>0$, and hence 
\[
L \cdot C=(M+F) \cdot C = M \cdot C + F \cdot C \geq F \cdot C >0. 
\]
Thus (7) holds. 
\end{proof}


\begin{thm}\label{t-K3-Bs-cano}
We work over a field $\kappa$ of characteristic $p>0$. 
Let $S$ be a geometrically integral canonical K3-like surface. 
Let $L$ be an ample Cartier divisor on $S$ such that $|L|$ is not base point free. 
Then the following hold. 
\begin{enumerate}
\item $\Bs\,|L|$ is irreducible. 
\item $\dim \Bs\,|L|=0$ or $\dim \Bs\,|L|=1$. 
\item If $\dim \Bs\,|L| =1$, then $S$ is smooth around $\Bs\,|L|$. 
\item If $\dim \Bs\,|L| =0$, then $\Bs\,|L|$ is {\cred scheme-theoretically equal to} a reduced point.  
\end{enumerate}
\end{thm}

\begin{proof}
Let $\mu : T \to S$ be the minimal resolution of $S$. 
We then have the following equation of sets: 
\[
\mu^{-1}(\Bs\,|L|) = \Bs\,|\mu^*L| = F, 
\]
where $F$ is a prime divisor on $T$ (Theorem \ref{t-K3-Bs-reg}). 
By the following equation of sets 
\[
\Bs\,|L| =\mu(\mu^{-1}(\Bs\,|L|)) = \mu(F),
\]
$\Bs\,|L|$ is an irreducible closed subset whose dimension is zero or one. 
Thus (1) and (2) hold. 

Let us show (3). 
Assume $\dim \Bs\,|L| =1$. 
In this case, we obtain $\mu^*L \cdot F>0$. 
It follows from  Theorem \ref{t-K3-Bs-reg}(7)  that $\mu^*L \cdot C >0$ for any curve $C$ on $T$ with $C \cap F \neq \emptyset$. 
In other words, $\Ex(\mu) \cap F = \emptyset$, i.e., $\mu : T \to S$ is isomorphic around $F$. 
Since $\mu(F)$ is a smooth Cartier divisor on $S$ (Theorem \ref{t-K3-Bs-reg}), 
$S$ is smooth around $\mu(F)$. 
Thus (3) holds. 

Let us show (4). 
Taking the base change to the separable closure of $\kappa$, 
we may assume that $\kappa$ is separably closed. 
Assume $\dim \Bs\,|L| =0$. 
In this case, $P := \mu(F) = (\Bs\,|L|)_{\red}$ is one point by (1). 
Hence we get $\mu^*L \cdot F=0$. 
By Theorem \ref{t-K3-Bs-reg}, we have $L^2 = 2$, $g:= \frac{1}{2}L^2 +1 =2$, 
$\mu^*L = M+F \sim E_1 + E_2 +F$,  {\cred and} 
\[
\varphi_{|M|} : T \xrightarrow{\pi} B =\P^1_{\kappa} \hookrightarrow \P^2_{\kappa}, 
\]
where $F$ is a section of $\pi$ and $E_1$ and $E_2$ are fibres of $\pi$ over general ${\cred \kappa}$-rational points of $B$. 
Set $E_{1, S} := \mu_* E_1$ and $E_{2, S} := \mu_*E_2$, so that we get $L \sim E_{1, S} + E_{2, S}$. 
We have the induced morphisms $\mu|_{E_i} : E_i \to E_{i, S}$. 
The remaining proof is divided into the following five steps. 
\begin{enumerate}
\renewcommand{\labelenumi}{(\roman{enumi})}
\item For each $i \in \{1, 2\}$, 
it holds that $(\mu|_{E_i})_* \MO_{E_i} = \MO_{E_{i, S}}$, 
i.e., 
$\mu|_{E_i} : E_i \to E_{i, S}$ is an isomorphism. 
\item Each $E_{i, S}$ is a projective Gorenstein curve such that 
$E_{i, S}$ is regular around $P$ and 
$\dim_{\kappa} H^0(E_{i, S}, \MO_{E_{i, S}})= \dim_{\kappa} {\cred H^1(E_{i, S}, \MO_{E_{i, S}})}=1$. 
\item The scheme-theoretic equation $E_{1, S} \cap E_{2, S} = P$ holds. 
\item The scheme-theoretic equation $\Bs\,|L|_{E_{1, S} + E_{2, S}}| \cap E_{i, S}  = \Bs\,|L|_{E_{i, S}} | = P$ 
holds for each $i \in \{1, 2\}$. 
\item The scheme-theoretic equation $\Bs\,|L| =P$ holds. 
\end{enumerate}

Let us show (i). 
Fix $i \in \{1, 2\}$. 
It suffices to show $R^1\mu_*\MO_T(-E_i)=0$. 
This follows from the relative Kawamata--Viehweg vanishing theorem \cite[Theorem 10.4]{Kol13}, 
which can be applied by  
\[
\left(-E_i -\left(K_T +\frac{2}{3} F\right)\right) \cdot F = -1 + 0 + \frac{4}{3}  >0. 
\]
Thus (i) holds. 

Let us show (ii). 
Fix $i \in \{1, 2\}$. 
Note that fibres over general $\kappa$-rational points of $B = \P^1_{\kappa}$ are irreducible. 
By $E_i \cdot F=1$, $E_i$ is a prime divisor and 
$E_i \cap F$ is a $\kappa$-rational point, which implies that $H^0(E_i, \MO_{E_i})=\kappa$. 
By the adjunction formula $\omega_{E_i} \simeq \MO_T(K_T+E_i)|_{E_i} \simeq \MO_{E_i}$, 
we obtain $H^1(E_i, \MO_{E_i}) \simeq \kappa$ by Serre duality. 
It follows from $E_i \cdot F=1$ that $E_i$ is regular around $E_i \cap F$, 
and hence $E_{i, S}$ is regular around $P$ by (i). 
Thus (ii) holds.

Let us show (iii). 
By $2 = L^2  = L \cdot (E_{1, S}+E_{2, S})$ and the ampleness of $L$, 
it holds that $L \cdot E_{1, S}=L \cdot E_{2, S}=1$. 
We have $L \sim 2E_{2, S}$ and $E_{1, S} \cdot (2E_{2, S}) = 1$, 
which implies the scheme-theoretic equation $E_{1, S} \cap (2E_{2, S}) = P$. 
Therefore, we obtain the following scheme-theoretic inclusions: 
\[
\{ P \} \subset E_{1, S} \cap E_{2, S} \subset E_{1, S} \cap (2E_{2, S})= \{P\}, 
\]
and hence $E_{1, S} \cap E_{2, S} = P$. 
Thus (iii) holds. 

Let us show (iv). By (iii), we have the following exact sequence: 
\[
0 \to H^0(E_{1, S} +E_{2, S}, L|_{E_{1, S} +E_{2, S}}) \to H^0(E_{1, S}, L|_{E_{1, S}}) \oplus H^0(E_{2, S}, L|_{E_{2, S}}) \to H^0(P, L|_P).
\]
Note that $L|_{E_{1, S}} \simeq \MO_T(2E_{2, S})|_{E_{1, S}} \simeq \MO_{E_{1, S}}(P)$ by $L \cdot E_{1, S}=1$. 
By the Riemann--Roch theorem, we obtain $h^0(E_{1, S}, \MO_{E_{1, S}}(P)) =1$ (note that $E_{1, S}$ is regular around $P$ by (ii)). 
Hence we get the isomorphism: 
\[
H^0(E_{1, S}, \MO_{E_{1, S}}) \xrightarrow{\simeq} H^0(E_{1, S}, \MO_{E_{1, S}}(P)), 
\]
which implies that $H^0(E_{1, S}, L|_{E_{1, S}}) \to H^0(P, L|_P)$ is zero. 
Similarly, $H^0(E_{2, S}, L|_{E_{2, S}}) \to H^0(P, L|_P)$ is zero. 
Therefore, {$\Bs\,|L|_{E_{i, S}}| = P$ for each $i \in \{1, 2\}$ and} we get the induced isomorphism: 
\[
 H^0(E_{1, S}  +E_{2, S}, L) \xrightarrow{\simeq} H^0(E_{1, S}, L) \oplus H^0(E_{2, S}, L). 
\]
Then each projection $H^0(E_{1, S}  +E_{2, S}, L|_{E_{1, S} +E_{2, S}}) \to H^0(E_{i, S}, L|_{E_{i, S}})$ is surjective. 
Therefore, (iv) holds. 

Let us show (v). 
We have an exact sequence
\[
H^0(S, L) \to H^0(E_{1, S}+E_{2, S}, L|_{E_{1, S}+E_{2, S}}) \to H^1(S, \MO_S)=0, 
\]
and hence 
\[
\Bs\,|L| = \Bs\,|L|_{E_{1, S}+E_{2, S}}|. 
\]
For the ideals $\m_P \subset \MO_{E_{1, S}+E_{2, S}}$ and 
$I_{\Bs |L|_{E_{1, S}+E_{2, S}}|} \subset \MO_{E_{1, S}+E_{2, S}}$ 
of $P$ and $\Bs |L|_{E_{1, S}+E_{2, S}}|$ respectively, 
it suffices to show that $\m_P = I_{\Bs |L|_{E_{1, S}+E_{2, S}}|}$. 
By the induced injection 
\[
\rho := \rho_1 \oplus \rho_2 : \MO_{E_{1, S}+E_{2, S}} \hookrightarrow \MO_{E_{1, S}} \oplus \MO_{E_{2, S}} 
\quad \text{where}\quad 
\rho_i: \MO_{E_{1, S}+E_{2, S}} \to \MO_{E_{i, S}},  
\]
it is enough to check $\rho(\m_P) = \rho( I_{\Bs |L|_{E_{1, S}+E_{2, S}}|})$, 
which is equivalent to 
$\rho_1(\m_P) = \rho_1( I_{\Bs |L|_{E_{1, S}+E_{2, S}}|})$ 
and 
$\rho_2(\m_P) = \rho_2( I_{\Bs |L|_{E_{1, S}+E_{2, S}}|})$. 
These follow from (iv). 
Thus (v) holds. 
\end{proof}

\begin{cor}\label{c-K3-to-Fano}
We work over an algebraically closed field $k$ of characteristic $p>0$. 
Let $X$ be a Fano threefold. 
Assume that the generic member $S$ of $|-K_X|$ is a geometrically integral canonical K3-like surface. 
Then $S$ is regular. 
\end{cor}

\begin{proof}
Set $\kappa := K(\P(H^0(X, -K_X)))$, $X_{\kappa} :=X \times_k \kappa$, 
and $L := (-K_{X_{\kappa}})|_S$, which is an ample Cartier divisor on $S$. 
Note that $S$ is a prime divisor on $X_{\kappa}$ and 
$S$ is regular outside $\Bs\,|-K_{X_{\kappa}}|$ (Theorem \ref{t-Bertini-generic}). 
{\cred It holds that $\Bs\,|-K_{X_{\kappa}}| = \Bs\,|L|$, because we have an exact sequence  
\[
H^0(X_{\kappa}, \MO_{X_{\kappa}}(-K_{X_{\kappa}})) \to H^0(S, L) \to H^1(X_{\kappa}, \MO_{X_{\kappa}}) 
\]
and $H^1(X_{\kappa}, \MO_{X_{\kappa}}) \simeq H^1(X, \MO_X) \otimes_k \kappa =0$ (Theorem \ref{t-kawakami}).}  
By Theorem \ref{t-K3-Bs-cano}(2), 
there are the following two cases: $\dim \Bs\,|L|=1$ and $\dim \Bs\,|L|=0$. 
If $\dim \Bs\,|L|=1$, then $S$ is smooth around $\Bs\,|L|$ (Theorem \ref{t-K3-Bs-cano}(3)), and hence $S$ is regular. 
We may assume that $\dim \Bs\,|L|=0$. Then $\Bs\,|-K_{X_{\kappa}}| = \Bs\,|L|$ is a reduced point $P$  (Theorem \ref{t-K3-Bs-cano}(4)). 
By $\Bs\,|-K_X| \times_k \kappa \simeq \Bs\,|-K_{X_{\kappa}}|$, 
also $\Bs\,|-K_X|$ is a reduced point $Q$, which is nothing but the image of $P$. 
Therefore, general members of $|-K_X|$ are smooth at $Q$  by Lemma \ref{l-Bs-red-pt} below. 
Hence $S$ is smooth around $\Bs\,|L|$, as required. 
\end{proof}

\begin{lem}\label{l-Bs-red-pt}
We work over an algebraically closed field $k$. 
Let $X$ be a smooth projective variety and fix a closed point $Q \in X$. 
Let $L$ be a Cartier divisor such that the scheme-theoretic equation 
\[
\Bs\,|L| = Q \amalg W
\]
holds for some closed subscheme $W$ of $X$ with $Q \not\in W$. 
Then general members of $|L|$ are smooth at $Q$. 
\end{lem}

\begin{proof}
Let $\m_Q \subset \MO_{X, Q}$ be the maximal ideal. 
By $\Bs\,|L| = Q \amalg W$, 
there exists an effective divisor $D \in |L|$ such that $f_D \in \m_Q$ and 
$f_{\cred D} \not\in \m_Q^2$, 
where $f_D \in \MO_{X, Q}$ is an defining element of $D$ which is uniquely determined 
up to $\MO_{X, Q}^{\times}$. 
Hence $D$ is smooth at $Q$. Therefore, also general members of $|L|$ are smooth at $Q$. 
\end{proof}

\section{Generic elephants}\label{s-generic-elephant}

\subsection{The case $\dim \Im\,\varphi_{|-K_X|} =1$}

\begin{prop}\label{p-composite-pencil}
We work over an algebraically closed field {\cred $k$} of characteristic $p>0$. 
Let $X$ be a Fano threefold. 
Assume that $\dim (\Im\,\varphi_{|-K_X|}) =1$. 
Then the following hold. 
\begin{enumerate}
    \item The mobile part $|D|$ of $|-K_X|$ is base point free. 
    \item For  the image $Y := \Im\,\varphi_{|D|}$ and the morphism $\psi : X \to Y$ 
    induced by $\varphi_{|D|}$, it holds that $Y \simeq \P^1$ and $\psi_*\MO_X = \MO_Y$. 
    \item If $|-K_X|$ is not base point free, then $(-K_X)^2 \cdot E=1$ for a fibre $E$ of $\psi$. 
\end{enumerate}
\end{prop}

\begin{proof}
Take the decomposition 
\[
|-K_X| = D_0 +|D|
\]
into the fixed part $D_0$ of $|-K_X|$ and the mobile part $D$ of $|-K_X|$.  
Set 
\[
g := \frac{1}{2} (-K_X)^3 +1 = \chi(X, -K_X) -2\qquad \text{and} \qquad g' := h^0(X, -K_X) -2. 
\]
By $H^2(X, -K_X)=H^3(X, -K_X)=0$ (Theorem \ref{t-kawakami}), we have 
\[
2 \leq g = \chi(X, -K_X) - 2 = h^0(X, -K_X) -h^1(X, -K_X) - 2 \leq h^0(X, -K_X) - 2 = g'. 
\]
We then have $\varphi := \varphi_{|-K_X|}:X \dashrightarrow \mathbb P^{g'+1}$. 
Let $\sigma: \widetilde{X} \to X$ be a desingularisation of the resolution of the indeterminacies of $\varphi_{|-K_X|}$, so that we have $\widetilde \varphi : \widetilde X \to \mathbb P^{g'+1}$. 
For the largest open subset $X'$ of $X$ on which $\varphi_{|-K_X|}$ is defined, 
it holds that $\dim (X \setminus X') \leq 1$ and 
$\widetilde X' := \sigma^{-1}(X') \to X'$ can be assumed to be an isomorphism. 
Take the decomposition 
\[
\sigma^*D = \widetilde D + F
\]
into the movile part $\widetilde{D}$ and the fixed part $F$, where $|\widetilde D|$ is base point free. 
We then get 
\begin{equation}\label{e1-composite-pencil}
H^0(\widetilde X, \widetilde D) \simeq
H^0(\widetilde X, \sigma^*D) \simeq H^0(X, D) \simeq H^0(X, -K_X)
\end{equation}
and $\varphi_{|\widetilde D|}: \widetilde X \to \mathbb P^{g'+1}$ with 
$\MO_{\widetilde X}(\widetilde D) \simeq \varphi_{|\widetilde D|}^*\MO_{\P^{g'+1}}(1)$. 
Set $Y := \overline{\varphi(X')} = \widetilde \varphi(\widetilde X)$. 
Let 
\[
\widetilde{\varphi} : \widetilde X \xrightarrow{\psi}Y \hookrightarrow \P^{g'+1} 
\]
be the induced morphisms. 

We now show that 
\begin{enumerate}
\item[(a)] $\Delta(Y, \MO_{\P^{g'+1}}(1)|_Y)=0$, 
\item[(b)] $Y \simeq \P^1$, 
\item[(c)] $\deg (\MO_{\P^{g'+1}}(1)|_Y) = g'+1$, and 
\item[(d)] $\psi_* \MO_{\widetilde X} = \MO_Y$. 
\end{enumerate}
Let $\psi : \widetilde X \to Z \to Y$ be the Stein factorisation of $\psi$. By 
\[
H^1(Z, \MO_Z) \hookrightarrow 
H^1(\widetilde X, \MO_{\widetilde X}) \simeq H^1(X, \MO_X)=0, 
\]
we obtain $H^1(Z, \MO_Z) =0$, which in turn implies $Z \simeq \mathbb P^1$. 
In particular, we get 
\begin{equation}\label{e2-composite-pencil}
\deg (\MO_{\P^{g'+1}}(1)|_Z) = h^0(Z, \MO_{\P^{g'+1}}(1)|_Z) -1 =h^0(\widetilde X, \widetilde D) -1 = g'+1, 
\end{equation}
where the first equality holds by the Riemann--Roch theorem, 
the second one follows from Lemma \ref{l-H^0-image}, and 
the last one holds by (\ref{e1-composite-pencil}). 
Set $\ell := \deg (Z \to Y) \in \Z_{>0}$. 
It holds that 
\begin{align*}
0 
&\leq \Delta(Y, \MO_{\P^{g'+1}}(1)|_Y) \\
&= \dim Y + \deg (\MO_{\P^{g'+1}}(1)|_Y) -h^0(Y, \MO_{\P^{g'+1}}(1)|_Y)\\
&= 1 + \frac{\deg (\MO_{\P^{g'+1}}(1)|_Z)}{\ell} - (g'+2)\\
&= \frac{(g'+1)(1-\ell)}{\ell}, 
\end{align*}
where the second equality follows from Lemma \ref{l-H^0-image} and (\ref{e1-composite-pencil}), 
and the third one holds by (\ref{e2-composite-pencil}). 
Then $\ell = 1$ and $\Delta(Y, \MO_{\P^{g'+1}}(1)|_Y) =0$, 
which deduces that $Z \to Y$ is birational and $Y \simeq \P^1$. 
Thus (a) and (b) hold. 
Since $Y$ is normal and $Z \to Y$ is a finite birational morphism, $Z \to Y$ is an isomorphism. 
Thus (d) holds. 
Finally, (c) follows from (\ref{e2-composite-pencil}). 
This completes the proof of (a)--(d).

Fix a closed point $P \in Y$ and set $\widetilde E := \psi^*P$ and $E := \sigma_*\widetilde E$. 
We then obtain $\widetilde{D} \sim (g'+1) \widetilde{E}$ and 
\[
D = \sigma_*\sigma^*D = \sigma_*(\widetilde D +F) = \sigma_*\widetilde D \sim (g'+1)E. 
\]
By $-K_X \sim D_0 +D \sim D_0 +(g'+1)E$, we obtain 
\begin{eqnarray*}
2g -2 
&=& (-K_X)^3 \\
&=& (D_0 +(g'+1)E) \cdot (-K_X)^2\\
&=& D_0 \cdot (-K_X)^2  +(g'+1)E \cdot (-K_X)^2\\
&=& D_0 \cdot (-K_X)^2  +(g'+1)E \cdot (D_0 +(g'+1)E) \cdot (-K_X)\\
&=& D_0 \cdot (-K_X)^2  +(g'+1)E \cdot D_0 \cdot (-K_X) + (g'+1)^2E^2 \cdot (-K_X).
\end{eqnarray*}
Since $D_0$ is effective, $-K_X$ is ample, and 
$|E|$ is base point free outside a closed subset of codimension two,
it holds that 
\[
D_0 \cdot (-K_X)^2 \geq 0, \qquad 
E \cdot D_0 \cdot (-K_X) \geq 0, \qquad 
E^2 \cdot (-K_X) \geq 0. 
\]

We now show that $E^2 \cdot (-K_X) =0$. 
Suppose the contrary, i.e., $E^2 \cdot (-K_X) >0$. 
The above equation deduces 
\[
2g > 2g -2 \geq (g'+1)^2E^2 \cdot (-K_X) \geq (g'+1)^2 \geq 2g',  
\]
which contradicts $g' \geq g$. Hence, we have $E^2 \cdot (-K_X) =0$. 

Pick two general members $\widetilde{E}_1, \widetilde{E}_2 \in |\widetilde{E}|$, which are distinct prime divisors \cite[Corollary 7.3]{Bad01}. 
Set $E_i := \sigma_*\widetilde{E}_i$. 
Then $E_1$ and $E_2$ are distinct prime divisors. 
We get 
\[
E \sim E_1 \sim E_2. 
\]
Therefore, we obtain 
\[
E_1 \cdot E_2 \cdot (-K_X) = E^2 \cdot (-K_X) =0, 
\]
which implies that $E_1 \cap E_2 =\emptyset$. 
Thus $|E|$ is base point free, and hence so is $|(g'+1)E| = |D|$. 
{\cred In particular, $\sigma: \widetilde X \xrightarrow{\simeq} X$.} 
Hence (1) and (2) holds. 

Let us show (3). 
Assume that $|-K_X|$ is not base point free, 
which is equivalent to $D_0 \neq 0$. 
Since $-K_X$ is ample, we get
$D_0 \cdot (-K_X)^2 >0$ and $E \cdot (-K_X)^2 >0$. 
The latter one implies 
\[
E \cdot D_0 \cdot (-K_X) = E \cdot  (D_0 +(g'+1)E) \cdot (-K_X) = E \cdot (-K_X)^2 >0. 
\]
We have
\[
2g-2 = D_0 \cdot (-K_X)^2 + (g'+1) E \cdot D_0 \cdot (-K_X). 
\]
By $g' \geq g$, we get $E \cdot D_0 \cdot (-K_X) =1$, and hence $E \cdot (-K_X)^2 =1$. 
Thus (3) holds. 
\qedhere 


\end{proof}


\begin{cor}\label{c-composite-pencil}
We work over an algebraically closed field {\cred $k$} of characteristic $p>0$. 
Let $X$ be a Fano threefold with $\rho(X)=1$. 
Then $\dim (\Im\,\varphi_{|-K_X|}) \geq 2$. 
\end{cor}

\begin{proof}
Suppose that $\dim (\Im\,\varphi_{|-K_X|}) \leq 1$. 
By $h^0(X, -K_X)  \geq \chi(X, -K_X) = \frac{(-K_X)^3}{2} + 3 \geq 2$, 
we get $\dim (\Im\,\varphi_{|-K_X|}) = 1$. 
It follows from Proposition \ref{p-composite-pencil} that there is a surjective morphism $X \to {\cred Y}$ to a curve ${\cred Y}$. 
This contradicts $\rho(X)=1$. 
\end{proof}

\subsection{The case $\dim \Im\,\varphi_{|-K_X|} \geq 2$}

\begin{nasi}[Mumford pullback]\label{n-Mumford}
We work over a field $\kappa$. 
Let 
\[
\mu : T^{\min} \to T
\] 
be the minimal resolution of a normal surface $T$ and let $E_1, ..., E_n$ be all the $\mu$-exceptional prime divisors, i.e., $\Ex(\mu) = E_1 \cup \cdots \cup E_n$.  
\begin{enumerate}
\item Given a $\Q$-divisor $D$ on $T$, we define $\mu^*D$ as a unique $\Q$-divisor on $T^{\min}$ satisfying the following properties {\cred (a) and (b)}. 
\begin{enumerate}
\item $\mu_*\mu^*D = D$. 
\item $\mu^*D \cdot E_1 = \cdots = \mu^*D \cdot E_n =0$. 
\end{enumerate}
{\cred Equivalently, for the proper transform $D'$ of $D$ on $T^{\min}$, 
$\mu^*D$ is given by $\mu^*D := D' + \sum_{i=1}^n e_iE_i$, 
where the rational numbers $e_1, ..., e_n$ are uniquely determined by (b), 
because the $n \times n$  matrix $(E_i \cdot E_j)$ is invertible \cite[Theorem 10.1]{Kol13}.}
\item There exists an effective $\Q$-divisor $\Delta$ such that 
\[
K_{T^{\min}} + \Delta  = \mu^*K_T. 
\]
Furthermore, $T$ is canonical if and only if $\Delta=0$ {\cred \cite[Theorem 4.13(1)]{Tan18}}. 
\end{enumerate}
\end{nasi}

\begin{thm}\label{t-generic-ele}
We work over an algebraically closed field {\cred $k$} of characteristic $p>0$. 
Let $X$ be a Fano threefold. 
Assume that (I) or (II) holds. 
\begin{enumerate}
\item[(I)] $\rho(X) =1$. 
\item[(II)] $\dim ({\rm Im}\,\varphi_{|-K_X|}) \geq 2$. 
\end{enumerate}
Let 
\[
|-K_X| = {\cred |M|} +F
\]
be the decomposition into the mobile part $M$ of $|-K_X|$ and the fixed part $F$ of $|-K_X|$. 
Then $F=0$ and the generic member $S$ of $|-K_X|$ is a geometrically integral regular K3-like surface. 
\end{thm}


\setcounter{step}{0}
\begin{proof}
Since (I) implies (II) (Corollary \ref{c-composite-pencil}), we may assume that (II) holds. 
Then general members of $|M|$ are prime divisors (Proposition \ref{p-Bertini-integral}). 
Replacing $M$ by a general member of $|M|$, 
the problem is reduced to the case when $M$ is a prime divisor. 
By \cite[Proposition 4.2]{CP08}, there exists a sequence of blowups 
\begin{equation}\label{e1-generic-ele}
\sigma: \widetilde X := X_{\ell} \xrightarrow{\sigma_{\ell}} X_{\ell-1} \xrightarrow{\sigma_{\ell-1}} 
\cdots 
\xrightarrow{\sigma_2} X_1 
\xrightarrow{\sigma_1} X_0 := X 
\end{equation}
such that 
\begin{itemize}
\item for each $i$, $\sigma_i: X_i \to X_{i-1}$ is a blowup along either a point or a smooth curve, 
\item the centre $\sigma_i(\Ex(\sigma_i))$ is contained in the proper transform $M_i$ of $M$, and 
\item the mobile part $\widetilde M$ of $|\sigma^*M|$ is base point free. 
\end{itemize}
Set $k' := K(\P(H^0(X, -K_X)))$, 
which is a purely transcendental extension over $k$ of finite degree. 
In what follows, we set $Y' := Y \times_k k'$ for any $k$-scheme $Y$. 
Applying the base change $(-) \times_k k'$ to (\ref{e1-generic-ele}), we get a sequence of blowups (note that blowups commute with flat base changes \cite[Section 8, Proposition 1.12(c)]{Liu02}): 
\[
\sigma': \widetilde X' := X'_{\ell} \xrightarrow{\sigma'_{\ell}} X'_{\ell-1} \xrightarrow{\sigma'_{\ell-1}} 
\cdots 
\xrightarrow{\sigma'_2} X'_1 
\xrightarrow{\sigma'_1} X'_0 = X'. 
\]
%
Let $S$ and $\widetilde {S}$ be the generic members of $|M|$ and $|\widetilde M|$, respectively. 
By 
\[
H^0(X, -K_X) \simeq H^0(X, M) \simeq H^0(\widetilde X, \sigma^*M) \simeq H^0(\widetilde X, \widetilde M),
\]
we have $S \sim M' := M \times_k k'$, 
$\widetilde{S} \sim \widetilde M' := \widetilde M \times_k k'$, and 
$\sigma'_* \widetilde{S} = S$. 
As general members of each of $|M|$ and $|\widetilde M|$ are prime divisors 
(Proposition \ref{p-Bertini-integral}), 
both $S$ and $S'$ are geometrically integral prime divisors 
{\cred (indeed, for $V := H^0(X, M)$ and the family 
$\alpha: X^{\univ}_{M, V} \to \P(V)$ parametrising all the members of $|M|$, 
general members of $|M|$ coincide with fibres over closed points of $\alpha$ and 
the generic member $S$ of $|M|$ is nothing but the generic fibre of $\alpha$ 
(cf. Notation \ref{n-generic-member}). 
Since general fibres of $\alpha$ are geometrically integral, 
so is the generic fibre $S$ \cite[Th\'eor\`eme 12.2.1(x)]{EGAIV3}. 
The same argument implies that also $S'$ is geometrically integral)}. 
Furthermore, $\widetilde S$ is regular, since $|\widetilde M|$ is base point free 
(Theorem \ref{t-Bertini-generic}).  
Let $E_i$ be the prime divisor on $\widetilde{X}$ that arises as the $i$-th blowup, 
i.e., $E_i \subset \widetilde X$ is the proper transform of $\Ex(\sigma_i) \subset X_i$. 
There exist $n_i \in \Z_{>0}$ and $a_i \in \Z_{>0}$ such that the following hold. 
\begin{enumerate}
\item $\sigma^*M = \widetilde M +\sum_i n_iE_i$. 
\item $\sigma'^*S = \widetilde S +\sum_i n_iE'_i$.  
\item $K_{\widetilde X} = \sigma^*K_{X} + \sum_i a_i E_i$. 
\item $K_{\widetilde X'} = \sigma'^*K_{X'} + \sum_i a_i E'_i$. 
\end{enumerate}

\begin{step}\label{s1-generic-ele}
If $i$ is an index such that the image $\sigma(E_i)$ on $X$ is one-dimensional, then 
it holds that $n_i \geq a_i$. 
\end{step}

\begin{proof}[Proof of Step \ref{s1-generic-ele}]
Fix such $i$. 
In order to prove $n_i \geq a_i$, 
we may work with an open neighbourhood of the generic point $\xi_i$ of $\sigma(E_i)$. 
Over a suitable open neighbourhood $U$ of $\xi_i$, 
$\sigma|_{\sigma^{-1}(U)} : \sigma^{-1}(U) \to U$ is obtained by a sequence of blowups along smooth curves. 
We then inductively obtain 
\[
K_{X_i} +M_j + (\text{effective divisor}) = \sigma_j^*(K_{X_{j-1}}+M_{j-1}), 
\]
where $M_0 := M$ and $M_j$ denotes the proper transform of $M$. 
This implies $n_i \geq a_i$, which completes the proof of Step \ref{s1-generic-ele}. 
\end{proof}

Let $\nu : S^N \to S$ be the normalisation of $S$ and 
let $\mu : S^{\min} \to S^N$ be the minimal resolution of $S^N$. 
Since $\widetilde S$ is regular, 
$\sigma'|_{\widetilde S} : \widetilde S \to S$ factors through $S^{\min}$:  
\[
\sigma'|_{\widetilde S} : \widetilde S \xrightarrow{\lambda} S^{\min} \xrightarrow{\mu} S^N \xrightarrow{\nu} S. 
\]

\begin{step}\label{s2-generic-ele}
The following hold. 
\begin{enumerate}
\setcounter{enumi}{4}
\item $K_{S^N} \sim -D$ for some effective $\Z$-divisor $D$ on $S^N$. 
\item $\kappa(\widetilde S, K_{\widetilde S}) = \kappa(S^{\min}, K_{S^{\min}}) \leq \kappa(S^N, K_{S^N}) \leq 0$. 
\end{enumerate}
\end{step}

\begin{proof}[Proof of Step \ref{s2-generic-ele}]
By (2) and (4), we get 
\begin{equation}\label{e2-generic-ele}
K_{\widetilde X'} +\widetilde S = {\cred \sigma'}^*(K_{X'}+S) + \sum_i (a_i -n_i) E'_i 
=-\sigma'^*F' + \sum_i (a_i -n_i) E'_i. 
\end{equation}
Taking the pushforward of $(\ref{e2-generic-ele})|_{\widetilde S}$ to $S^N$ by $\mu \circ \lambda : \widetilde S \to S^N$, we obtain 
\[
K_{S^N}  \sim (\mu \circ \lambda)_*K_{\widetilde{S}} \sim 
-(\mu \circ \lambda)_*({\cred \sigma'}^*F|_{\widetilde S}) +\sum_i (a_i-n_i) (\mu \circ \lambda)_*(E_i|_{\widetilde S}) \leq 0, 
\]
where ${\cred \sigma'}^*F|_{\widetilde S}$ is clearly effective and 
$\sum_i (n_i-a_i) (\mu \circ \lambda)_*(E_i|_{\widetilde S})$ is effective  by Step \ref{s1-generic-ele}. 
Thus (5) holds. 

Let us show (6). 
It follows from (5) that $\kappa(S^N, K_{S^N}) \leq 0$. 
We have 
\[
H^0(S^{\min}, mK_{S^{\min}}) \hookrightarrow H^0(S^N, mK_{S^N}), 
\]
which implies $\kappa(S^{\min}, K_{S^{\min}}) \leq \kappa(S^N, K_{S^N})$. 
It is well known that $\kappa(\widetilde{S}, K_{\widetilde{S}}) = \kappa(S^{\min}, K_{S^{\min}})$. 
Thus (6) holds. 
This completes the proof of Step \ref{s2-generic-ele}. 
\end{proof}

Set $A := {\cred \sigma'}^*(-K_{X'})|_{\widetilde S}$ and 
\[
L := K_{\widetilde S} + A = 
\left(K_{\widetilde X'}+\widetilde S +\sigma'^*(-K_{X'})\right)\Big|_{\widetilde S} 
=\left(\widetilde S + \sum_i a_i E'_i\right)\Big|_{\widetilde S}. 
\]
In particular, 
\begin{equation}\label{e3-generic-ele}
h^0(\widetilde S, L) \geq h^0(\widetilde S, \MO_{\widetilde{X}'}(\widetilde S)|_{\widetilde S}). 
\end{equation}

\begin{step}\label{s3-generic-ele}
It holds that $h^0(\widetilde S,A) \geq 2$. 
\end{step}

\begin{proof}[Proof of Step \ref{s3-generic-ele}]
It is enough to check the following inequalities:  
\[
h^0(\widetilde S,A) = h^0(\widetilde S, {\cred \sigma'}^*\MO_{X'}(-K_{X'})|_{\widetilde S}) 
\geq 
h^0(S, \MO_{X'}(-K_{X'})|_{S}) \geq h^0(S, \MO_{X'}(S)|_S) \geq 2. 
\]
It suffices to show the last inequality, since the other ones are obvious. 
By $H^1(X', \MO_{X'})=0$ and the exact sequence 
\[
0 \to \MO_{X'} \to \MO_{X'}(S) \to \MO_{X'}(S)|_S \to 0, 
\]
we get  
\[
h^0(S, \MO_{X'}(S)|_S) = h^0(X', S) -h^0(X', \MO_{X'}) = h^0(X', -K_{X'}) -1
\]
\[
\geq \chi(X, -K_X) -1= \frac{(-K_X)^3}{2} +2\geq 2. 
\]
This completes the proof of Step \ref{s3-generic-ele}. 
\end{proof}

\begin{step}\label{s4-generic-ele}
For $g := \frac{(-K_X)^3}{2} +1$ and $g' := h^0(X, -K_{X})-2$, it holds that 
\[
g -1 \geq  
h^0(\widetilde S, L) -1 - h^2(\widetilde S, \MO_{\widetilde S}) 
\geq 
g' -h^2(\widetilde S, \MO_{\widetilde S}). 
\]
Furthermore, if $g=g'$ and $h^2(\widetilde S, \MO_{\widetilde S}) =1$, 
then $F=0$. 
\end{step}

\begin{proof}[Proof of Step \ref{s4-generic-ele}]
{\cred In order to show the required inequalities,} 
it is enough to prove the following inequalities: 
\begin{align*}
2g -2 &\overset{(\alpha)}{\geq}  
\sigma'^*(-K_{X'}) \cdot \left(\widetilde S+ \sum_i a_i E'_i\right) \cdot \widetilde S\\
&\overset{(\beta)}{\geq}  2h^0(\widetilde S, L) -2 - 2h^2(\widetilde S, \MO_{\widetilde S}) \\
&\overset{(\gamma)}{\geq}  2g' -2h^2(\widetilde S, \MO_{\widetilde S}). 
\end{align*}

Let us show $(\alpha)$. It holds that 
\begin{eqnarray*}
2g -2 
&=& \sigma'^*(-K_{X'})^3\\
&=& \sigma'^*(-K_{X'})^2 \cdot \left(\sigma'^*F' + \widetilde S + \sum_i n_i E'_i\right)\\
&\geq & \sigma'^*(-K_{X'})^2 \cdot \widetilde S\\
&= & \sigma'^*(-K_{X'}) \cdot \left(\sigma'^*F' + \widetilde S + \sum_i n_i E'_i\right) \cdot \widetilde S\\
&\geq & \sigma'^*(-K_{X'}) \cdot \left(\widetilde S + \sum_i a_i E'_i\right) \cdot \widetilde S, 
\end{eqnarray*}
where the last inequality holds by $\sigma'^*(-K_{X'}) \cdot \sigma'^*F' \cdot \widetilde S \geq 0$ and 
the following facts. 
\begin{itemize}
\item If $\dim \sigma'(E'_i)=0$, then $\sigma'^*(-K_{X'}) \cdot E'_i \cdot \widetilde S =0$.
\item If $\dim \sigma'(E'_i)=1$, then $n_i \geq a_i$ (Step \ref{s1-generic-ele}). 
\end{itemize}
Thus $(\alpha)$ holds. 

Let us show $(\beta)$. 
By the Riemann--Roch theorem, we get  
\begin{align*}
\chi(\widetilde S, L) 
&= \chi (\widetilde S, \MO_{\widetilde S}) + \frac{1}{2}L \cdot (L-K_{\widetilde S}) \\
&= \chi (\widetilde S, \MO_{\widetilde S}) + \frac{1}{2}A \cdot (K_{\widetilde S} + A ) \\
&= \chi (\widetilde S, \MO_{\widetilde S}) + \frac{1}{2}\sigma'^*(-K_{X'})|_{\widetilde S} \cdot (K_{\widetilde X'} + \widetilde S + \sigma'^*(-K_{X'}) )|_{\widetilde S}  \\
&=\chi (\widetilde S, \MO_{\widetilde S}) + \frac{1}{2}\sigma'^*(-K_{X'}) \cdot \left( \widetilde S + \sum_i a_i E'_i\right) \cdot \widetilde S. 
\end{align*}
As $A = \sigma'^*(-K_{X'})|_{\widetilde S}$ is big, 
we get 
\[
h^2(\widetilde S, L) = h^2(\widetilde S, K_{\widetilde S}+A) = h^0(\widetilde S, -A)=0.
\]
Since $\widetilde S$ is geometrically integral and $h^0(\widetilde S, A) \geq 2$ (Step \ref{s3-generic-ele}), 
it follows from Theorem \ref{t-Mumford-vanishing} that $h^1(\widetilde S, L) =h^1(\widetilde S, -A) \leq h^1(\widetilde S, \MO_{\widetilde S})$. 
Hence  
\begin{align*}
&\,\,\frac{1}{2}\sigma'^*(-K_{X'}) \cdot \left( \widetilde S + \sum_i a_i E'_i\right) \cdot \widetilde S \\
=&\,\,\chi(\widetilde S, L) - \chi (\widetilde S, \MO_{\widetilde S})\\
=&\,\,h^0(\widetilde S, L) -h^1(\widetilde S, L) 
-h^0(\widetilde S, \MO_{\widetilde S})+h^1(\widetilde S, \MO_{\widetilde S})-h^2(\widetilde S, \MO_{\widetilde S})\\
\geq&\,\,h^0(\widetilde S, L)  -1-h^2(\widetilde S, \MO_{\widetilde S}). 
\end{align*}
Note that we have $h^0(\widetilde S, \MO_{\widetilde S})=1$, since $\widetilde S$ is geometrically  integral. 
Thus $(\beta)$ holds. 

Let us show $(\gamma)$. 
By $H^1(\widetilde X', \MO_{\widetilde X'})=0$ and the following exact sequence 
\[
0 \to \MO_{\widetilde X'} \to \MO_{\widetilde X'}(\widetilde S) \to \MO_{\widetilde X'}(\widetilde S)|_{\widetilde S} \to 0, 
\]
we obtain 
\[
h^0(\widetilde S, L) \geq 
h^0(\widetilde S, \MO_{\widetilde{X}'}(\widetilde S)|_{\widetilde S}) = h^0(\widetilde X', \widetilde S) -1 = h^0(X', -K_{X'}) -1
= g'+1.  
\]
where the first inequality is guaranteed by (\ref{e3-generic-ele}). 
Thus $(\gamma)$ holds. 

Finally, assuming that $g=g'$ and $h^2(\widetilde S, \MO_{\widetilde S}) =1$, 
it is enough to show that $F=0$. 
In this case, all the inequalities $(\alpha), (\beta), (\gamma)$ are equalities. 
By the proof of $(\alpha)$, we get $\sigma'^*(-K_{X'})^2 \cdot \sigma'^*F'=0$, 
which implies $F'=0$, i.e., $F=0$. 
This completes the proof of Step \ref{s4-generic-ele}. 
\end{proof}


\begin{step}\label{s5-generic-ele}
The following hold. 
\begin{enumerate}
\setcounter{enumi}{6}
\item $g=g'$. 
\item $h^2(\widetilde S, \MO_{\widetilde S})=1$. 
\item $S$ is a geometrically integral canonical K3-like surface. 
\end{enumerate}
\end{step}

\begin{proof}[Proof of Step \ref{s5-generic-ele}]
It follows from (6) that $h^2(\widetilde S, \MO_{\widetilde S})= h^0(\widetilde S, K_{\widetilde S}) \leq 1$. 
By $g' \geq g$ and Step \ref{s4-generic-ele}, we obtain 
\[
 g \geq h^0(\widetilde S, L)  -h^2(\widetilde S, \MO_{\widetilde S}) \geq g'+1- h^2(\widetilde S, \MO_{\widetilde S}) \geq g + 1 - 1 = g. 
\]
Then all the inequalities in this equation are equalities, and hence  
\[
g =g'\qquad \text{and} \qquad h^2(\widetilde S, \MO_{\widetilde S})= h^0(\widetilde S, K_{\widetilde S}) = 1. 
\]
Thus (7) and (8) hold. 

Let us show (9). 
By (5), we have $K_{S^N} \sim -D$ for some effective $\Z$-divisor $D$. 
On the other hand, we get 
$h^0(S^{\min}, K_{S^{\min}}) = h^0(\widetilde S, K_{\widetilde S}) =1$, which implies 
$K_{S^{\min}} \sim E$ for some effective $\Z$-divisor $E$. 
There exists an effective $\Q$-divisor $\Delta$ on $S^{\min}$ such that 
$K_{S^{\min}} + \Delta = \mu^*K_{S^N}$  (\ref{n-Mumford}). 
We then get 
\[
E + \Delta \equiv K_{S^{\min}} +\Delta  = \mu^*K_{S^N} \equiv \mu^*(-D), 
\]
and hence $E + \Delta + \mu^*D \equiv 0$. 
Since all of $E, \Delta, \mu^*D$ are effective $\Q$-divisors, 
we obtain $E = \Delta = \mu^*D =0$ by taking the intersection with an ample divisor on $S^{\min}$. 
In particular, $D  = 0$, $K_{S^N} \sim 0$, and $K_{S^{\min}} \sim 0$. Then $K_{S^{\min}} \sim \mu^*K_{S^N}$, which implies that $S^N$ is canonical (\ref{n-Mumford}).



We have 
\[
\omega_{S} \simeq \MO_{X'}(K_{X'}+S)|_S \simeq \MO_{X'}(-F')|_S. 
\]
Since $S$ is Gorenstein, we obtain $\nu^*\omega_S \simeq \MO_{S^N}(K_{S^N} + C)$ 
for the conductor $C$, which is an effective $\Z$-divisor on $S^N$ such that $\Supp\,C= \Ex(\nu)$. 
We get 
\[
\nu^*(\MO_{X'}(-F')|_S) \simeq \nu^*\omega_S \simeq \MO_{S^N}(K_{S^N} + C) 
\simeq \MO_{S^N}(C).
\]
As $C$ is effective, it holds that   $C=0$, i.e., $S$ is normal. 
Hence $S$ is normal and has at worst canonical singularities. 
By $K_S \sim 0$ and $H^1(S, \MO_S) =0$, $S$ is a geometrically integral canonical K3-like surface. 
This completes the proof of Step \ref{s5-generic-ele}. 
\end{proof}
By Step \ref{s4-generic-ele} and Step \ref{s5-generic-ele}, we obtain $F=0$. 
Then the generic member $S$ of $|-K_X|  = |M|$ is a geometrically integral canonical K3-like surface by Step \ref{s5-generic-ele}. 
It follows from Corollary \ref{c-K3-to-Fano} that $S$ is regular. 
This completes the proof of Theorem \ref{t-generic-ele}. 
\end{proof}

\begin{cor}\label{c-generic-ele}
We work over an algebraically closed field {\cred $k$} of characteristic $p>0$. 
Let $X$ be a Fano threefold. 
Assume that (I) or (II) holds. 
\begin{enumerate}
\item[(I)] $\rho(X) =1$. 
\item[(II)] $\dim ({\rm Im}(\varphi_{|-K_X|})) \geq 2$. 
\end{enumerate}
Then the following hold. 
\begin{enumerate}
\item $H^i(X, nK_X)=0$ for all $i\in \{1, 2\}$ and $n \in \Z$. 
\item 
For all $i>0$ and $m \in \Z_{\geq 0}$, 
it holds that 
$H^i(X, -mK_X) = 0$ and 
\[
h^0(X, -mK_X)=  \frac{1}{12} m(m+1)(2m+1)(-K_X)^3 + 2m + 1. 
\]
In particular, $h^0(X, -K_X) = \frac{1}{2} (-K_X)^3 +3 = g+2$ 
for $g := \frac{1}{2} (-K_X)^3 +1$. 
\end{enumerate}
\end{cor}

\begin{proof}
Since (2) follows from (1) and Corollary \ref{c-kawakami}, 
let us show (1). 
Fix $n \in \Z$. 
Set $k' := K(\P(H^0(X, -K_X)))$ and $X' := X \times_k k'$. 
Let $S$ be the generic member of $|-K_X|$, which is a geometrically integral regular K3-like surface (Theorem \ref{t-generic-ele}). 
{\cred For $\ell \in \Z$,} we have the following exact sequence: 
\[
0 \to \MO_{X'}(nK_{X'}- {\cred \ell} S) \to \MO_{X'}(nK_{X'} 
 -(\ell-1)S) \to \MO_{X'}(nK_{X'} -(\ell-1)S)|_S \to 0. 
\]
By 
{\cred $H^1(S, \MO_{X'}(nK_{X'}-(\ell-1)S)|_S) \simeq 
H^1(S, \MO_{X'}( (n+\ell-1)K_{X'}))=0$} (Theorem \ref{t-K3-vanishing}), 
we have a surjection for any $\ell >0$: 
\[
H^1(X', \MO_{X'}(nK_{X'} - \ell S))  \to H^1(X', \MO_{X'}(nK_{X'})). 
\]
For $\ell \gg 0$, the Serre vanishing theorem implies 
$H^1(X', \MO_{X'}(nK_{X'}))=0$. 
By Serre duality, we get 
$h^2(X', \MO_{X'}(nK_{X'})) = h^1(X', \MO_{X'}( (1-n)K_{X'}))=0$. Thus (1) holds. 
\end{proof}

\section{The case when $|-K_X|$ is not base point free}

\begin{prop}\label{p-non-bpf1}
We work over an algebraically closed field $k$ of characteristic $p>0$. 
Let $X$ be a Fano threefold. 
Assume that $\Bs\,|-K_X| \neq \emptyset$ 
and $\dim ({\rm Im}\,\varphi_{|-K_X|}) \geq 2$. 
Then the following hold. 
\begin{enumerate}
\item The base scheme $Z := \Bs\,|-K_X|$ is isomorphic to $\P^1_k$. 
\item For a general member $T$ of $|-K_X|$, $T$ is a prime divisor 
which is smooth around $Z$. 
\item $K_X \cdot Z = 2- g$. 
\item There is an exact sequence 
\[
0 \to \MO_Z(g-2) \to N_{Z/X} \to \MO_Z(-2) \to 0, 
\]
where $\MO_Z(\ell)$ denotes the invertible sheaf of degree $\ell$ on $Z \simeq \P^1_k$. 
In particular, $\deg N_{Z/X} = g-4$. 
\end{enumerate}
\end{prop}

\begin{proof}
Set $k' := K(\P(H^0(X, -K_X)))$. 
Let $S$ be the generic member of $|-K_X|$. 
By {\cred Theorem \ref{t-generic-ele}}, 
$S$ is a geometrically integral regular K3-like surface on $X' := X \times_k k'$ with $S \sim -K_{X'}$. 

Let us show (1). 
We have an exact sequence 
\[
0 \to \MO_{X'}(-K_{X'}-S) \to \MO_{X'}(-K_{X'}) \to \MO_{X'}(-K_{X'})|_S \to 0. 
\]
By $S \in |-K_{X'}|$ and $H^1(X', \MO_{X'}(-K_{X'}-S)) \simeq H^1(X', \MO_{X'})=0$, 
it holds that $\Bs\,|-K_{X'}| = \Bs\,|-K_{X'}|_S|$. 
We then obtain 
\[
\Bs\,|-K_X| \times_k k' \simeq \Bs\,|-K_{X'}| = \Bs\,|-K_{X'}|_S| \simeq \P^1_{k'}, 
\]
where the last isomorphism follows from Theorem \ref{t-K3-Bs-reg}, 
{\cred which is applicable because $S$ has a $k'$-rational point  \cite[Proposition 5.15(3)]{Tana}}. 
Then $\Bs\,|-K_{X'}|_S|$ is a smooth projective curve over $k'$ with 
\[
H^1(\Bs\,|-K_{X'}|_S|, \MO_{\Bs\,|-K_{X'}|_S|})=0.
\]
These properties descend via the base change $(-) \times_k k'$, and hence 
$\Bs\,|-K_X|$ is a smooth projective curve over $k$ with $H^1(\Bs\,|-K_X|, \MO_{\Bs\,|-K_X|})=0$, 
which implies $\Bs\,|-K_X| \simeq \P^1_k$. 
Thus (1) holds. 

Let us show (2). 
Recall that we have the universal family 
\[
\rho : \mathcal S \hookrightarrow X \times_k \P(H^0(X, -K_X)) \xrightarrow{{\rm pr}_2} 
\P(H^0(X, -K_X))
\]
parameterising the members of $|-K_X|$. 
Here $S$ is the generic fibre of $\rho$ and a general member $T$ is a fibre of $\rho$ over a general closed point of $\P(H^0(X, -K_X))$. 
We have the following two inclusions, each of which is a closed immersion: 
\[
\Bs\,|-K_X| \times_k \P(H^0(X, -K_X)) \subset \mathcal S\subset X \times_k \P(H^0(X, -K_X)). 
\]
After taking the generic fibres, these inclusions become  
\[
\Bs\,|-K_X| \times_k k' \subset S \subset X \times_k k' =X'. 
\]
Note that $\Bs\,|-K_X| \times_k k' (\simeq \P^1_{k'})$ is an effective Cartier divisor on a regular surface $S$. 
Therefore, $\Bs\,|-K_X|$ is an effective Cartier divisor on $T$ for a general member $T \in |-K_X|$. 
Since $\Bs\,|-K_X|$ is smooth by (1), $T$ is smooth around $\Bs\,|-K_X|$. 
Thus (2) holds.

Let us show (3) and (4). 
By applying Theorem \ref{t-K3-Bs-reg} for $L := \MO_{X'}(-K_{X'})|_S$, we get 
\[
(Z\times_k k'\text{ in }S)^2  =-2 \qquad \text{and}\qquad \MO_{X'}(-K_{X'})|_S \cdot (Z \times_k k') = g-2. 
\]
It holds that 
\[
-K_X \cdot Z = \MO_{X'}(-K_{X'}) \cdot (Z \times_k k')  =  \MO_{X'}(-K_{X'})|_S \cdot (Z \times_k k') = g-2. 
\]
Thus (3) holds. 
Since $T$ is smooth around $Z$, the following sequence is exact (cf. \cite[Ch. II, Theorem 8.17]{Har77}): 
\[
0\to N_{Z/T} \to N_{Z/X} \to  N_{T/X}|_Z \to 0. 
\]
We have $(Z\text{ in }T)^2 = (Z\times_k k'\text{ in }S)^2  =-2$ and $T \cdot Z = -K_X \cdot Z = g-2$. 
Therefore, $N_{Z/T} \simeq \MO_Z({\cred -2})$ and $N_{T/X}|_Z \simeq \MO_Z({\cred g-2})$. 
We get 
\[
\deg N_{Z/X} =\deg (\det N_{Z/X}) = \deg N_{Z/X} + \deg (N_{T/X}|_Z)= g-4.
\]
Thus (4) holds. 
\end{proof}

\begin{prop}\label{p-non-bpf2}
We work over an algebraically closed field {\cred $k$} of characeteristic $p>0$. 
Let $X$ be a Fano threefold. 
Assume that $\Bs\,|-K_X| \neq \emptyset$ 
and $\dim ({\rm Im}\,\varphi_{|-K_X|}) \geq 2$. 
Let $\sigma : X' \to X$ be the blowup along $Z := \Bs\,|-K_X| (\simeq \P^1_k)$ (cf. Proposition \ref{p-non-bpf1}) and let 
\begin{equation}\label{e1-non-bpf2}
\varphi_{|-K_{X'}|} : X' \xrightarrow{\psi} Y \hookrightarrow \P^{g+1}
\end{equation}
be the induced morphisms, where $Y := 
\varphi_{|-K_{X'}|}(X')$ and the latter morphism $Y \hookrightarrow \P^{g+1}$ is the induced closed immersion. 
Set $Z' := \Ex(\sigma)$. 
Then the following hold. 
\begin{enumerate}
\item $\dim Y=2$. 
\item There exists a non-empty open subset $Y'$ of $Y$ such that $\chi(X'_y, \MO_{X'_y})=0$ for any point $y \in Y'$. 
\item $\psi|_{Z'} : Z' \to Y$ is birational. Furthermore, either 
\begin{enumerate}
\item $\psi|_{Z'}$ is an isomorphism, or 
\item $\psi|_{Z'}$ is the contraction of the curve $\Gamma$ on $Z'$ with $\Gamma^2 <0$.
\end{enumerate}
\item {\cred $\psi_*\MO_{X'} = \MO_Y$}. 
\item {\cred There exists a non-empty open subset $Y''$ of $Y$ such that 
the scheme-theoretic fibre $X'_y$ is geometrically integral for every point $y \in Y''$.} 
\end{enumerate}
\end{prop}

\begin{proof}
By
\[
K_{X'} = \sigma^* K_X + Z' \qquad \text{and} \qquad  \sigma^*(-K_X) =-K_{X'} + Z', 
\]
$-K_{X'}$ coincides with the mobile part of $|\sigma^*(-K_X)|$, 
which is base point free by construction. 
Hence we obtain the induced morphisms (\ref{e1-non-bpf2}). 

Let us show (1). 
By $Y = {\rm Im}\,\varphi_{|-K_X|}$ and the assumption $\dim ({\rm Im}\,\varphi_{|-K_X|}) \geq 2$, we have $\dim Y \geq 2$. 
In order to show $\dim Y =2$, it suffices to show $(-K_{X'})^3 =0$. 
For a fibre $F'$ of $Z' \to Z$ over a closed point of $Z$, 
we have $\sigma^*K_X \cdot Z' = (K_X \cdot Z) F'$, $Z' \cdot F' = -1$, and $Z'^3 = -\deg_Z(N_{Z/X})$ \cite[Lemma 2.11]{Isk77}, 
which implies the following. 
\begin{itemize}
\item $K_X^3  = 2-2g$. 
\item $(\sigma^*K_X)^2 \cdot Z' = 0$. 
\item $\sigma^*K_X \cdot Z'^2 = (K_X \cdot Z) F' \cdot Z' = - K_X \cdot Z=g-2$  (Proposition \ref{p-non-bpf1}). 
\item $Z'^3 = -\deg_Z(N_{Z/X}) = 4-g$  (Proposition \ref{p-non-bpf1}). 
\end{itemize}
It holds that 
\[
K_{X'}^3 = (\sigma^*K_X + Z')^3 =K_X^3 + 3 (\sigma^*K_X)^2 \cdot Z' + 3 \sigma^*K_X \cdot Z'^2 + 
Z'^3
\]
\[
=(2-2g) + 0 +3(g-2) +(4-g) =0. 
\]
Thus (1) holds. 

Let us show (2). 
By the generic flatness, it suffices to show that $\chi(X'_K, \MO_{X'_K})=0$, 
where $K :=K(Y)$ and $X'_K$ denotes the generic fibre of $\psi : X' \to Y$. 
Note that $X'_K$ is a regular projective curve over $K$. 
Thus it is enough to check $\omega_{X'_K/K} \simeq \MO_{X'_K}$. 
Note that $\omega_{X'}$ is a pullback of an invertible sheaf on $Y$. 
For some non-empty smooth open subset $Y_1$ of $Y$ and its inverse image $X'_1 := \psi^{-1}(Y_1)$, we get 
\[
\MO_{X'_K} \simeq \omega_{X'_1/Y_1}|_{X'_K} \simeq \omega_{X'_K/K}, 
\]
where the latter isomorphism follows from \cite[Theorem 3.6.1]{Con00}. 
Thus (2) holds. 

Let us show (3). 
The induced moprhism $\psi|_{Z'} : Z' \to Y$ is surjective, because 
\begin{eqnarray*}
(-K_{X'})^2 \cdot Z' 
&=& (\sigma^*K_{X} +Z')^2 \cdot Z'\\
&=& (\sigma^*K_{X})^2 \cdot Z' + 2 \sigma^*K_X \cdot Z'^2 + Z'^3\\
&=& 0 + 2 (g-2) +(4-g) \\
&=& g  = \frac{(-K_X)^3}{2} +1\neq 0. 
\end{eqnarray*}
For $\ell := \deg (\psi|_{Z'} : Z' \to Y)$, the following holds: 
\begin{align*}
0 &\leq \Delta(Y, \MO_{\P^{g+1}(1)}|_Y) \\
&= \dim Y + (\MO_{\P^{g+1}}(1)|_Y)^2 -h^0(Y, \MO_{\P^{g+1}}(1)|_Y)\\
&= 2+ \frac{(-K_{X'})^2 \cdot Z'}{\ell} -h^0(X', -K_{X'})\\
&= 2+ \frac{g}{\ell} -(g+2)\\
&= \frac{g(1-\ell)}{\ell}, 
\end{align*}
where the second equality follows from Lemma \ref{l-H^0-image}. 
Therefore, $\ell =1$ and $\Delta(Y, \MO_{\P^{g+1}(1)}|_Y) =0$. 
In particular, $\psi|_{Z'} : Z' \to Y$ is birational and $Y$ is normal. 
Thus (3) holds. 

{\cred 
Let us show (4). 
Let $\psi : X'  \xrightarrow{\widetilde \psi} \widetilde Y \xrightarrow{\theta} Y$ be the Stein factorisation of $\psi$. 
Since the composite morphism $\psi|_{Z'} : Z' \xrightarrow{\widetilde \psi|_{Z'}} \widetilde Y \xrightarrow{\theta} Y$ is birational by (3), 
both $\widetilde \psi|_{Z'}$ and $\theta$ are birational. 
Then $\theta : \widetilde Y \to Y$ is a finite birational morphism of normal projective surfaces by (3), and hence $\theta$ is an isomorphism. 
Thus (4) holds.

Let us show (5). 
Fix a general closed point $y \in Y$. 
It is enough to show that the fibre $X'_y$ of $\psi : X' \to Y$ over $y$ is an integral scheme. 
By (4), the generic fibre of $\psi : X' \to Y$ is geometrically irreducible \cite[Lemma 2.2]{Tan18-b}, 
and hence so is a general fibre $X'_y$ \cite[Proposition 9.7.8]{EGAIV3}. It is clear that $\dim X'_y =1$. 
Since $\psi|_{Z'} : Z' \to Y$ is birational, we get $Z' \cdot X'_y = 1$. 
Hence $X'_y$ generically reduced \cite[Lemma 1.18]{Bad01}, i.e., 
the local ring $\MO_{X'_y, \gamma}$ is a field for the generic point $\gamma$ of $X'_y$. 
 Since $X'_y$ is CM, $X'_y$ is $S_1$ and $R_0$, i.e., $X'_y$ is reduced. 
Hence $X'_y$ is integral. 
Thus (5) holds. 
}
\end{proof}

\begin{thm}\label{t-non-bpf-main}
We work over an algebraically closed field {\cred $k$} of characeteristic $p>0$. 
Let $X$ be a Fano threefold with $\rho(X)=1$. 
Then $|-K_X|$ is base point free. 
\end{thm}

\setcounter{step}{0}
\begin{proof}
Suppose that $\Bs\,|-K_X| \neq \emptyset$. 
Let us derive a contradiction. 
By $\rho(X) =1$, we have $\dim ({\rm Im}\,\varphi_{|-K_X|}) \geq 2$ (Corollary \ref{c-composite-pencil}). 
Then we may apply Proposition \ref{p-non-bpf2} and 
we shall use the same notation as in the statement of Proposition \ref{p-non-bpf2}. 
We have $Z' \simeq \F_m = \P_{\P^1}( \MO_{\P^1} \oplus \MO_{\P^1}(m))$ for some $m \geq 0$. 
By $\rho(X) =1$, we get $\rho(X')=2$ and $\rho(Y)=1$. 
It follows from Proposition \ref{p-non-bpf2} that $Y \simeq \P(1, 1, m)$, 
which is nothing but the birational contraction 
of the curve $\Gamma \subset {\cred Z'}$ with $\Gamma^2 <0$, where $m := -\Gamma^2 >0$. 
Set $y_0 := \psi(\Gamma) \in Y$. 

For a general fibre $B_{Z'}$ of $\sigma|_{Z'} : Z'= \F_m \to \P^1 =Z$, 
we set $B := \psi(B_{Z'}) \simeq \P^1_k$. 
Then 
a general fibre of $\psi^{-1}(B) \to B$ is one-dimensional and {\cred integral  (Proposition \ref{p-non-bpf2})}. 
Furthermore, $\psi^{-1}(B^{\circ}) \to B^{\circ}$ is flat 
for some non-empty open subset $B^{\circ} \subset B$. 
Hence $\psi^{-1}(B^{\circ})$ is {\cred integral} and two-dimensional. 
Then its closure $V := \overline{\psi^{-1}(B^{\circ})}$, equipped with the reduced scheme structure, is a prime divisor on $X'$. 
We have $V \cdot X'_y =0$ for a general closed point $y \in Y$. 

\begin{claim*}
$V = \psi^{-1}(B)_{\red}$. 
\end{claim*}

\begin{proof}[Proof of Claim]
The inclusion $V \subset \psi^{-1}(B)_{\red}$ is clear. 
Suppose that $V \not\supset \psi^{-1}(B)_{\red}$. 
Then it holds that $\psi^{-1}(b)_{\red} \not\subset V$ for some $b \in B$. 
Since $\psi^{-1}(b)$ is connected  {\cred (Proposition \ref{p-non-bpf2})} and $\psi^{-1}(b) \cap V \neq \emptyset$, 
we can find a curve $C \subset \psi^{-1}(b)$ such that 
$V \cap C \neq \emptyset$ and $C \not\subset V$. 
In particular, $V \cdot C >0$. 
In order to derive a contradiction, let us prove that $\rho(X') \geq 3$. 
It is enough to show that $H_{X'}, \psi^*H_Y, V$ are $\Z$-linearly independent  in $\Pic\,X'/{\equiv}$, 
where $H_{X'}$ and $H_Y$ are ample Cartier divisors on $X'$ and $Y$, respectively. 
For $\alpha, \beta, \gamma \in \Z$, assume that 
\[
\alpha H_{X'} + \beta\psi^*H_{Y} + \gamma  V \equiv 0. 
\]
By taking the intersection with a general fibre of $\psi$, we obtain $\alpha =0$. 
Then consider the intersection with the above curve $C$, which deduces $\gamma =0$ and hence $\beta=0$. 
This completes the proof of Claim. 
\end{proof}

Let $\psi_V : V \to B$ be the morphism induced by $\psi : X' \to Y$. 
Note that $\psi_V$ is a flat morphism from a projective surface $V$ to $B \simeq \P^1$. 
{\cred Since $\psi^{-1}(B^{\circ})$ is an integral scheme, 
$\psi_V : V \to B$ is a compatification of $\psi^{-1}(B^{\circ}) \to B^{\circ}$.} 
Therefore, we obtain 
\begin{equation}\label{e1-t-non-bpf-main}
\chi(V_y, \MO_{V_y}) = \chi(X'_y, \MO_{X'_y}) =0, 
\end{equation} 
\begin{equation}\label{e2-t-non-bpf-main}
\chi( V_{y_0}, \MO_{V_{y_0}}) = \chi( V_{y}, \MO_{V_{y}}),\quad \text{and}\quad 
\MO_{X'}(Z')|_V \cdot V_{y_0} = \MO_{X'}(Z')|_V \cdot V_{y} 
\end{equation}
for a general closed point $y \in B$  (Proposition \ref{p-non-bpf2}).  
It holds that 
\begin{equation}\label{e3-t-non-bpf-main}
\MO_{X'}(Z')|_V \cdot V_{y} 
= Z' \cdot X'_y =  1. 
\end{equation}
On the other hand, $V_{y_0} = \psi_V^{-1}(y_0)$ is an effective Cartier divisor on $V$, 
and hence pure one-dimensional. 
Let $(V_{y_0})_{\red} = \bigcup_{i=1}^s \Gamma_i$ be the irreducible decomposition with $\Gamma_1 :=\Gamma$. 
We then obtain 
\begin{equation}\label{e4-t-non-bpf-main}
\MO_{X'}(Z')|_V \cdot V_{y_0} = \MO_{X'}(Z')|_V \cdot \left( \sum_{i=1}^s \ell_i \Gamma_i \right) = \sum_{i=1}^s \ell_i \MO_{X'}(Z') \cdot \Gamma_i,  
\end{equation}
where  $\ell_i :={\rm length}_{\MO_{V_{y_0}, \gamma_i}} \MO_{V_{y_0}} >0$ 
for the generic point $\gamma_i$ of $\Gamma_i$ (cf. \cite[Lemma 1.18]{Bad01}). 
It follows from (\ref{e2-t-non-bpf-main})--(\ref{e4-t-non-bpf-main}) that 
\[
\sum_{i=1}^s \ell_i \MO_{X'}(Z') \cdot \Gamma_i = 1. 
\]

We now show that $Z' \cdot \Gamma_i >0$ for any $i$. 
Otherwise, there exists $i$ such that $Z' \cdot \Gamma_i \leq 0$. 
Fix ample Cartier divisors $H_{X'}$ and $H_Y$  on $X'$ and $Y$,  respectively. 
It suffices to show that $H_{X'}, \psi^*H_Y, Z'$ are $\Z$-linearly independent in $\Pic\,X'/\equiv$, as this implies $\rho(X') \geq 3$. 
Assume 
\[
\alpha H_{X'} + \beta \psi^*H_Y + \gamma  Z' \equiv 0
\]
for some $\alpha, \beta, \gamma \in \Z$. 
Taking the intersections with $\Gamma_i$ and a general fibre of $\psi$, 
we obtain $\alpha = \gamma =0$ and hence $\beta =0$. 
Therefore, $Z' \cdot \Gamma_i >0$ for any $i$.

We then get $s=1$ and $\ell_1 =1$. 
In other words, $(V_{y_0})_{\red} = \Gamma$ and $V_{y_0}$ is generically reduced, 
i.e., the local ring $\MO_{V_{y_0}, \gamma_1}$ is a field. 
{\cred 
Since $V$ is a prime divisor on a smooth threefold $X'$, $V$ is CM, and hence so is $V_{y_0}$. 
Therefore, $V_{y_0}$ is reduced, i.e., $V_{y_0}= \Gamma$. 
}
Hence
\[
0 = \chi(V_y, \MO_{V_y}) = \chi(V_{y_0}, \MO_{V_{y_0}}) = 
{\cred h^0({\cred \Gamma}, \MO_{{\cred \Gamma}}) - h^1({\cred \Gamma}, \MO_{{\cred \Gamma}})
=h^0({\cred \Gamma}, \MO_{{\cred \Gamma}})
= 1,}
\]
where the first equality follows from (\ref{e1-t-non-bpf-main}) and 
the second one holds by (\ref{e2-t-non-bpf-main}). 
This is a contradiction. 
\end{proof}

\section{The case when $|-K_X|$ is base point free}\label{s-bpf-case}

Throughout this section, we work over an algebraically closed field $k$ of characteristic $p>0$. 

\subsection{Birational case}

\begin{lem}\label{l-birat-or-hyperell}
Let $X$ be a Fano threefold 
such that $|-K_X|$ is base point free. 
For $Y := \Im (\varphi_{|-K_X|} : X \to \P^{h^0(X, -K_X)-1})$, 
let $\psi: X \to Y$ be the morphism induced by $\varphi_{|-K_X|}$. 
Then $\deg \psi =1$ or $\deg \psi =2$. 
\end{lem}


\begin{proof}
By construction, we can find an ample Cartier divisor $H_Y$ on $Y$ such that $\psi^*H_Y \sim -K_X$ and $\psi^* : H^0(Y, H_Y) \xrightarrow{\simeq} H^0(X, -K_X)$ 
(Lemma \ref{l-H^0-image}). 
It holds that 
\[
0 \leq \Delta(Y, H_Y) = \dim Y +H_Y^3 -h^0(Y, H_Y) 
=3 + \frac{(-K_X)^3}{ \deg \psi} -h^0(X, -K_X). 
\]
By $(-K_X)^3 = 2g -2$ and $h^0(X, -K_X)= g+2$ (Corollary \ref{c-generic-ele}), we obtain  
\[
g+2 = h^0(X, -K_X)\leq 
3+ \frac{(-K_X)^3}{\deg \psi} = 3+ \frac{2g -2}{\deg \psi}, 
\]
which implies $0 \leq (g-1)(2 - \deg \psi)$. 
Hence we obtain $\deg \psi =1$ or $\deg \psi =2$.  
\end{proof}

\begin{thm}\label{t-Fano-birat}
Let $X$ be a Fano threefold such that $|-K_X|$ is base point free and 
$\varphi_{|-K_X|} : X \to \P^{h^0(X, -K_X) -1}$ is birational onto its image. 
Then 
\[
\bigoplus_{m=0}^{\infty} H^0(X, \MO_X(-mK_X))
\]
is generated by $H^0(X, \MO_X(-K_X))$ as a $k$-algebra. 
In particular, $|-K_X|$ is very ample. 
\end{thm}

\begin{proof}
Set $Y := \varphi_{|-K_X|}(X)$. We have the following induced morphisms: 
\[
\varphi_{|-K_X|} : X \xrightarrow{\psi} Y \hookrightarrow  \P^{h^0(X, -K_X) -1}. 
\]
Fix general hyperplane secionts $H_Y, H'_Y \subset Y$. 
Set 
\[
H_X := \psi^{-1}(H_Y), \quad H'_X := \psi^{-1}(H'_Y),\quad C_Y := H_Y \cap H'_Y, \quad 
C_X := H_X \cap H'_X. 
\]
Note that all $H_X, H'_X, C_X, H_Y, H'_Y, C_Y$ are integral schemes (Proposition \ref{l-Bertini-birat}).
By applying Proposition \ref{p-wt-lift} twice, 
it suffices to show that 
\[
R(C_X, \MO_X(-K_X)|_{C_X})  
\]
is generated by $H^0(C_X, \MO_X(-K_X)|_{C_X})$. 
It follows from the adujction formula that $C_X$ is a projective Gorenstein curve with 
\[
\MO_X(-K_X)|_{C_X} \simeq \MO_X( K_X -2K_X)|_{C_X} \simeq 
(\omega_X \otimes \MO_X(H_X) \otimes \MO_X(H'_X))|_{C_X} \simeq \omega_{C_X}. 
\]
Since $|-K_X|$ is base point free, also $|-K_X|_{C_X}| = |\omega_{C_X}|$ is base point free. 
By construction, $\varphi_{|\omega_{C_X}|}$ is birational onto its image $C_Y$. 
Therefore, $R(C, \omega_C)$ is generated by $H^0(C, \omega_C)$ 
\cite[Theorem (A1) in page 39]{Fuj83}. 
\end{proof}

\subsection{Hyperelliptic case}

\begin{dfn}\label{d-hyperell}
We say that a Fano threefold $X$ is {\em  hyperelliptic} 
if $X$ is of index one, $|-K_X|$ is base point free, and 
the morphism 
$\varphi : X \to \varphi_{|-K_X|}(X)$, induced by $\varphi_{|-K_X|} : X \to \P^{h^0(X, -K_X)-1}$, 
is of degree two. 
\end{dfn}

\begin{nota}\label{n-hyperell}
Let $X$ be a hyperelliptic Fano threefold and 
let $\varphi : X \to Y := \varphi_{|-K_X|}(X) \subset \P^{g+1}$ be the morphism induced by $\varphi_{|-K_X|}$, where $g := h^0(X, -K_X)-2$. 
Since $-K_X$ is ample, $\varphi$ is a finite surjective morphism of projective threefolds 
such that $[K(X):K(Y)]=2$. 
Set $\MO_Y(\ell) := \MO_{\P^{g+1}}(\ell)|_Y$ for any $\ell \in \Z$. 
\end{nota}

\begin{thm}\label{t-hyperell-classify}
We use Notation \ref{n-hyperell}. 
Then the following hold. 
\begin{enumerate}
\item $\Delta(Y, \MO_Y(1)) =0$. 
\item $Y$ is smooth. 
\item One of the following holds. 
\begin{enumerate}
\renewcommand{\labelenumii}{(\roman{enumii})}
\item $Y \simeq \P^3$. 
\item $Y$ is a smooth quadric hypersurface in $\P^4$. 
\item $\rho(X) \geq 2$. 
\end{enumerate}
\end{enumerate}
\end{thm}


\begin{proof}
Let us show (1). 
For the double cover $\varphi : X \to Y$, 
it holds that 
\[
2g-2 = (-K_X)^3  = (\deg \varphi) \cdot (\deg Y) = 2 \cdot \deg Y,
\]
which implies $\deg Y = g-1$. 
Therefore, we get 
\[
\Delta(Y, \MO_Y(1)) = \dim Y + \MO_Y(1)^3 - h^0(Y, \MO_Y(1)) 
= 3 + \deg Y - (g+2) = 0, 
\]
where the second equality follows from $h^0(Y, \MO_Y(1)) = h^0(X, -K_X) = g+2$ 
(Lemma \ref{l-H^0-image}, Corollary \ref{c-generic-ele}). 
Thus (1) holds. 

We now show that (1) and (2) imply (3). 
Assume that none of (i) and (ii) holds. 
By (1), (2), and Theorem \ref{t-Delta0-sm}, we have $Y =\P_{\P^1}(E)$ for some vector bundle $E$ on $\P^1$ of rank $3$. 
We then obtain $\rho(X) \geq \rho(Y) =2$, and hence (iii) holds. 
Thus (1) and (2) imply (3).

It is enough to prove (2). 
Suppose that $Y$ is singular. 
Let us derive a contradiction. 
By Remark \ref{r-Delta0-sing}, there are the following two cases. 
\begin{enumerate}
\item[(a)] $Y$ is a cone over a smooth polarised surface $(Z, \MO_Z(1))$ 
with $\Delta(Z, \MO_Z(1))=0$. 
\item[(b)] $Y$ is a cone over a smooth polarised curve $(Z, \MO_Z(1))$  
with $\Delta(Z, \MO_Z(1))=0$. 
\end{enumerate}
Fix Cartier divisors $H_Y$ and $H_Z$ on $Y$ and $Z$ respectively such that $\MO_Y(H_Y) \simeq \MO_Y(1)$ and $\MO_Z(H_Z) \simeq \MO_Z(1)$.

Assume (a). By Remark \ref{r-Delta0-sing}, we have 
\[
\begin{CD}
W=\P_Z(\MO_Z \oplus \MO_Z(1))  @>\mu >> Y\\
@VV\pi V\\
Z, 
\end{CD}
\]
where $\pi$ is the induced projection and $\mu$ is the birational morphism such that 
$v := \mu(\widetilde{Z})$ is the vertex of $Y$ for the section $\widetilde Z \subset W$ of $\pi$ corresponding to the surjection: $\MO_Z \oplus \MO_Z(1) \to \MO_Z, \,\, (a, b) \mapsto a$. 
By Theorem \ref{t-Delta0-sm}, 
 $(Z, \MO_Z(1)) = (\P^2, \MO_{\P^2}(2))$ or $Z = \F_n$ for some $n \geq 0$. 
 For the latter case: $Z = \F_n$, 
 $Y$ is not $\Q$-factorial (Proposition \ref{p-RNS-Qfac}). 
However, this contradicts the fact that $Y$ is the image of a finite morphsim $\varphi : X \to Y$ from a smooth variety $X$ \cite[Lemma 5.16]{KM98}. 


We may assume that $(Z, \MO_Z(1)) = (\P^2, \MO_{\P^2}(2))$. 
We then have  $H_Z \sim 2 H'_Z$ for some Cartier divisor $H'_Z$ on $Z$, 
which implies 
\[
H_Y \sim \mu_*\pi^*H_Z \sim 2 \mu_* \pi^*H'_Z = 2 H'_Y, 
\]
where $H'_Y := \mu_*\pi^*H'_Z$ is a Weil divisor, which is not necessarily Cartier. 
Therefore, we get 
\[
-K_X = \varphi^*H_Y \sim 2 \varphi^*H'_Y, 
\]
where the linear equivalence can be checked after removing $v$ and $\varphi^{-1}(v)$. 
This contradicts the assumption that $X$ is of index $1$. 
This completes the case (a). 

Assume (b). By Remark \ref{r-Delta0-sing}, we have $(Z, \MO_Z(1)) = (\P^1, \MO_{\P^1}(m))$ with $m >0$ and $Y \subset \P^{m+2}$ is obtained from $Z \subset \P^m$ by applying cone construction:  
\[
\begin{CD}
W=\P_Z(\MO_Z \oplus \MO_Z \oplus \MO_Z(1))  @>\mu >> Y\\
@VV\pi V\\
Z. 
\end{CD}
\]
Furthermore, we get $m \geq 2$, as otherwise the cone 
$Y$ would be smooth, 
i.e., $Y \simeq \P^3$. 
Hence it holds that 
$H_Y \sim \mu_*\pi^*(mP)$  for a closed point $P \in \P^1$. 
This contradicts 
the assumption that $X$ is of index $1$. 
\end{proof}

\begin{prop}\label{p-hyperell-P^3}
We use Notation \ref{n-hyperell}. 
Assume $Y = \P^3$. 
Then the following hold. 
\begin{enumerate}
    \item $\MO_X(-K_X) \simeq \varphi^*\MO_{\P^3}(1)$. 
    \item $(-K_X)^3 = 2$. 
    \item If $p \neq 2$, then $\varphi : X \to Y = \P^3$ is a double cover ramified 
along a smooth prime divisor $S \subset \P^3$ of degree $6$. 
\end{enumerate}
\end{prop}

\begin{proof}
Since $X$ is of index one, we get (1). 
It is clear that (1) implies (2). 
For the branch divisor $S \subset \P^3$, we have 
\[
K_X = \psi^*\left( K_Y + \frac{1}{2}S \right), 
\]
which deduces (3). 
\end{proof}

\begin{prop}\label{p-hyperell-Q}
We use Notation \ref{n-hyperell}. 
Assume that $Y$ is a smooth quadric hypersurface in $\P^4$.  
Then the following hold. 
\begin{enumerate}
    \item $\MO_X(-K_X) \simeq \varphi^*(\MO_{\P^4}(1)|_Y)$. 
    \item $(-K_X)^3 = 4$. 
    \item If $p \neq 2$, then $\varphi : X \to Y$ is a double cover ramified 
along a smooth prime divisor $S \subset Y$ which is a complete intersection of 
$Y$ and a hypersurface of degree $4$. 
\end{enumerate}
\end{prop}

\begin{proof}
The same argument as in Proposition \ref{p-hyperell-P^3} works. 
\end{proof}

Although the following proposition will not be used in the rest of this paper, we shall later need it in order to classify Fano threefolds with $\rho=2$. 

\begin{prop}\label{p-g=2-Pic1}
Let $X$ be a Fano threefold such that $(-K_X)^3 =2$. 
Assume that $\dim (\Im\,\varphi_{|-K_X|}) \geq 2$. 
Then the following hold. 
\begin{enumerate}
\item 
$X$ is isomorphic to a weighted hypersurface in $\P(1, 1, 1, 1, 3)$ of degree $6$. 
\item $\rho(X) = r_X=1$, where $r_X$ denotes the index of $X$. 
\item 
$h^0(X, -K_X) =4$ and the induced morphism $\varphi_{|-K_X|} : X \to \P^3$ is a double cover (cf. Proposition \ref{p-hyperell-P^3}). 
\end{enumerate}
\end{prop}

\begin{proof}
First of all, we show that $|-K_X|$ is base point free. 
Suppose that $|-K_X|$ is not base point free. 
Let us  derive a contradiction. 
By $\dim (\Im\,\varphi_{|-K_X|}) \geq 2$, 
the generic member $S$ of $|-K_X|$ is a geometrically integral regular K3-like surface 
(Theorem \ref{t-generic-ele}). 
Set $k' := K(\P ( H^0(X, -K_X)))$ and $X' := X \times_k k'$.  Then 
the restriction map 
\[
H^0(X', -K_{X'}) \to H^0(S, -K_{X'}|_S) 
\]
is surjective, because $H^1(X', -K_{X'} -S) \simeq H^1(X', \MO_{X'}) \simeq H^1(X, \MO_X) \otimes_k k' =0$. 
Hence $|-K_{X'}|_S|$ is not base point free. 
It follows from Theorem \ref{t-K3-Bs-reg}(5) that 
the base locus $F$ of $|-K_{X'}|_S|$ is a prime divisor on $S$ satisfying $(-K_{X'}|_S) \cdot F = g-2$ for $g := \frac{(-K_{X'}|_S)^2}{2}+1$. 
We have $g =  \frac{(-K_{X'}|_S)^2}{2}+1 = \frac{(-K_X)^3}{2}+1 = 2$, which implies $(-K_{X'}|_S) \cdot F = g-2 =0$. This is absurd, because $(-K_{X'}|_S)$ is ample. 
This completes the proof of the base point freeness of $|-K_X|$. 

\medskip 

By assuming (1), we now finish the proof. 
It is known that (1) implies $\rho(X)=1$ \cite[Theorem 3.7]{Mor75}. 
We have $r_X=1$, because $r_X \geq 2$ implies $(-K_X)^3 \geq 2^3 =8$. 
By $h^0(X, -K_X) = \frac{1}{2}(-K_X)^3 +3 = 4$ (Corollary \ref{c-generic-ele}), we obtain a finite surjective morphism $\varphi_{|-K_X|} : X \to \P^3$. 
By $r_X =1 \neq 4$, $\varphi_{|-K_X|}$ is a double cover (Lemma \ref{l-birat-or-hyperell}), and hence $X$ is hyperelliptic (Definition \ref{d-hyperell}). 

\medskip 

It is enough to show (1). 
Set $\MO_X(\ell) := \MO_X(-\ell K_X)$ for $\ell \in \Z$. 
Note that $\MO_X(-K_X)$ is $3$-regular with respect to an ample globally generated invertible sheaf $\MO_X(-K_X)$ \cite[Section 1.8, especially Example 1.8.24]{Laz04}, \cite[Subsection 5.2]{FGI05}, i.e., 
the following holds (Corollary \ref{c-generic-ele}): 
\[
H^i(X, \MO_X(-(3 -i)K_X)) =0 \qquad \text{for}\qquad i>0. 
\]
Then the induced $k$-linear map 
\[
H^0(X, -K_X) \otimes_k H^0(X, -rK_X) \to H^0(X, -(r+1)K_X)
\]
is surjective for every $r \geq 3$ \cite[Lemma 5.1(a)]{FGI05}. 
Then, as a $k$-algebra, $\bigoplus_{m=0}^{\infty}H^0(X, -mK_X)$ is generated by 
\[
H^0(X, -K_X) \oplus 
H^0(X, -2K_X) \oplus 
H^0(X, -3K_X).  
\]
It follows from Corollary \ref{c-generic-ele} that 
\[
h^0(X, -mK_X) = \frac{1}{12}m(m+1)(2m+1)(-K_X)^3 + 2m +1 
\]
\[
= \frac{1}{6}m(m+1)(2m+1) + 2m +1, 
\]
which implies the following: 
\begin{itemize}
\item $h^0(X, -K_X)= \frac{1}{6}\cdot 1\cdot 2 \cdot 3 + 3 = 4$. 
\item $h^0(X, -2K_X)= \frac{1}{6}\cdot 2\cdot 3 \cdot 5 + 5 = 10$. 
\item $h^0(X, -3K_X)= \frac{1}{6}\cdot 3\cdot 4 \cdot 7 + 7 = 21$. 
\item $h^0(X, -4K_X)= \frac{1}{6}\cdot 4\cdot 5 \cdot 9 + 9 = 39$. 
\item $h^0(X, -5K_X)= \frac{1}{6}\cdot 5\cdot 6 \cdot 11 + 11 = 66$. 
\item $h^0(X, -6K_X)= \frac{1}{6}\cdot 6\cdot 7 \cdot 13 + 13 = 104$. 
\end{itemize}
Fix a $k$-linear basis: $H^0(X, -K_X) = \bigoplus_{0 \leq i \leq 3} kx_i$. 
Note that the subset  
$\{ x_ix_j \,|\, 0 \leq i \leq j \leq 3\}$ 
(resp. $\{ x_ix_jx_k \,|\, 0 \leq i \leq  j\leq k \leq 3\}$) 
of $H^0(X, -2K_X)$ (resp. $H^0(X, -3K_X)$) is linearly independent over $k$, 
as otherwise the induced morphism $\varphi_{|-K_X|} : X \to \P^3$ would factor through 
a hypersurface $Z \subset \P^3$ of degree $2$ (resp. $3$) which is defined by the corresponding linear dependence equation. 
Since the subspace generated by these elements is of dimension $10$ (resp. $20$), 
we can find an element $x_4 \in H^0(X, -3K_X)$ such that 
$H^0(X, - {\cred 3} K_X) = S^3(H^0(X, -K_X)) \oplus kx_4$. 
Let $k[y_0, y_1, y_2, y_3, y_4]$ be the polynomial ring with $\deg y_0 = \deg y_1 = \deg y_2 =\deg y_3 = 1$ and $\deg y_4 =3$. 
We obtain a surjective graded $k$-algebra homomorphism: 
\[
k[y_0, y_1, y_2, y_3, y_4] \to \bigoplus_{m=0}^{\infty} H^0(X, -mK_X), \qquad y_i \mapsto x_i,  
\]
which induces a closed immersion $X \hookrightarrow \P(1, 1, 1, 1, 3) = \Proj\,k[y_0, y_1, y_2, y_3, y_4]$. 
Hence $X$ is a weighted hypersurface in $\P(1, 1, 1, 1, 3)$.

For $d \in \Z_{>0}$, $k[y_0, y_1, y_2, y_3, y_4]_d$ 
denotes the $k$-linear supspace of $k[y_0, y_1, y_2, y_3, y_4]$ consisting of all the homogeneous elements of degree $d$. 
We now see how to compute $\dim_k k[y_0, y_1, y_2, y_3, y_4]_6$. 
Pick a monomial $y_0^{d_0}y_1^{d_1}y_2^{d_2}y_3^{d_3}y_4^{d_4}$ with 
$d_0+d_1 + d_2+d_3 + {\cred 3} d_4 =6$.  
If $d_4 =2$, then there is one solution. 
If $d_4=1$, then the number of the solutions of $\sum_{i=0}^3 d_i =3$ is $\binom{6}{3}=20$. 
If $d_4=0$, then the number of the solutions of $\sum_{i=0}^3 d_i =6$ is $\binom{9}{3}=84$. 
Hence $\dim_k k[y_0, ..., y_4]_6 = 1 + 20 + 84 =105$. 
Similarly, $\dim_k k[y_0, ..., y_4]_4 = \binom{4}{3}+\binom{7}{3} = 39$ 
and $\dim_k k[y_0, ..., y_4]_5 = \binom{5}{3}+\binom{8}{3} = 66$. 
Therefore, we obtain $\deg X = 6$. Thus (1) holds. 
\end{proof}

\section{Intersection of quadrics}

{\cred Throughout this section, we work over an algebraically closed field $k$ of characteristic $p>0$.}


\subsection{Anti-canonically embedded Fano threefolds} 

\begin{dfn}\label{d-anti-emb}
We say that $X \subset \P^{g+1}$ is an {\em anti-canonically embedded Fano threefold} 
if $X$ is a smooth projective threefold such that 
$-K_X$ is very ample, $X$ is a closed subscheme of $\P^{g+1}$, and 
the induced closed immersion $X \hookrightarrow \P^{g+1}$ is given by the complete linear system $|-K_X|$. 
In this case, 
we have that $h^0(X, -K_X) = g+2$ and $\deg X = (-K_X)^3 = 2g -2$ 
(Corollary \ref{c-generic-ele}). 
\end{dfn}


\begin{rem}\label{r-anti-emb}
Let $X \subset \P^{g+1}$ be an anti-canonically embedded Fano threefold. 
Then a general hyperplane section $X \cap \P^g$ is a smooth K3 surface and 
the intersection $X \cap \P^{g-1}$ with a general $(g-1)$-dimensional  linear subvariety {\cred $\P^{g-1}$} is a canonical curve of genes $g$ (Definition \ref{d-cano-curve}). 
In particular, $g \geq 3$. 
\end{rem}

\begin{dfn}\label{d-cano-curve}
We say that $C \subset \P^{g-1}$ is a {\em canonical curve} 
if $C$ is a smooth projective curve of genus $g$ such that 
$|K_C|$ is very ample, $C$ is a closed subscheme of $\P^{g-1}$, and 
the induced closed immersion 
$C \hookrightarrow \P^{g-1}$ is given by the complete linear system $|K_C|$. 
\end{dfn}

\subsection{Trigonal case}


\begin{nota}\label{n-trigonal}
Let $X \subset \P^{g+1}$ be an anti-canonically embedded Fano threefold 
with $g \geq 5$. 
Assume that $X$ is not an intersection of quadrics. 
Set 
\[
W  := \bigcap_{\substack{X \subset Q,\\ Q:\text{quadric}}} Q,   
\]
where the right hand side denotes the scheme-theoretic intersection of all the quadric hypersurfaces of $\P^{g+1}$ containing $X$. 
\end{nota}

\begin{thm}[The Noether--Enriques--Petri theorem]\label{t-NEP}
Let $C \subset \P^{g-1}$ be a canonical curve. 
Assume $g \geq 5$ and $C$ is not an intersection of quadrics. 
Set 
\[
W_C := \bigcap_{\substack{C \subset Q,\\ Q:\text{quadric}}} Q 
\]
to be the scheme-theoretic intersection of all the quadric hypersurfaces in $\P^{g-1}$ containing $C$. 
Then the following hold. 
\begin{enumerate}
\item $W_C$ is a smooth projective surface with $\deg W_C = g-2$. 
\item One of the following holds. 
\begin{enumerate}
    \item $W_C$ is a rational scroll.
    \item $g=6$ and $W_C \subset \P^5$ is a Veronese surface.  
\end{enumerate}  
\end{enumerate}
\end{thm}

\begin{proof}
By \cite[Theorem 4.7 and Lemma 4.8]{SD73}, $(W_C)_{\red}$ is irreducible, two-dimensional, and $\deg (W_C)_{\red} = g-2$. 
Then $(W_C)_{\red}$ is a variety of minimal degree, and hence 
$(W_C)_{\red}$ is an intersection of quadrics 
{\cred \cite[Section 1 and Proposition 1.5(ii)]{SD74} (cf. \cite[Lemma 2.5 and Remark 2.6]{Isk77})}. 
Therefore, we obtain 
\[
W_C = \bigcap_{\substack{C \subset Q,\\ Q:\text{quadric}}} Q 
\subset 
\bigcap_{\substack{(W_C)_{\red} \subset Q,\\ Q:\text{quadric}}} Q =(W_C)_{\red}, 
\]
which implies $W_C= (W_C)_{\red}$, as required. 
Furthermore, $W_C$ is smooth by \cite[Lemma 4.10]{SD73}. 
\end{proof}

\begin{lem}\label{l-int-quad-IOA}
Let $X \subset \P^N$ be a projective variety with $\dim X \geq 2$. 
Fix a hyperplane $\P^{N-1} \subset \P^N$ and set $Y := X \cap \P^{N-1}$. 
Assume that 
 \begin{enumerate}
    \item[(i)] $H^0(\P^N, \MO_{\P^N}(1)) \xrightarrow{\simeq} H^0(X, \MO_{\P^N}(1)|_X)$ and  
     \item[(ii)] $\bigoplus_{m\geq 0} (X, \MO_{\P^N}(m)|_X)$ is generated by $H^0(X, \MO_{\P^N}(1)|_X)$. 
 \end{enumerate}
Then the following hold. 
\begin{enumerate}
\item 
For 
\begin{itemize}
    \item $A :=$ the set of the quadric hypersurfaces of $\P^N$ containing $X$, and 
    \item $B:=$ the set of the quadric hypersurfaces of $\P^{N-1}$ containing $Y$, 
\end{itemize}
the following map is bijective: 
\[
A \to B, \qquad Q \mapsto Q \cap \P^{N-1}. 
\]
\item $X \subset \P^N$ is an intersection of quadrics if and only if 
$Y \subset \P^{N-1}$ {\cred is} an intersection of quadrics. 
\end{enumerate}
\end{lem}

\begin{proof}
See \cite[Lemma 2.10]{Isk77}. 
\end{proof}

\begin{prop}\label{p-W-minimal}
We use Notation \ref{n-trigonal}. 
Then the following hold. 
\begin{enumerate}
    \item $W$ is a $4$-dimensional projective normal variety such that $\Delta(W, \MO_{\P^{g+1}}(1)|_W)=0$ and $\deg W = g-2$. 
    \item $\dim \Sing\,W \leq 1$. 
\end{enumerate}
\end{prop}

\begin{proof}
Let us show (1). 
Let $\P^g \subset \P^{g+1}$ be a hyperplane such that $S := X \cap \P^g$ is smooth, and hence a K3 surface. 
Let $\P^{g-1} \subset \P^g$ be a hyperplane of $\P^g$ such that 
$C:= S \cap \P^{g-1}$ is smooth, and hence a canonical curve. 
{\cred Set $\MO_S(1) := \MO_{\P^g}(1)|_S$ and $\MO_C(1) := \MO_{\P^{g-1}}(1)|_C$.} 
Note that $R(S, \MO_S(1))$ and $R(C, \MO_C(1))$ are generated by 
$H^0(S, \MO_S(1))$ and $H^0(C, \MO_C(1))$, respectively. 
In particular, the following 3 sets are bijectively corresponding via restriction {\cred (Lemma \ref{l-int-quad-IOA}(1))}: 
\begin{itemize}
    \item The set of the quadric hypersurfaces of $\P^{g+1}$ containing $X$.     
    \item The set of the quadric hypersurfaces of $\P^{g}$ containing $S$.
    \item The set of the quadric hypersurfaces of $\P^{g-1}$ containing $C$.
\end{itemize}
Recall that 
\[
W := \bigcap_{\substack{X \subset Q,\\ Q:\text{quadric}}} Q  
\]
is the scheme-theoretic intersection of the quadric hypersurfaces containing $X$. 
Similarly, we set 
\[
W_S := \bigcap_{\substack{S \subset Q,\\ Q:\text{quadric}}} Q,\qquad \text{and}\qquad 
W_C := \bigcap_{\substack{C \subset Q,\\ Q:\text{quadric}}} Q. 
\]
The above bijective correspondence deduces the following equality of closed subschemes of $\P^{g+1}$: 
\[
W_C = W \cap \P^{g-1}. 
\]
By Theorem \ref{t-NEP}, $W_C \subset \P^{g-1}$ is a surface of minimal degree, and hence $W_C$ is an intersection of quadrics. 
Therefore, we obtain a decomposition 
\[
W_{\red} = W_1 \cup W_2,
\]
into reduced closed subschemes, where 
$W_1$ is a $4$-dimensional projective variety and $\dim W_2 \leq 1$ ($W_2$ is possibly reducible). 
We have $X \subset W_1$ and 
\[
W_C = W \cap \P^{g-1} = W_1 \cap \P^{g-1}. 
\]
Then $\deg W_1 = \deg W_C = g-2$, and hence $W_1 \subset \P^{g+1}$ 
is of minimal degree. 
In particular, $W_1$ is an intersection of quadrics. 
Therefore, 
\[
W_1 \subset W_{\red} \subset W = \bigcap_{\substack{X \subset Q,\\ Q:\text{quadric}}} Q \subset \bigcap_{\substack{W_1 \subset Q,\\ Q:\text{quadric}}} Q  =W_1, 
\]
which implies $W =W_1$. Thus (1) holds. 

Let us show (2). 
Suppose $\dim \Sing\,W \geq 2$. 
As $W$ is normal by (1), we obtain $\dim \Sing\,W = 2$. 
We take two general hyperplane sections and its intersection: $\P^{g-1} = \P^g \cap \P^g$. 
Then $C := X \cap \P^{g-1}$ and $W_C = W \cap \P^{g-1}$ are smooth 
(Theorem \ref{t-NEP}), which is a contradiction. 
Hence $\dim {\rm Sing} W \leq 1$. 
Thus (2) holds. 
\end{proof}






\begin{prop}\label{p-trigonal-Wsm}
We use Notation \ref{n-trigonal}. 
If $W$ is smooth, then $\rho(X) \geq 2$. 
\end{prop}

\begin{proof}
We have the induced closed embeddings: $X \subset W \subset \P^{g+1}$. 
Recall that $\deg W = g -2 \geq 5-2 =3$. 
Since $W$ is smooth, we have $(W, \MO_{\P^{g+1}}(1)|_W) = (\P^4, \MO_{\P^4}(1)), (W, \MO_{\P^{g+1}}(1)|_W) = (Q^4, \MO_{\P^5}(1)|_{Q^4})$, or $W$ is a $\P^3$-bundle over $\P^1$ (Theorem \ref{t-Delta0-sm}). 
If $(W, \MO_{\P^{g+1}}(1)|_W) = (\P^4, \MO_{\P^4}(1))$ or $(W, \MO_{\P^{g+1}}(1)|_W) =(Q^4, \MO_{\P^5}(1)|_{Q^4})$, then we have $\deg W \leq 2$, which is a contradiction. 
Hence $W$ is a $\P^3$-bundle over $\P^1$. 
Let $\pi: W \to \P^1$ be the projection. 
If $\pi(X) = \P^1$, then $\rho(X) \geq 2$. 

Therefore, we may assume that $\pi(X)$ is a point. 
In this case, $X$ is contained in a fibre of a $\P^3$-bundle $\pi : W \to \P^1$, and hence $X \simeq \P^3$. 
Then $X \subset \P^{g+1}$ is a Veronese variety, which is an intersection of quadrics (Lemma \ref{l-Veronese-N1}). 
Hence this case is excluded, because 
we assume, in Notation \ref{n-trigonal}, that $X$ is not an intersection of quadrics. 
\end{proof}

\begin{lem}\label{l-Veronese-N1}
For positive integers $n$ and $d$, set $V$ to be the image of the closed immersion 
$\varphi_{|\MO_{\P^n}(d)|} : \P^n \hookrightarrow \P^N$, where 
$N := \binom{n+d}{n}-1$. 
Then $V \subset \P^N$ is an intersection of quadrics. 
\end{lem}

The variety $V$ as above is called a {\em Veronese variety}. 

\begin{proof}
See \cite[Chapter 1, Subsection 4.4, Example 1.28]{Sha13} 
or \cite[Theorem IV.25 and Appendix D]{Her06}. 
\end{proof}


Note that the singular locus $\Sing\,W$ of $W$ is a linear subvariety of $\P^{g+1}$. In particular, the following hold. 
\begin{itemize}
    \item $\dim \Sing\,W=0$ if and only if $\Sing\,W$ is a point. 
    \item $\dim \Sing\,W=1$ if and only if $\Sing\,W$ is a line. 
\end{itemize}

\begin{prop}\label{p-tri-0dimv}
We use Notation \ref{n-trigonal}. 
Assume $\rho(X)=1$. 
Then $\dim \Sing\,W \neq 0$. 
\end{prop}


\begin{proof}
Suppose that $\dim \Sing\,W =0$, i.e., $P:=\Sing\,W$ is a point. 
Let us derive a contradiction. 
Note that $W \subset \P^{g+1}$ is a cone over a smooth threefold $Z \subset \P^g$ of minimal degree. 
By $\deg Z = \deg W = g-2 \geq 3$, $Z$ is a $\P^2$-bundle over $\P^1$ {\cred (Theorem \ref{t-Delta0-sm})}. 
In particular, $\rho (Z) =2$. 

For the blowup $\widetilde{\mu} : W' \to W$ of $W$ at $P$, we get the induced $\P^1$-bundle $\widetilde{\pi}: W' \to Z$. 
For the proper transform $X' := \mu_*^{-1}(X)$ of $X$, 
we have the induced morphisms $\mu : X' \to X$ and $\pi := \widetilde \pi|_{X'} : X' \to Z$. 
\[
\begin{tikzcd}
W' \arrow[d, "\widetilde{\pi}"'] \arrow[dr, "\widetilde{\mu}"]\\
Z & W
\end{tikzcd} \qquad \qquad 
\begin{tikzcd}
X' \arrow[d, "\pi"'] \arrow[dr, "\mu"]\\
Z & X.
\end{tikzcd}
\]
We treat the following two cases (I) and (II) separately: 
\[
{\rm (I)}\,P \not\in X \hspace{20mm} 
{\rm (II)}\,P \in X. 
\]


(I) Assume $P \not\in X$. 
In this case, we get $\mu : X' \xrightarrow{\simeq} X$, and hence $\rho(X') = \rho(X) =1$. 
Since any fibre of $\widetilde{\pi}$ is one-dimensional, 
it follows from $\rho(X')=1$ that $\pi: X'  \to Z$ is surjective. 
Then the inequality $1 =\rho(X) \geq \rho(Z)$ implies $\rho(Z)=1$. 
However, this contradicts $\rho(Z)=2$. 

(II) Assume $P  \in X$. 
In this case, $\mu : X' \to X$ is the blowup of $X$ at $P$. 
In particular, $X'$ is a smooth projective threefold with $\rho(X')=2$. 
Since any fibre of $\widetilde{\pi}$ is one-dimensional, 
it holds that $\dim \pi(X') =2$ or $\dim \pi(X')=3$.

Let us show $\dim \pi(X') \neq 2$. 
Suppose $\dim \pi(X') =2$. 
In order to apply Lemma \ref{l-prim-P3}, 
let us confirm its assumptions, i.e., every fibre of $\pi|_{X'} : X' \to \pi(X')$ is $\P^1$ and $\pi(X') = \P^2$. 
For a closed point $z \in \pi(X')$, we obtain a scheme-theoretic inclusion $X'_z \subset W'_z =\P^1$. 
Hence any fibre $X'_z$ of $\pi|_{X'}  : X' \to \pi(X')$ is isomorphic to $\P^1$. 
It holds that $\pi(X')$ is a fibre of the $\P^2$-bundle $Z \to \P^1$, 
as otherwise 
the composite morphism $X' \to Z \to \P^1$ is surjective, which leads to the following contradiction for the normalisation $\pi (X')^{N}$ of $\pi (X')$: 
\[
2 = \rho(X') > \rho( \pi(X')^N) > \rho(\P^1) =1. 
\]
Therefore, $\pi(X') \simeq \P^2$, and hence we may apply Lemma \ref{l-prim-P3}, 
so that $X \simeq \P^3$. 
This contradicts Lemma \ref{l-Veronese-N1}, which completes the proof of $\dim \pi(X') \neq 2$.



Hence we obtain $\dim \pi(X') = 3$, and hence $\pi: X' \to Z$ is surjective. 
We get a surjection to $\P^1$: $g: X' \xrightarrow{\pi} Z \xrightarrow{\rho} \P^1$, 
where $\rho$ denotes the $\P^2$-bundle. 
Taking the Stein factorisation of $g$, we obtain a morphism 
\[
h: X' \to B
\]
to a smooth projective curve $B$ with $h_*\MO_{X'} = \MO_B$. 
Set $E: =\Ex(\mu) = \P^2$. 
Then the composite morphism 
\[
\P^2 = E \hookrightarrow X' \xrightarrow{h} B
\]
is trivial. 
We then obtain the following factorisation: 
\[
h : X' \xrightarrow{\mu} X  \xrightarrow{q} B. 
\]
However, this is a contradiction, 
because $q$ is surjective and $\rho(X)=1$. 
\qedhere


\end{proof}

The following lemma has {\cred been} already used in the above proof. 
We shall establish a more general result of this lemma in \cite{AT}. 

\begin{lem}\label{l-prim-P3}
Let $X$ be a Fano threefold with $\rho(X)=1$. 
Let $\mu: Y \to X$ be a blowup at a point $P \in X$. 
Assume that there exists a morphism $\pi : Y \to \P^2$ such that 
every fibre of $\pi$ is isomorphic to $\P^1$. 
Then $X \simeq \P^3$. 
\end{lem}

\begin{proof}
We now show that $r_X=2$ or $r_X =4$. 
Suppose $r_X=1$ or $r_X=3$. 
Let us derive a contradiction. 
We can write $-K_X \sim r_X H_X$ for some ample Cartier divisor $H_X$ on $X$ satisfying $\Pic\,X = \Z H_X$. 
In particular, we obtain $\Pic\,Y = \Z(\mu^*H_X) \oplus \Z E$ for $E := \Ex(\mu)$. 
Let $L$ be a line on $\P^2$ and set $D := \pi^*L$. 
In particular, $D$ is a $\P^1$-bundle over $\P^1$. 
By $\MO_Y(E)|_E \simeq \MO_E(-1)$ and $K_Y=\mu^*K_X +2E$, we obtain  
\[
K_Y^2 \cdot E = 4, \qquad 
K_Y \cdot E^2 = 2, \qquad E^3 =1. 
\]
We also have 
\[
\mu^*K_X \cdot K_Y^2 = \mu^*K_X \cdot (\mu^*K_X +2E)^2 = K_X^3
\]
\[
(\mu^*K_X)^2 \cdot K_Y = (\mu^*K_X)^2 \cdot (\mu^*K_X +2E) = K_X^3. 
\]
It holds that 
\[
D \equiv -a \mu^*K_X -bE
\]
for some  $a \in \frac{1}{3}\Z$ and $b \in \Z$. 
Since $\kappa(Y, D) =2$ and $\kappa(Y, \mu^*(-K_X))=3$, 
it holds that $a>0$ and $b>0$. 
Since $D$ is a $\P^1$-bundle over $\P^1$, it holds that 
\[
D \cdot K_Y^2 = D \cdot (K_Y +D -D)^2 = K_D^2 - 2K_D \cdot (\text{a fibre of }D \to L)= 8 - 2 (-2)=12. 
\]
We obtain 
\[
12 = (-a \mu^*K_X -bE) \cdot K_Y^2 =-a K_X^3 -4b. 
\]
By $D^2 \cdot K_Y = \deg(K_Y|_{\text{a fibre of }\pi})= -2$ and $\mu^*K_X \cdot E =0$, we obtain 
\[
-2 = D^2 \cdot K_Y = (-a \mu^*K_X -bE)^2 \cdot K_Y 
= a^2 K_X^3 +2b^2. 
\]
Hence we get $a(4b+12) =-a^2 K_X^3 = 2b^2 +2$, which implies 
$2a(b+3) = b^2 +1$. 
We can write $a = a'/3$ for some $a' \in \Z$, so that  
\[
2a'(b+3) = 3(b^2 +1). 
\]
Then $b$ is odd, and hence we have $b =2b'+1$ for some $b' \in \Z_{\geq 0}$, 
which implies 
\[
2a'(b'+2) = 3(2b'^2 +2b' +1). 
\]
This is a contradiction, because the left (resp. right) hand side is even (resp. odd). 
Therefore, $r_X =2$ or $r_X =4$. 

{\cred Note that $-K_Y$ is ample by Kleiman's criterion.} 
Then it follows from $K_Y = \mu^*K_X +2E$ that $r_Y =2$. 
By $\rho(Y)=2$, there are only two possibilities (Theorem \ref{t-index-2}): 
$Y \in |\MO_{\P^2 \times \P^2}(1, 1)|$ or $Y$ is a blowup of $\P^3$ at a point. 
For the former case, the two projections $\pi_1, \pi_2: Y \to \P^2$ correspond to the distinct extremal rays. 
Therefore, we have that $Y$ is a blowup of $\P^3$ at a point. 
In this case, we obtain $X \simeq \P^3$. 
\qedhere
\end{proof}




\begin{prop}\label{p-tri-1dimv}
We use Notation \ref{n-trigonal}. 
Assume $\rho(X)=1$. 
Then $\dim \Sing\,W \neq  1$. 
\end{prop}

\begin{proof}
Suppose $\dim \Sing\,W =1$, i.e., $\ell :=\Sing\,W$ is a line on $\P^{g+1}$. 
Let us derive a contradiction. 
In this case, $W \subset \P^{g+1}$ is a cone over a smooth surface $Z \subset \P^{g-1}$ of minimal degree. 
{\cred By $\deg Z = \deg W = g-2 \geq 3$,} either $Z$ is a Hirzebruch surface or 
 $Z \subset \P^{g-1}$ {\cred is} a Veronese surface  with $g=6$ 
 {\cred (Theorem \ref{t-Delta0-sm})}. 
By Theorem \ref{t-Veronese-main}, a cone $W$ over the Veronese surface $Z$ does not contain a smooth prime divisor. 
Hence $Z$ is a Hirzebruch surface. 
In particular, $\rho(Z)=2$.

For the blowup $\widetilde{\mu} : W' \to W$ of $W$ along $\ell$, 
we get the induced $\P^2$-bundle $\widetilde{\pi}: W' \to Z$. 
For the proper transform $X' := \mu_*^{-1}(X)$ of $X$, 
we have the induced morphisms $\mu : X' \to X$ and $\pi := \widetilde \pi|_{X'} : X' \to Z$. 
\[
\begin{tikzcd}
W' \arrow[d, "\widetilde{\pi}"'] \arrow[dr, "\widetilde{\mu}"]\\
Z & W
\end{tikzcd} \qquad \qquad 
\begin{tikzcd}
X' \arrow[d, "\pi"'] \arrow[dr, "\mu"]\\
Z & X.
\end{tikzcd}
\]
In what follows, we treat the following two cases separately: 
\[
{\rm (I)}\,\,\ell \not\subset X\hspace{20mm} {\rm (II)}\,\,\ell \subset X. 
\]

(I) 
Assume $\ell \not\subset X$. 
In this case, $X \cap \ell$ is either empty or zero-dimensional. 
Let $S \subset X \times_k \kappa$ be the generic hyperplane section of $X$. 
Note that we have $\rho(S) =\rho(X)=1$ \cite[Proposition 5.17]{Tana}. 
Since blowups commute with flat base changes, 
$\widetilde{\mu} \times_k \kappa : W' \times_k \kappa \to W \times_k \kappa$ is the blowup along  $\ell \times_k \kappa$. 
Then $S$ is disjoint from the blowup centre $\ell \times_k \kappa$. 
By $\rho(S)=1$, either $S$ is contained in a fibre of 
$\widetilde \pi \times_k \kappa : W'\times_k \kappa \to Z \times_k \kappa$ or 
the induced morphism $\pi_S : S \to Z \times_k \kappa$ is surjective. 
The former case is impossible, because $\kappa(S) =0 \neq -\infty=\kappa(\P^2)$ (note that any fibre of $\widetilde \pi \times_k \kappa$ is $\P^2$). 
Then $\pi_S : S \to Z \times_k \kappa$ is surjective. 
Since $Z$ is a Hirzebruch surface, we obtain the following contradiction: 
\[
1 = \rho(S) \geq \rho(Z \times_k \kappa) \geq \rho(Z)=2. 
\]

(II) Assume $\ell \subset X$. Set $E := \Ex(\mu)$, which is a $\P^1$-bundle over $\ell = \P^1$. 

Let us show that $E \not\simeq \P^1 \times \P^1$. 
Recall that $\ell$ is a line on $\P^{g+1}$. 
It follows from \cite[Proposition 2.2.14]{IP99}  
that $\deg N_{\ell/X} = 2g(\ell) -2 - K_X \cdot \ell = -1$. 
We have $E \simeq \P_{\ell}(N_{\ell/X})$. 
Since $\deg N_{\ell/X}$ is an odd number, we have that $E \not\simeq \P^1 \times \P^1$.


We now show that  $\pi|_E : E \to Z$ is surjective. 
Otherwise, its image is either a point or a curve. 
For the former case, we get the   factorisation: $\pi : X' \xrightarrow{\mu} X \to Z$, 
which contradicts $\rho(X)=1$. 
Suppose the latter case, i.e., the image of $\pi|_E : E \to Z$ is a curve. 
Since $E$ is a Hirzebruch surface with $E \not\simeq \P^1 \times \P^1$, 
the induced morphism $\mu|_E : E \to \ell$ 
is the unique morphism to a curve satisfying $(\mu|_E)_*\MO_E = \MO_{\ell}$. 
Hence we again get the   factorisation: $\pi : X' \xrightarrow{\mu} X \to Z$, 
which contradicts $\rho(X)=1$. 
This completes the proof of the surjectivity of $\pi|_E : E \to Z$. 

Since $\pi|_E : E \to Z$ is surjective, also $\pi : X' \to Z$ is surjective. 
This leads to the following contradiction: $2=\rho(X') > \rho(Z)=2$. 
\end{proof}

We are ready to prove the main result of this section.

\begin{thm}\label{t-int-of-quad}
Let $X \subset \P^{g+1}$ be an anti-canonically embedded Fano threefold 
with $g \geq 5$ and $\rho(X)=1$. 
Then $X$ is an intersection of quadrics. 
\end{thm}

\begin{proof}
Suppose that $X$ is not an intersection of quadrics. 
We use Notation \ref{n-trigonal}. 
By Proposition \ref{p-W-minimal}, 
$W$ is a $4$-dimensional normal variety with $\dim \Sing\,W \leq 1$. 
Hence the assertion follows from 
Proposition \ref{p-trigonal-Wsm}, 
Proposition \ref{p-tri-0dimv}, and Proposition \ref{p-tri-1dimv}. 
\end{proof}

\appendix

\section{Singular varieties of minimal degree}


\subsection{Rational normal scrolls}

We summarise some basic  properties on rational normal scrolls. 
In order to treat rational normal scrolls and 
cones over the Veronese surface simultaneously, 
we consider a slightly generalised setting, i.e., 
we consider a projective space bundle over $\P^n$ instead of over $\P^1$. 
Most of the arguments in this subsection are based on \cite{EH87}.

\begin{nota}\label{n-rat-scroll}
Fix $n \in \Z_{>0}$. 
Let $0 \leq a_0 \leq a_1 \leq \cdots \leq a_d$ be integers with $a_d >0$. 
For 
\[
E := \MO_{\P^n}(a_0) \oplus \cdots \oplus \MO_{\P^n}(a_d), 
\]
we set 
\[
\P_{\P^n}(a_0, ..., a_d) := \P(E) := \P_{\P^n}(E) =  \P_{\P^n}( \MO_{\P^n}(a_0) \oplus \cdots \oplus \MO_{\P^n}(a_d)),
\]
which is a $\P^d$-bundle over $\P^n$.  
Let $S_{\P^n}(a_0, ..., a_d)$ be its image by $\MO_{\P(E)}(1)$. 
Let 
\[
\varphi_{|\MO_{\P(E)}(1)|} : \P_{\P^n}(a_0, ..., a_d) \to S_{\P^n}(a_0, ..., a_d) \subset \P^{-1+ \sum_i h^0(\P^n, \MO_{\P^n}(a_i))} 
\]
be the morphism induced by the complete linear system $|\MO_{\P(E)}(1)|$, 
which is base point free (cf. Theorem \ref{t-RNS-bpf}(1)). 
Fix a hyperplane $H$ on $\P^n$ and set $F := \pi^*H$, 
where $\pi : \P_{\P^n}(E) \to \P^n$ denotes the projection. 
For each $0 \leq i \leq d$, we have 
\begin{itemize}
    \item the section of $\pi$ 
    \[
    \Gamma_i :=\P_{\P^n}(\MO_{\P^n}(a_i)) \subset \P_{\P^n}(E)
    \]
    corresponding to the projection $E \to \MO_{\P^n}(a_i)$, and 
    \item the $\P^{d-1}$-bundle over $\P^n$ 
    \[
    D_i := \P_{\P^n}(\MO_{\P^n}(a_0) \oplus \cdots \oplus \MO_{\P^n}(a_{i-1}) \oplus 
    \MO_{\P^n}(a_{i+1}) \oplus \cdots \oplus \MO_{\P^n}(a_{d})) \subset \P_{\P^n}(E) 
    \]
    corresponding to the projection 
    \[
    E \to \MO_{\P^n}(a_0) \oplus \cdots \oplus \MO_{\P^n}(a_{i-1}) \oplus 
    \MO_{\P^n}(a_{i+1}) \oplus \cdots \oplus \MO_{\P^n}(a_{d}).
    \]
\end{itemize} 
\end{nota}

\begin{rem}\label{r-rat-scroll}
We use Notation \ref{n-rat-scroll}. 
By construction, the following hold {\cred for each $0 \leq i \leq d$}. 
\begin{enumerate}
    \item We have $D_i \cap \Gamma_i = \emptyset$. 
    \item $\MO_{\P(E)}(1)|_{\Gamma_i} = \MO_{\P(E)}(1)|_{\P(\MO_{\P^n}(a_i))}=
\MO_{\P(\MO_{\P^n}(a_i))}(1) = \MO_{\P^n}(a_i)$.  
\item 
By (2), 
$\varphi_{|\MO_{\P(E)}(1)|}(\Gamma_i)$ is a point if and only if $a_i=0$. 
\end{enumerate}
\end{rem}


\begin{lem}\label{l-RNS-lin-eq}
We use Notation \ref{n-rat-scroll}. For each $0 \leq i \leq d$, it holds that 
\[
\MO_{\P(E)}(1) \sim D_i + a_iF. 
\]
\end{lem}

\begin{proof}
We can write $D_i \sim \MO_{\P(E)}(1) + r F$ for some $r \in \Z$. 
{\cred Via $\pi|_{\Gamma_i} : \Gamma_i \xrightarrow{\simeq} \P^n$, it holds that 
\[
\MO_{\P^n} \simeq \MO_{\P(E)}(D_i)|_{\Gamma_i} \simeq 
(\MO_{\P(E)}(1) \otimes  \MO_{\P(E)}(rF))|_{\Gamma_i}
\simeq \MO_{\P^n}(a_i + r).  
\]
Hence $a_i +r =0$.} 
We are done. 
\end{proof}

\begin{thm}\label{t-RNS-bpf}
We use Notation \ref{n-rat-scroll}. 
\begin{enumerate}
    \item $|\MO_{\P(E)}(1)|$ is base point free. 
    \item $|\MO_{\P(E)}(1)|$ is very ample if and only if $0<a_0$, i.e., all of $a_0, ..., a_d$ are positive. 
    \item $K_{\P(E)} + \sum_{i=0}^d D_i \sim \pi^*K_{\P^n}$. 
    \item $S(0, ..., 0, a_0, ..., a_d) \subset \P^{s{\cred -1} +\sum_i h^0(\P^n, \MO_{\P^n}(a_i))}$ is the $s$-th cone over $S(a_0, ..., a_d) \subset \P^{{\cred -1+} \sum_i h^0(\P^n, \MO_{\P^n}(a_i))}$, where $s$ is the number of $0$ appearing in $S(0, ..., 0, a_0, ..., a_d)$. 
\end{enumerate}
\end{thm}

\begin{proof}
The assertion (1) follows from Lemma \ref{l-RNS-lin-eq}, 
$D_0 \cap \cdots \cap D_d =\emptyset$, 
and {\cred the fact that} $|F|$ is base point free. 

Let us show (2). 
Assume that $\MO_{\P(E)}(1)$ is very ample. 
Then all of $a_0, ..., a_d$ are positive by Remark \ref{r-rat-scroll}(3). 
Let us prove the opposite implication. 
Assume that all of $a_0, ..., a_d$ are positive. 
As $\P_{\P^n}(E)$ is a smooth projective toric variety, 
it is enough to show that $\MO_{\P(E)}(1)$ is ample 
{\cred \cite[Theorem 6.1.15]{CLS11}}. 
Fix a curve $C$ on $\P_{\P^n}(E)$. 
Since $|\MO_{\P(E)}(1)|$ is base point free by (1), 
it suffices to show that $\MO_{\P(E)}(1) \cdot C >0$. 
If $\pi(C)$ is a point, then this follows from the fact that 
$\MO_{\P(E)}(1)$ is $\pi$-ample. 
Hence we may assume that $\pi(C)$ is a curve. 
{\cred By $D_0 \cap \cdots \cap D_d =\emptyset$,} 
there exists $0 \leq i \leq d$ such that $C \not\subset D_i$. 
Then 
\[
\MO_{\P(E)}(1) \cdot C = D_i \cdot C  + a_i F\cdot C \geq a_i F \cdot C >0. 
\]
Thus (2) holds.

Let us show (3). 
Fix a fibre $\P^d$ of $\pi : \P(E) \to \P^n$. 
Then we obtain
\[
(K_{\P(E)}+\sum_{i=0}^d D_i)|_{\P^d} = K_{\P^d} + \sum_{i=0}^d (D_i|_{\P^d}) \sim 0, 
\]
because $D_0|_{\P^d}, ..., D_d|_{\P^d}$ are hyperplanes. 
Hence we get $K_{\P(E)}+\sum_{i=0}^d D_i \sim \pi^*L$ for some Cartier divisor {\cred $L$} on $\P^n$. 
By taking the restriction to the section $\Gamma_0 =\P^n$ of $\pi$, we get 
\[
L \sim (K_{\P(E)}+\sum_{i=0}^d D_i)|_{\Gamma_0} = 
(K_{\P(E)}+\sum_{i=1}^d D_i)|_{\Gamma_0} \sim K_{\Gamma_0}, 
\]
where the latter linear equivalence follows from adjunction formula 
{\cred and $D_1 \cap \cdots \cap D_n =\Gamma_0$}. 
Thus (3) holds. 

Concerning (4), we may apply the same argument as in \cite[page 6]{EH87}. 
\end{proof}

\begin{prop}\label{p-RNS-Qfac}
We use Notation \ref{n-rat-scroll}. 
Then the following are equivalent. 
\begin{enumerate}
    \item $S_{\P^n}(a_0, ..., a_d)$ is $\Q$-factorial. 
    \item 
    One of the following holds. 
    \begin{enumerate}
        \item $a_0>0$, i.e., all of $a_0, ..., a_d$ are positive. 
        \item $a_0 = \cdots = a_{d-1}=0$, i.e., only $a_d$ is positive. 
    \end{enumerate}
\end{enumerate}
Furthermore, if {\rm (b)} holds, then 
$\varphi_{|\MO_{\P(E)}(1)|}$ is a birational morphism with $\Ex(\varphi_{|\MO_{\P(E)}(1)|}) = D_d$.  
\end{prop}

\begin{proof}
Let us show (1) $\Rightarrow$ (2). 
Assume that (2) does not hold. 
Then $a_0=0$ and $a_{d-1}>0$. 
In particular, $a_d >0$. 
Let $C$ be a curve on $\P_{\P^n}(E)$ such that $\varphi_{|\MO_{\P(E)}(1)|}(C)$ is a point. 
In order to prove that $S(a_0, ..., a_d)$ is not $\Q$-factorial, it is enough to show that $C \subset D_{d-1} \cap D_d$. 
It follows that $\pi(C)$ is a curve, 
because $\pi$ and $\varphi_{|\MO_{\P(E)}(1)|}$ correspond to distinct extremal rays. 
For $i \in \{ d-1, d\}$, we have 
\[
0 = \MO_{\P(E)}(1) \cdot C = ( D_i + a_i F) \cdot C = D_i \cdot C + a_i F \cdot C >D_i \cdot C, 
\]
which implies $C \subset D_i$. This completes the proof of  (1) $\Rightarrow$ (2). 

Let us show (2) $\Rightarrow$ (1). 
Assume (2). 
It suffices to prove that $S(a_0, ..., a_d)$ is $\Q$-factorial. 
If (a) holds, then $\varphi_{|\MO_{\P(E)}(1)|}$ is a closed immersion 
(Theorem \ref{t-RNS-bpf}), 
so that we get $\P(a_0, ..., a_d) \simeq S(a_0, ..., a_d)$, and hence $S(a_0, ..., a_d)$ is $\Q$-factorial. 
We may assume that $a_0 = \cdots = a_{d-1}=0$. 
We then have $E=\MO_{\P^n}^{\oplus d} \oplus \MO_{\P^n}(a_d)$ with $a_d >0$. 
By 
\[
D_d = \P_{\P^n}(\MO_{\P^n}^{\oplus d}) \simeq \P^n \times \P^{d-1}, 
\]
we see that $\dim D_d > \dim \varphi_{|\mathcal{O}_{\mathbb{P}(E)}(1)|}(D_d)$. 
Hence $\varphi_{|\MO_{\P(E)}(1)|} : \P_{\P^n}(a_0, ..., a_d) \to S_{\P^n}(a_0, ..., a_d)$ is a divisorial contraction, i.e., $\Ex(\varphi_{|\MO_{\P(E)}(1)|})$ is a prime divisor. 
Since $\P(a_0, ..., a_d)$ is toric, $S(a_0, ..., a_d)$ is $\Q$-factorial {\cred \cite[Proposition 15.4.5]{CLS11}}. 
\end{proof}

{\cred 
We shall need the following lemma in the next subsection. 
Although this result is well known, we give a proof for the sake of completeness. 

\begin{lem}\label{l-cone-blouwp}
Let $X \subset \P^n =\Proj\,k[x_0, ..., x_n]$ be a projective variety. 
Take $s \in \Z_{>0}$ and let $Y \subset \P^{n+s} = \Proj\,k[x_0, ..., x_n, y_1, ..., y_s]$ be the $s$-th cone of $X$.  
Set $L$ to be the vertex, i.e., $L$ is the reduced closed subscheme of $\P^{n+s}$ whose closed points are 
$\{ [x_0 : \cdots : x_n: y_1: ...: y_s] \in \P^{n+s}(k) \,|\, x_0 = ...=x_n=0, y_1, \cdots, y_s \in k\}$. 
Then, for $\MO_X(1) := \MO_{\P^n}(1)|_X$ and $Z := \P_X( \MO^{\oplus s}_X \oplus \MO_X(1))$,  there exist the following morphisms 
\[
\begin{tikzcd}
& Z= \P_X( \MO_X^{\oplus s} \oplus \MO_X(1)) \arrow[ld, "\pi"'] \arrow[rd, "\sigma"] &\\
X && Y,
\end{tikzcd}
\]
where $\pi : Z= \P_X( \MO_X^{\oplus s} \oplus \MO_X(1)) \to X$ denotes the induced $\P^s$-bundle and 
$\sigma : Z \to Y$ is the blowup along the vertex $L$. 
\end{lem}

\begin{proof}
The proof consists of the following two steps: 
\begin{enumerate}
\item[(I)] The case when $X = \P^n$. 
\item[(II)] The general case. 
\end{enumerate}

(I) We first treat the case when $X= \P^n$. 
In this case, we have $Y = \P^{n+s}$ and $Z = \P_{\P^n}(\MO_{\P^n}^{\oplus s} \oplus \MO_{\P^n}(1))$. 
By Theorem \ref{t-RNS-bpf}(4), we get the above diagram, 
where $\pi$ is the induced $\P^s$-bundle and $\sigma$ is a birational morphism. 
It is enough to show that $\sigma: Z \to Y$ is the blowup along the vertex $L$. 
Note that $L$ is scheme-theoretically equal to the base scheme of the linear system 
that induces the dominant rational map 
\[
Y =\P^{n+s} \dashrightarrow 
\P^n =X,\qquad 
[x_0: \cdots :x_n:y_1: \cdots :y_s] \mapsto [x_0: \cdots :x_n].
\]
Then the blowup $\sigma' : Z' \to Y$ along $L$ coincides with the resolution of its indeterminacies. 
By the universal property of blowups \cite[Ch. II, Proposition 7.14]{Har77}, 
we obtain a factorisation $\sigma : Z \xrightarrow{\theta} Z' \xrightarrow{\sigma'} Y$. 
Then $\theta : Z \to Z'$ is an isomorphism, because $\theta : Z \to Z'$ 
is a birational morphism of normal projective varieties satisfying 
$\rho(Z) = \rho(Z')$. 

(II) Let us treat the general case. 
Set 
$\widetilde X := \P^n, \widetilde Y := \P^{n+s},$ and $\widetilde Z := \P_{\P^n}(\MO_{\P^n}^{\oplus s} \oplus \MO_{\P^n}(1))$. 
By the case (I), we have induced morphisms
\[
\begin{tikzcd}
& \widetilde{Z}= \P_{\P^n}( \MO_{\P^n}^{\oplus s} \oplus \MO_{\P^n}(1)) \arrow[ld, "\widetilde \pi"'] \arrow[rd, "\widetilde \sigma"] &\\
\widetilde{X}=\P^n && \widetilde{Y}=\P^{n+s}, 
\end{tikzcd}
\]
where $\widetilde \pi$ is the induced $\P^s$-bundle and $\widetilde \sigma$ is the blowup along $L$. 
Note that we have the induced closed embeddings $X \subset \widetilde{X}$ and $Y \subset \widetilde Y$. 
We also have a natural closed embedding $Z \subset \widetilde Z$, 
because there is the following cartesian diagram 
\[
\begin{tikzcd}
Z \arrow[r, hook]\arrow[d, "\pi"] & \widetilde Z\arrow[d, "\widetilde{\pi}"]\\
X \arrow[r, hook] & \widetilde X, 
\end{tikzcd}
\]
where the vertical arrows are the induced $\P^s$-bundles and 
the horizontal ones are closed immersions. 
In particular, $Z = \widetilde{\pi}^{-1}(X)$. 
For the blowup $Z'$ of $Y$ along the vertex $L$ (note that $L \subset Y$), 
we have a closed embeddding $Z' \subset \widetilde Z$ \cite[Ch. II, Corollary 7.15]{Har77}.
By construction, we have a dominant rational map $Y \dashrightarrow X$ compatible with $\widetilde Y \dashrightarrow \widetilde X$, 
i.e., the induced morphisms 
\[
Y \setminus L \to X 
\quad \text{and} 
\quad 
\widetilde Y \setminus L \to \widetilde X
\]
\[
\text{are given by}\quad 
[x_0: \cdots : x_n : y_1: \cdots : y_s] \mapsto [x_0: \cdots :x_n]. 
\]
Hence we get the  morphism $\pi' : Z' \to X$ induced by $Z' \hookrightarrow \widetilde Z \to \widetilde X$. 
For every closed point $x \in X$, 
it holds that 
$\pi'^{-1}(x) = \widetilde \pi^{-1}(x)$, because 
we have 
$\dim \pi'^{-1}(x) \geq \dim Z' -\dim X = \dim Y -\dim X = s$ and 
a scheme-theoretic inclusion 
$\pi'^{-1}(x) \subset \widetilde \pi^{-1}(x) \simeq \P^s$. 
Therefore, we get the following set-theoretic equation: 
\[
Z(k) = \bigcup_{x \in X(k)} \widetilde{\pi}^{-1}(x) = Z'(k), 
\]
where $Z(k)$ (resp. $Z'(k)$) denotes the subset of $Z$ (resp. $Z'$) 
consisting of all the closed points. 
Since both $Z$ and $Z'$ are reduced closed subschemes of $\widetilde Z$, 
we get a scheme-theoretic equation $Z=Z'$. 
\end{proof}
}




\subsection{Cones over the Veronese surface}

The purpose of this subsection is to prove that the $d$-th cone over the Veronese surface does not have a smooth prime divisor when $d \geq 2$ 
({\cred Theorem \ref{t-Veronese-main}}). 
Throughout this subsection, we shall use the following notation.

\begin{nota}\label{n-Vcone}
We work over an algebraically closed field $k$. 
Let $T \subset \P^5$ be the Veronese surface, i.e., $T$ is the image of 
the Veronese embedding: 
\[
\nu : \P^2 \hookrightarrow \P^5, \qquad [x:y:z] \mapsto [x^2:y^2:z^2:xy:yz:zx]. 
\]
Fix an integer $d \geq 1$. 
Let $V \subset \P^{d+5}$ be the $d$-th cone over $T$. 
We then obtain the following diagram
\[
\begin{tikzcd}
W:=\P_{\P^2}(\MO_{\P^2}^{\oplus {\cred d}} \oplus \MO_{\P^2}(2)) \arrow[d, "\pi"'] \arrow[dr, "\sigma"]\\
T \simeq \P^2 & V,
\end{tikzcd}
\]
where $\pi$ denotes the projection and ${\cred \sigma}$ is the blowup along the vertex 
{\cred (Lemma \ref{l-cone-blouwp})}. 
We set  
\[
a_0 :=0,\,\,a_1:=0,\,\,...,\,\,a_{d-1} = 0,\,\,a_d :=2, \quad \text{and}\quad 
E := \MO_{\P^2}^{\oplus d} \oplus \MO_{\P^2}(2) = \bigoplus_{i=0}^d \MO_{\P^2}(a_i). 
\]
We may use Notation \ref{n-rat-scroll} for $n:=2$. 
In particular, 
\begin{enumerate}
    \item $H$ is a line on $T=\P^2$ and $F := \pi^*H$. 
    \item $\sigma : W \to V$ is a birational morphism with $\Ex({\cred \sigma}) = D_d$ {\cred (Proposition \ref{p-RNS-Qfac})}. 
    \item For each $0 \leq i \leq d$, we have a section $\Gamma_i$ of $\pi$ 
    and a $\P^{d-1}$-bundle $D_i \subset W$ over $T=\P^2$, which are mutually disjoint. 
\end{enumerate}
For every $\ell \in \Z$, let $\MO_W(\ell) := \MO_{\P(E)}(\ell)$. 
We set $F_V := \sigma_*F$. 
Note that $D_d \simeq \P^2 \times \P^{d-1}$ and the induced morphism 
$\pi|_{D_d} : D_d \to T =\P^2$ coincides with the first projection. 
Fix a closed point $Q \in \P^{d-1}$ and 
a line $\zeta$ on $\P^2 \times \{ Q\}$. 
Hence 
\[
\zeta \subset \P^2 \times \{ Q\} \subset \P^2 \times \P^{d-1} = D_d.
\]
\end{nota}


We now summarise some basic properties in Proposition \ref{p-Vcone-Cl} and 
Proposition \ref{p-Vcone-discrep}. 

\begin{prop}\label{p-Vcone-Cl}
We use Notation \ref{n-Vcone}. 
Then the following hold. 
\begin{enumerate}
    \item $\Cl\,V = \Z F_V$. 
    \item $F \cdot \zeta =1$ and $D_d \cdot \zeta =-2$. 
    \item $\sigma^*F_V = F + \frac{1}{2} D_d$. Furthermore, $F_V$ is not Cartier. 
    \item For every Weil divisor $B$ on $V$, $2B$ is Cartier. 
\end{enumerate}
In particular, it holds that 
\[
{\rm CaCl}\,V = \Z(2F_V)
\]
for the subgroup ${\rm CaCl}\,V$ of $\Cl\,V$ consisting of the linear equivalence classes of the Cartier divisors. 
\end{prop}

\begin{proof}
The assertion (1) follows from $\Cl\,W = \Z D_d + \Z F$ and 
\[
\Cl\,V = {\cred \sigma}_*(\Cl\,W) = \Z {\cred \sigma}_*D_d + \Z {\cred \sigma}_*F  = \Z F_V. 
\]

Let us show (2). 
We have 
\[
F \cdot \zeta = F|_{\P^2 \times \{ Q\}} \cdot \zeta = (\text{line}) \cdot (\text{line}) = 1. 
\]
Since $\sigma(\zeta)$ is a point, it follows from Lemma \ref{l-RNS-lin-eq} that 
\[
0 = \MO_W(1) \cdot \zeta = (D_d + 2F) \cdot \zeta. 
\]
In particular, $D_d \cdot \zeta = -2$. 
Thus (2) holds. 

Let us show (3). 
Since we can uniquely write $F+ c D_d = \sigma^*F_V$ for some $c \in \Q_{\geq 0}$, 
the equalities $(D_d + 2F) \cdot \zeta=0$ and $\sigma^*F_V \cdot \zeta =0$, 
together with ${\cred D_d} \cdot \zeta \neq 0$, imply $\sigma^*F_V  = F + \frac{1}{2} D_d$. 
In particular, $F_V$ is not Cartier. Thus (3) holds. 

Let us show (4). 
Since $\sigma : W \to V$ is given by the complete linear system 
$|\MO_W(1)|$, 
there exists an ample Cartier divisor $H_V$ on $V$ such that 
$\MO_W(1) \simeq \sigma^*\MO_V(H_V) \simeq \MO_V(\sigma^* H_V)$. By 
\[
\sigma^*H_V \sim D_d +2F = \sigma^*(2F_V), 
\]
we obtain $H_V =\sigma_*\sigma^*H_V \sim \sigma_*\sigma^*(2F_V) = 2F_V$. 
Hence $2F_V$ is Cartier. 
Thus (4) holds by (1). 
\end{proof}

\begin{prop}\label{p-Vcone-discrep}
We use Notation \ref{n-Vcone}. 
Then the following hold. 
\begin{enumerate}
    \item $K_W + \sum_{i=0}^d D_i\sim \pi^*K_{\P^2}$. 
    \item $-K_W \sim (d+1)D_d + (2d+3)F$. 
    \item $-K_V \sim (2d+3)F_V$. 
    \item $K_W = \sigma^*K_V + \frac{1}{2} D_d$. 
    \item $V$ is $\Q$-facotiral and terminal. 
\end{enumerate}
\end{prop}

\begin{proof}
The assertion (1) follows from Theorem \ref{t-RNS-bpf}. 
By $\MO_W(1) \sim D_0 \sim \cdots \sim D_{d-1} \sim D_d + 2F$ (Lemma \ref{l-RNS-lin-eq}), 
we obtain 
\[
K_W + (d+1)D_d + 2dF \sim 
K_W + \sum_{i=0}^d D_i \sim \pi^*K_{\P^2} \sim -3F, 
\]
which implies 
\[
-K_W \sim (d+1)D_d + (2d+3)F \quad \text{and}\quad
-K_V \sim \sigma_*( (d+1)D_d + (2d+3)F) = (2d+3)F_V. 
\]
Thus (2) and (3) hold. 

Let us show (4) and (5). 
There exists $a \in \Q$ such that 
\[
K_W = \sigma^*K_V + a D_d. 
\]
Recall that we have 
$F \cdot \zeta = 1$ and $D_d \cdot \zeta =-2$ (Proposition \ref{p-Vcone-Cl}). 
These, together with (2), imply   
\[
-K_W \cdot \zeta = ( (d+1)D_d + (2d+3)F) \cdot \zeta 
= -2(d+1) + (2d+3) =1. 
\]
By $\sigma^*K_V \cdot \zeta =0$, we get $a = 1/2$, and hence (4) holds. 
In particular, $V$ is terminal. 
Note that $V$ is $\Q$-factorial by Proposition \ref{p-RNS-Qfac}. 
Thus (5) holds. 
\end{proof}

\begin{thm}\label{t-Veronese-main}
Let $V \subset \P^{d+5}$ be the $d$-th cone over the Veronese surface $T \subset \P^5$. 
If $d \geq 2$, then there exists no smooth prime divisor on $V$. 
\end{thm}

\begin{proof}
{\cred In what follows, we use Notation \ref{n-Vcone}.} 
Suppose that there is a smooth prime divisor $X$ on $V$. 
Let us derive a contradiction. 
We have $X \sim b F_V$ for some $b \in \Z_{>0}$ (Proposition \ref{p-Vcone-Cl}).  
Recall that $\dim \Sing\,V = d-1 \geq 1$.

\begin{claim}\label{cl-Veronese-main}
The following hold. 
\begin{enumerate}
 \item 
$b$ is odd. 
\item 
$\Sing\,V \subset X$. 
\end{enumerate}
\end{claim}

\begin{proof}[Proof of Claim \ref{cl-Veronese-main}]
Let us show (1). 
Suppose that $b$ is even. 
By $X \sim bF_V$, $X$ is an effective Cartier divisor  (Proposition \ref{p-Vcone-Cl}). 
Since $X$ is smooth, $V$ is smooth around $X$. 
As $2X$ is an ample effective Cartier divisor, 
$\Sing\,V$ must be (at most) zero-dimensional, 
which contradicts $\dim \Sing\,V = d-1 \geq 1$. 
Hence (1) holds. 

Let us show (2). 
Suppose $\Sing\,V \not\subset X$. 
Fix $v \in \Sing\,V$ with $v \not\in X$. 
By construction, we have $\Sing\,V \subset F_V$. 
For an open neighbourhood $U_v$ of $v \in V$ 
satisfying $U_v \cap X =\emptyset$ and $2F_V|_{U_v} \sim 0$ (cf. Proposition \ref{p-Vcone-Cl}), 
it follows from $X \sim bF_V$ that 
\[
0 =X|_{U_v} \sim bF_V|_{U_v} \sim F_V|_{U_v}, 
\]
where the latter linear equivalence holds by (1) and $2F_V|_{U_v} \sim 0$. 
Taking the pullback by 
the restriction $\sigma' : \sigma^{-1}(U_v) \to U_v$ of $\sigma$, we obtain 
\[
0 \sim  \sigma'^*(F_V|_{U_v}) = \sigma^*(F_V)|_{\sigma^{-1}(U_v)} 
= \left( F + \frac{1}{2} D_d\right)\Big{|}_{\sigma^{-1}(U_v)}. 
\]
This contradicts $D_d|_{\sigma^{-1}(U_v)} \neq \emptyset$. 
Thus (2) holds. 
This complete the proof of  Claim \ref{cl-Veronese-main}. 
\end{proof}

Recall that $\sigma : W \to V$ is the blowup along the vertex $\Sing\,V$ {\cred (Lemma \ref{l-cone-blouwp})}. 
Therefore, for the proper transform $Y := \sigma_*^{-1}X$, 
the induced birational morphism 
$\tau : Y \to X$ is the blowup along  $\Sing\,V$ \cite[Ch. II, Corollary 7.15]{Har77}. 
Hence also $Y$ is a smooth projective variety. 
By Claim \ref{cl-Veronese-main}(2), we obtain 
\[
Y + \alpha D_d = \sigma^*X
\]
for some $\alpha \in \Q_{>0}$. 
We have $\alpha \geq 1/2$ by Proposition \ref{p-Vcone-Cl}. 
By ${\rm codim}_V\,\Sing\,V=3$, 
we obtain $(K_V+X)|_X = K_X$, i.e., $\Diff_X(0)=0$ \cite[Proposition 4.5(1)]{Kol13}. 
This, together with $K_W = \sigma^*K_V + \frac{1}{2}D_d$ (Proposition \ref{p-Vcone-discrep}), implies 
\[
K_Y =\sigma^*K_X + \left( \frac{1}{2} -\alpha\right) (D_d|_Y). 
\]
By $\frac{1}{2} -\alpha \leq 0$, this is a contradiction, because $\sigma : Y \to X$ is the blowup along $\Sing\,V$, so that the coefficient of $D_s|_Y$ must be positive (note that the coefficient of $K_Y -\sigma^*K_X$ is uniquely determined by the negativity lemma \cite[Lemma 3.39]{KM98}). 
\end{proof}

\bibliographystyle{skalpha}
\bibliography{reference.bib}

\end{document}